\documentclass[11pt,a4paper,usenames,dvipsnames]{article}

\usepackage{url,amsmath,amssymb,latexsym,pstricks,mathrsfs,comment,amsthm,graphicx,tikz,tikz-cd,enumerate,accents,pgffor,cite,wrapfig,multicol,float,calc,geometry,tikz-3dplot,hhline}

\usepackage[colorlinks]{hyperref}
\usepackage[subnum]{cases}
\usepackage{xcolor}
\usepackage{enumitem}
\usepackage[T1]{fontenc}
\usepackage{ wasysym }
\usepackage{stmaryrd}

\usepackage{tikz}
\usetikzlibrary{decorations.markings}
\usetikzlibrary{arrows,matrix}
\usepgflibrary{arrows}
\usetikzlibrary{matrix,arrows,positioning,automata}
\pgfdeclarelayer{background layer}
\pgfsetlayers{background layer,main}

\tikzset{->-/.style={decoration={
  markings,
  mark=at position #1 with {\arrow{>}}},postaction={decorate}}}

\usetikzlibrary{intersections,through,hobby}

\usetikzlibrary{calc}

\newcommand{\nc}{\newcommand}
\nc{\rnc}{\renewcommand}

\nc\PV{\mathcal{PV}}
\rnc\L{\mathcal L}
\nc\R{\mathcal R}
\nc\PO{\P\O}
\nc\OI{\O\I}
\nc\IB{\mathcal{IB}}
\nc\TL{\mathcal T\!\mathcal L}
\nc\V{\mathcal V}
\nc\C{\mathcal C}
\nc\PT{\P\T}
\rnc\O{\mathcal O}
\nc\B{\mathcal B}
\nc\T{\mathcal T}
\nc\PlP{\mathscr P\mathcal P}
\nc\I{\mathcal I}
\renewcommand{\P}{\mathcal P} 
\renewcommand{\S}{\mathbb{S}}
\renewcommand{\SS}{\mathcal{S}}

\newcommand{\N}{\mathbb{N}}
\newcommand{\PP}{\mathbb{P}}
\nc\RRR{\mathbb R}
\nc\kk\Bbbk

\nc\ve\varepsilon
\nc\io\iota
\nc\al\alpha
\nc\be\beta
\nc\ga\gamma
\nc\de\delta
\nc\si\sigma
\nc\lam\lambda
\nc\ze\zeta
\rnc\th\theta

\nc\Ga\Gamma
\nc\De\Delta
\nc\Om\Omega

\nc\sib{\overline\si}
\nc\veb{\overline\ve}
\nc\taub{\overline\tau}
\nc\lamb{\overline\lam}
\nc\rhob{\overline\rho}
\nc\iob{\overline\io}
\nc\etab{\overline{\eta}}
\nc\mub{\overline{\mu}}

\nc\bd{{\bf d}}
\nc\br{{\bf r}}
\nc\s{\mathfrak s}
\nc\U{{\raise1.6 ex\hbox{\rotatebox{180}{$U$}}}}

\nc{\set}[2]{\{#1:#2\}}
\nc{\bigset}[2]{\big\{#1:#2\big\}}
\nc{\pres}[2]{\la#1:#2\ra}
\nc\Mod[1]{\ (\operatorname{mod}\ #1)}

\newcommand{\dom}{\operatorname{dom}} 
\nc\im{\operatorname{im}}

\nc\AND{\qquad\text{and}\qquad}
\nc\ANDSIM{\qquad\text{and similarly}\qquad}
\nc{\COMMA}{,\qquad}
\nc{\COMMa}{,\quad}

\newcommand{\sm}{\setminus}
\newcommand{\restr}{{\restriction}}
\rnc\emptyset\varnothing
\nc\op\oplus
\nc\sub\subseteq
\nc\mt\mapsto
\nc\la\langle
\nc\ra\rangle
\rnc\iff{\ \Leftrightarrow\ }
\rnc\implies{\ \Rightarrow\ }
\nc\tsharp{_{\scriptscriptstyle{\op}}^\sharp}

\nc\bit{\begin{itemize}}
\nc\eit{\end{itemize}}
\nc\ben{\begin{enumerate}[label=\textup{(\roman*)},leftmargin=7mm]}
\nc\bena{\begin{enumerate}[label=\textup{(\alph*)},leftmargin=7mm]}
\nc\een{\end{enumerate}}

\nc\pf{\begin{proof}}
\nc\epf{\end{proof}}
\nc\epfres{\hfill\qed}
\nc\epfreseq{\tag*{\qed}}
\let\oldproofname=\proofname
\renewcommand{\proofname}{\rm\bf{\oldproofname}}

\nc{\pfitem}[1]{\medskip\noindent #1.}
\nc{\firstpfitem}[1]{#1.}
\nc{\pfcase}[1]{\medskip\noindent {\bf Case #1.}}

\newcommand{\uv}[1]{\fill (#1,2)circle(.17);}
\newcommand{\lv}[1]{\fill (#1,0)circle(.17);}
\newcommand{\uvs}[1]{{\foreach \x in {#1} { \uv{\x}}}}
\newcommand{\lvs}[1]{{\foreach \x in {#1} { \lv{\x}}}}
\newcommand{\darcx}[3]{\draw(#1,0)arc(180:90:#3) (#1+#3,#3)--(#2-#3,#3) (#2-#3,#3) arc(90:0:#3);}
\newcommand{\darc}[2]{\darcx{#1}{#2}{.4}}
\newcommand{\uarcx}[3]{\draw(#1,2)arc(180:270:#3) (#1+#3,2-#3)--(#2-#3,2-#3) (#2-#3,2-#3) arc(270:360:#3);}
\newcommand{\uarc}[2]{\uarcx{#1}{#2}{.4}}
\newcommand{\stline}[2]{\draw(#1,2)--(#2,0);}
\newcommand{\uvc}[2]{\fill[#2] (#1,2)circle(.17);}
\newcommand{\lvc}[2]{\fill[#2] (#1,0)circle(.17);}
\newcommand{\uvcs}[2]{{\foreach \x in {#1} { \uvc{\x}{#2}}}}
\newcommand{\lvcs}[2]{{\foreach \x in {#1} { \lvc{\x}{#2}}}}
\newcommand{\stlinec}[3]{\draw[#3](#1,2)--(#2,0);}
\newcommand{\darcxc}[4]{\draw[#4](#1,0)arc(180:90:#3) (#1+#3,#3)--(#2-#3,#3) (#2-#3,#3) arc(90:0:#3);}
\newcommand{\darcc}[3]{\darcxc{#1}{#2}{.4}{#3}}
\newcommand{\uarcxc}[4]{\draw[#4](#1,2)arc(180:270:#3) (#1+#3,2-#3)--(#2-#3,2-#3) (#2-#3,2-#3) arc(270:360:#3);}
\newcommand{\uarcc}[3]{\uarcxc{#1}{#2}{.4}{#3}}
\nc\udotted[2]{\draw[dotted](#1+.5,2)--(#2-.5,2);}
\nc\ddotted[2]{\draw[dotted](#1+.5,0)--(#2-.5,0);}
\nc\understringx[4]{\draw[thick] (#1,#2) .. controls (#1,#2/2+#4/2) and (#3,#2/2+#4/2) .. (#3,#4);}
\nc\overstringx[4]{\draw[white,line width=2mm] (#1,#2) .. controls (#1,#2/2+#4/2) and (#3,#2/2+#4/2) .. (#3,#4); \draw[thick] (#1,#2) .. controls (#1,#2/2+#4/2) and (#3,#2/2+#4/2) .. (#3,#4);}
\nc\ststring[2]{\draw[white,line width=2mm] (#1,2)--(#2,0); \draw[thick] (#1,2)--(#2,0);}

\nc\ul\underline
\nc\ol\overline
\nc\wh\widehat
\nc\ulwh[1]{\ul{\hspace{0.7truemm}\wh{#1}\hspace{0.7truemm}}}

\DeclareMathSymbol{\widetildesym}{\mathord}{largesymbols}{"65}
\newcommand\lowerwidetildesym{%
  \text{\smash{\raisebox{-1.3ex}{%
    $\widetildesym$}}}}
\newcommand\fixwidetilde[1]{%
  \mathchoice
    {\accentset{\displaystyle\lowerwidetildesym}{#1}}
    {\accentset{\textstyle\lowerwidetildesym}{#1}}
    {\accentset{\scriptstyle\lowerwidetildesym}{#1}}
    {\accentset{\scriptscriptstyle\lowerwidetildesym}{#1}}
}
\nc\wt\fixwidetilde

\numberwithin{equation}{section}

\newtheorem{thm}[equation]{Theorem}
\newtheorem{lemma}[equation]{Lemma}

\theoremstyle{definition}

\newtheorem{rem}[equation]{Remark}

\newtheorem{ass}{Assumption}

\begin{document}

\title{\vspace{-1.5cm}Presentations for tensor categories}
\author{James East\footnote{The author is supported by ARC Future Fellowship FT190100632.  He thanks Richard Garner, Peter Hines, Steve Lack and Paul Martin for some helpful conversations and correspondence, and is very grateful to the anonymous referee for their valuable comments.}\\
{\it\small Centre for Research in Mathematics and Data Science,}\\
{\it\small Western Sydney University, Locked Bag 1797, Penrith NSW 2751, Australia.}\\
{\tt\small J.East\,@\,WesternSydney.edu.au}}
\date{}

\maketitle

\vspace{-.5cm}

\begin{abstract}
We introduce new techniques for working with presentations for a large class of (strict) tensor categories.  We then apply the general theory to obtain presentations for partition, Brauer and Temperley-Lieb categories, as well as several categories of (partial) braids, vines and mappings.



\end{abstract}

\thanks{Dedicated to the memory of Profs.~John H.~Conway and Mark V.~Sapir.}

\tableofcontents

\section{Introduction}\label{s:intro}

Presentations are important tools in many areas of mathematics.  They break complex structures into simple building blocks (generators) and local moves (relations) that determine equivalence of representations by generators.  An apt analogy are the Reidemeister moves from knot theory \cite{Reidemeister1927}, cf.~\cite{Birman1993,Alexander1923,Markov1935,Kauffman2005}, which are closely related to tensor presentations of tangle categories \cite{Halverson1996,Turaev1989,Kauffman1990,Conway1970}.  

Presentations of various kinds of algebraic structures feature prominently in many studies of knots, links, braids and related geometrical/topological objects; see for example 
\cite{
%Kock2004,
%Stroppel2005,
Kauffman1990,
%Kauffman1987,
Jones1987,
Jones1983_2, 
%Jones2001,
Birman1993,
%Baez1992,
Artin1947,
%Kauffman1987b, 
%JE2007a,
%JE2007b,
Lavers1997,
EL2004,
Paris2004, 
Dehornoy1994,
%DEEFHHL2015,
%DEEFHHLM2019,
Deligne2007
}.  
A key structure appearing in many such studies is the Temperley-Lieb category (and associated algebras and monoids), which also arises in many other areas of mathematics and science, especially theoretical physics \cite{Jones1994_a,Jones1994,TL1971}.  As with other so-called diagram categories, the Temperley-Lieb category has a basis consisting of certain kinds of set partitions, which are represented and composed diagrammatically.  The larger partition categories were introduced independently by Jones \cite{Jones1994} and Martin \cite{Martin1994} in the 1990s, again mostly motivated by problems in theoretical physics, but their applications extend in many other directions, including to representation theory, topology, knot theory, logic, combinatorics, theoretical computer science and more.  See for example 
\cite{
%ER2023,
%AV2020,
MM2013,
%BH2019,
%HR2005,
Brauer1937,
TL1971,
%Jones1994,
%Jones2001,
Jones1994_a,
%Jones1987,
%Jones1983_2,
GL1998, 
%Abramsky2008, 
%BS2011,
%Xi1999,
%Wenzl1988,
Wilcox2007,
Stroppel2005,
EMRT2018,
%EG2017, 
%Kauffman1990,
%Kauffman1987,
%BD1995,
%LZ2017,
%LZ2012,
LZ2015,
Martin1991,
Martin1994,
DP2013
} and references therein, and especially~\cite{Martin2008} for an excellent survey focussing particularly on connections with statistical mechanics.  Some of the ideas behind these categories go back to important early papers of Schur \cite{Schur} and Brauer~\cite{Brauer1937} on classical groups and invariant theory; cf.~\cite{LZ2012,LZ2015,LZ2017,BH2019}.
Presentations play an important role in many/most of the above studies.  

Presentations for the partition algebras were given by Halverson and Ram in \cite{HR2005}, along with a sketch of a proof; full proofs may be found in \cite{JEgrpm,JEgrpm2}.  The techniques used in \cite{JEgrpm,JEgrpm2} are semigroup-theoretic in nature; results concerning the algebras are deduced from corresponding results on associated diagram monoids, via the theory of twisted semigroup algebras~\cite{Wilcox2007}.  For other diagram monoids, see \cite{JE2021,BDP2002,MM2007,KM2006,HLP2013} and references therein.  Presentations are also crucial tools in Lehrer and Zhang's recent work on invariant theory \cite{LZ2012,LZ2015,LZ2017}.  For example, the first main result in \cite{LZ2015} is a presentation for the Brauer category.  Presentations for the Temperley-Lieb category may be found in a number places.  See especially \cite[Section 3]{Abramsky2008} for an interesting discussion; among other things, it is explained there that the Temperley-Lieb category is the free (so-called) pivotal category over one self-dual generator, a fact that follows from (the well-known) Theorem~\ref{t:TL2} below, and which is attributed in \cite{Abramsky2008} to \cite{DP2003,FY1989}.  A presentation for the partition category is given by Comes in \cite{Comes2020}, building on deep work of Abrams \cite{Abrams1996} and Kock \cite{Kock2004} on Frobenius algebras and cobordism categories; see also \cite{DP2013}, which treats a number of categories of relations.  Further discussion, and sketches of proofs in some cases, can be found in \cite{Hu2019}; see also \cite{Pilon2021} on oriented versions of some categories.

Many proofs have been given in the literature for presentations of the above-mentioned categories, algebras and monoids.  The techniques vary widely, as does the level of mathematical rigour in the arguments; for discussions of the latter point, see the introductions to \cite{BDP2002} or \cite{JEgrpm}.  One of the main goals of the current work is to introduce a systematic (and completely rigorous) framework for dealing with presentations for categories such as those mentioned above.  Although diagram categories were the original source of motivation for the current work, our methods apply to far wider classes of categories.  An additional family of applications come from the important class of braided monoidal categories \cite{BS2011, BJS2021, BD1995,JS1993,FY1989,BN1996}, and we are able to obtain presentations for a number of categories of (partial) braids and vines \cite{JE2007a,JE2007b,EL2004,Lavers1997}.  We also treat several categories of mappings, which are of course among the most important and foundational.  

We begin in Section \ref{s:C} by proving a number of general results concerning a large class of tensor categories, the key structural features of which are axiomatised in (subsets of) Assumptions \ref{a:1}--\ref{a:8}.  These results show how presentations for such categories can be built from presentations for endomorphism monoids and certain one-sided units (Theorem \ref{t:pres}), and then how to rewrite these into tensor presentations (Theorems \ref{t:pres2} and~\ref{t:pres3}).  Section \ref{s:DC} concerns diagram categories, and applies the general machinery of Section \ref{s:C} to quickly obtain presentations for the partition category (Theorems~\ref{t:P1} and~\ref{t:P2}), the Brauer category (Theorems~\ref{t:B1} and~\ref{t:B2}), and the Temperley-Lieb category (Theorems~\ref{t:TL1} and~\ref{t:TL2}); we also explain in Section \ref{ss:L} how to convert each such result into a presentation for the associated linear diagram category.  
In particular, we recover tensor presentations of Comes \cite{Comes2020} for the partition category (Theorem~\ref{t:P2}), and Lehrer and Zhang \cite{LZ2015} for the Brauer category (Theorem~\ref{t:B2}), with complete self-contained proofs using the uniform framework outlined in Section \ref{s:C}.  (To the best of our knowledge, all other presentations and tensor presentations are new, with the exception of the tensor presentation of the Temperley-Lieb category (Theorem \ref{t:TL2}), which is regarded as folklore.)
Section \ref{s:BC} treats a number of categories related to Artin braid groups: the partial vine category (Theorems~\ref{t:PV1} and~\ref{t:PV2}), the partial braid category (Theorems~\ref{t:IB1} and~\ref{t:IB2}), and the (full) vine category (Theorems~\ref{t:V1} and~\ref{t:V2}).  Sections \ref{ss:T} and \ref{ss:O} apply these results on partial braids and vines to many categories of transformations/mappings.  All in all, we give presentations for a dozen concrete categories, but due to the scope of the general results of Section \ref{s:C}, the potential for further applications is vast.

\section{Categories}\label{s:C}

This section develops a general theory of presentations for a broad class of (strict) tensor categories satisfying natural connectivity assumptions, and possessing certain one-sided units.  After fixing notation and gathering some basic ideas in Sections \ref{ss:cat} and \ref{ss:pre_pres}, we give the two main general results in Sections~\ref{ss:pres} and~\ref{ss:pres2}.  The first of these (Theorem \ref{t:pres}) shows how to construct a presentation for a category out of presentations for its endomorphism monoids, under the assumptions alluded to above.  The second (Theorem \ref{t:pres2}) then shows how to convert this to a tensor presentation.  In Section \ref{ss:pres3} we prove a variation of the second main result (see Theoem~\ref{t:pres3}), tailored to work for categories whose connectivity is more complicated, and motivated by three of the applications considered later in the paper.

Before we begin, we note that Higgins considered the relationship between presentations for (fundamental) groupoids and their vertex groups in \cite{Higgins1964}.  We also observe that there is a body of literature showing how to obtain presentations for monoids/semigroups constructed from other families by natural constructions, such as (semi)direct products or wreath products; see for example \cite{DR2009,HR1994,RRW1998}.  In a sense, our Theorem \ref{t:pres} could be thought of as belonging to this latter field of study.

\subsection{Preliminaries}\label{ss:cat}

All categories we consider will have object set $\N=\{0,1,2,\ldots\}$ or $\PP=\{1,2,3,\ldots\}$, and we write $\S$ to stand for either of these.  Moreover, the collections of morphisms $m\to n$ ($m,n\in\S$) will always be a set.  Thus, for simplicity throughout, ``category'' will always mean ``small category with object set $\S$ (being $\N$ or $\PP$)''.  We will always identify such a category $\C$ with its set of morphisms.  For $a\in \C$, we write $\bd(a)$ and $\br(a)$ for the domain and range of~$a$.  We compose morphisms left-to-right, and often supress the composition symbol, so that for $a,b\in \C$, $ab=a\circ b$ is defined if and only if $\br(a)=\bd(b)$, in which case $\bd(ab)=\bd(a)$ and $\br(ab)=\br(b)$.  For $m,n\in\S$ we write
\[
\C_{m,n}=\set{a\in \C}{\bd(a)=m,\ \br(a)=n}
\]
for the (possibly empty) set of all morphisms $m\to n$.  For $n\in\S$ we write $\C_n=\C_{n,n}$ for the endomorphism monoid at $n$.

As in \cite[p.~52]{MacLane1998}, a congruence on a category $\C$ is an equivalence $\xi$ on (morphisms of) $\C$ that preserves objects and is compatible with composition: i.e., for all $(a,b)\in\xi$ and all $c,d\in \C$, 
\bit
\item $\bd(a)=\bd(b)$ and $\br(a)=\br(b)$,
\item $(ac,bc)\in\xi$ and $(da,db)\in\xi$, whenever the stated products are defined.
\eit
The quotient category $\C/\xi$ consists of all $\xi$-classes under the induced composition.  For a set $\Om\sub \C\times \C$ of pairs satisfying $\bd(u)=\bd(v)$ and $\br(u)=\br(v)$ for all $(u,v)\in\Om$, we write $\Om^\sharp$ for the congruence on $\C$ generated by $\Om$: i.e., the smallest congruence containing $\Om$.

By a category morphism $\phi:\C\to\mathcal D$ we mean a functor that acts as the identity on objects: i.e., we have $\bd(a\phi)=\bd(a)$ and $\br(a\phi)=\br(a)$ for all $a\in \C$.  The kernel of $\phi$ is the congruence
\[
\ker(\phi) = \bigset{(a,b)\in \C\times \C}{a\phi=b\phi},
\]
and we have $\C/\ker(\phi)\cong\im(\phi)$.  A surmorphism is a surjective morphism.

The categories we consider are all strict tensor (a.k.a.~monoidal) categories, in the sense of \cite{JS1993}; see also \cite[Chapters~VII and XI]{MacLane1998}.  In fact, since all categories we consider have object set $\N$ or $\PP$, certain simplifications arise, and we take \cite[Section 24]{MacLane1965} as our reference.  

A (strict) tensor category is a category $\C$ (over $\S$) with an extra (totally defined) binary operation~$\op$ satisfying the following properties (writing $\io_n$ for the identity at~$n\in\S$):
\bit
\item $\bd(a\op b)=\bd(a)+\bd(b)$ and $\br(a\op b)=\br(a)+\br(b)$ for all $a,b\in \C$,
\item $a\op(b\op c)=(a\op b)\op c$ for all $a,b,c\in \C$,
\item $a\op \io_0=a=\io_0\op a$ for all $a\in \C$ (in the case that $\S=\N$),
\item $\io_m\op\io_n=\io_{m+n}$ for all $m,n\in\S$,
\item $(a\circ b)\op(c\circ d) = (a\op c)\circ(b\op d)$ for all $a,b,c,d\in \C$ with $\br(a)=\bd(b)$ and $\br(c)=\bd(d)$.
\eit
The symbol~$\otimes$ is often used in place of $\op$ in the literature, but since the object part of the operation is addition on $\S$, we prefer $\op$.

The term ``strict'' refers to the fact that we have equalities in the above axioms; in more general monoidal categories only isomorphism is required.  In what follows, we will generally drop the word ``strict'', as all tensor categories we consider are strict.  The following basic lemma will be used often.

\begin{lemma}\label{l:tc}
Let $a$, $b$ and $c$ be elements of a (strict) tensor category over $\N$.
\ben
\item \label{tc1} If $\br(a)=0$ and $\br(b)=\bd(c)$, then $a\op(b\circ c)=(a\op b)\circ c$.
\item \label{tc2} If $\bd(a)=0$ and $\br(b)=\bd(c)$, then $a\op(b\circ c)=b\circ(a\op c)$.
\item \label{tc3} If $\bd(c)=0$ and $\br(a)=\bd(b)$, then $(a\circ b)\op c=a\circ(b\op c)$.
\item \label{tc4} If $\br(c)=0$ and $\br(a)=\bd(b)$, then $(a\circ b)\op c=(a\op c)\circ b$.
\een
\end{lemma}

\pf
For \ref{tc1} we have $a\op(b\circ c) = (a\circ\io_0)\op(b\circ c) = (a\op b)\circ(\io_0\op c) = (a\op b)\circ c$.  The others are similar.
\epf

The next result can be deduced from the previous one, or can be similarly proved directly.

\begin{lemma}\label{l:ab}
If $a$ and $b$ are elements of a (strict) tensor category over $\N$, and if $\br(a)=\bd(b)=0$, then $a\circ b=a\op b=b\op a$.  \epfres
\end{lemma}

A tensor congruence on a tensor category $\C$ is a congruence $\xi$ (as above) such that for all $(a,b)\in\xi$, and for all $c\in\C$, we have $(a\op c,b\op c),(c\op a,c\op b)\in\xi$.  For a set $\Om\sub \C\times \C$ of pairs satisfying $\bd(u)=\bd(v)$ and $\br(u)=\br(v)$ for all $(u,v)\in\Om$, we write $\Om\tsharp$ for the tensor congruence on $\C$ generated by $\Om$.

Some of the tensor categories we consider have further structure.  Namely, the endomorphism monoids $\C_n=\C_{n,n}$ ($n\in\S$) contain natural copies of the symmetric groups $\SS_n$, which satisfy the so-called PROP axioms, as in \cite[Section 24]{MacLane1965}.  (The term PROP is a contraction of PROducts and Permutations.)  Since the PROP structures of our categories will not play an explicit role, we will not give the full details here.  But some rough ideas are worth noting.  First, for $m,n\in\S$, and for $f\in\SS_m\sub \C_m$ and $g\in\SS_n\sub \C_n$, the permutation $f\op g\in\SS_{m+n}\sub \C_{m+n}$ acts as expected: as $f$ on the ``first $m$ points'' and as $g$ (suitably translated) on the ``last $n$ points''.  Second, for each $m,n\in\S$ there is a permutation $f_{m,n}\in\SS_{m+n}$ such that for all $m,n,k,l\in\S$, and for all $a\in \C_{m,n}$ and $b\in \C_{k,l}$, $(a\op b)\circ f_{n,l} = f_{m,k}\circ(b\op a)$.  Other categories we consider are PROBS, which have permutations replaced by braids \cite{JS1993,Artin1947}.

\subsection{Free (tensor) categories and presentations}\label{ss:pre_pres}

Let $X$ be an alphabet (a set whose elements are called letters).  The free monoid over $X$ is the set $X^*$ of all words over $X$ under concatenation.  The empty word $\io$ is the identity of $X^*$.  Let $R\sub X^*\times X^*$ be a set of pairs of words over $X$.  We say a monoid $M$ has presentation $\pres XR$ if $M\cong X^*/R^\sharp$ (where, as above, $R^\sharp$ is the congruence on $X^*$ generated by $R$): i.e., if there exists a monoid surmorphism $X^*\to M$ with kernel $R^\sharp$.  If $\phi$ is such a surmorphism, we say $M$ has presentation $\pres XR$ via $\phi$.  Elements of $X$ and $R$ are called generators and relations, respectively.  A relation $(u,v)\in R$ is typically displayed as an equation: $u=v$.

For category presentations, digraphs and paths play the role of alphabets and words.  Let $\Ga$ be a digraph with vertex set $\S$ (being $\N$ or $\PP$), possibly with multiple/parallel edges, and possibly with loops.  We identify $\Ga$ with its edge set, and denote the source and target of an edge $x\in\Ga$ by $\bd(x)$ and~$\br(x)$, respectively.  The free category over $\Ga$ is the set~$\Ga^*$ of all paths in $\Ga$ under concatenation (where defined).  The empty path at $n\in\S$ will be denoted by~$\io_n$.  Every other path can be thought of as a word of the form $w=x_1\cdots x_k$, where $k\geq1$ and $x_i\in\Ga$ for all $1\leq i\leq k$, and where $\br(x_i)=\bd(x_{i+1})$ for all $1\leq i<k$.  For such a word/path, we have $\bd(w)=\bd(x_1)$ and $\br(w)=\br(x_k)$.  Now let $\Om\sub\Ga^*\times\Ga^*$ be a set of pairs of paths, such that $\bd(u)=\bd(v)$ and $\br(u)=\br(v)$ for all $(u,v)\in\Om$.  We say a category~$\C$ (over~$\S$) has presentation $\pres\Ga\Om$ if~$\C\cong\Ga^*/\Om^\sharp$: i.e., if there exists a surmorphism $\Ga^*\to \C$ with kernel~$\Om^\sharp$.  If $\phi$ is such a surmorphism, we say $\C$ has presentation $\pres \Ga\Om$ via $\phi$.  

There is also a natural notion of a (strict) tensor category presentation.  Let $\De$ be a digraph with vertex set $\S$, again identified with its edge set.  We will denote the free tensor category over~$\De$ by~$\De^\oast$.  It consists of all terms constructed in the following way:
\begin{enumerate}[label=\textup{(T\arabic*)},leftmargin=12mm]
\item \label{T1} All empty paths $\io_n$ ($n\in\S$) are terms, with $\bd(\io_n)=\br(\io_n)=n$.
\item \label{T2} All edges $x\in\De$ are terms, with $\bd(x)$ and $\br(x)$ the source and target of $x$, respectively.
\item \label{T3} If $s$ and $t$ are terms, and if $\br(s)=\bd(t)$, then the formal expression $s\circ t$ is a term, with $\bd(s\circ t)=\bd(s)$ and $\br(s\circ t)=\br(t)$.  If $\bd(s)=m$ and $\br(s)=n$, then $\io_m\circ s=s=s\circ\io_n$.
\item \label{T4} If $s$ and $t$ are terms, then the formal expression $s\op t$ is a term, with $\bd(s\op t)=\bd(s)+\bd(t)$ and $\br(s\op t)=\br(s)+\br(t)$.
\een
Note that \ref{T1}--\ref{T4} describe the elements of $\De^\oast$, while \ref{T3} and \ref{T4} also give the definition of the $\circ$ and $\op$ operations.  Now let $\Xi\sub\De^\oast\times\De^\oast$ be a set of pairs of terms, such that $\bd(u)=\bd(v)$ and $\br(u)=\br(v)$ for all $(u,v)\in\Xi$.  We say a tensor category $\C$ (over $\S$) has tensor presentation $\pres\De\Xi$ if~$\C\cong\De^\oast/\Xi\tsharp$: i.e., if there is a surmorphism $\De^\oast\to \C$ with kernel~$\Xi\tsharp$.  If $\phi$ is such a surmorphism, we say $\C$ has presentation $\pres\De\Xi$ via $\phi$.

An arbitrary term from $\De^\oast$, as in \ref{T1}--\ref{T4}, can be quite unwieldy.  However, every such term $w$ is equal in $\De^\oast$ to one of the form
\begin{equation}\label{e:term}
X_1\circ\cdots\circ X_m = (x_{1,1}\op\cdots\op x_{1,k_1}) \circ\cdots\circ (x_{m,1}\op\cdots\op x_{m,k_m}),
\end{equation}
where each $x_{i,j}$ is either an empty path $\io_n$ ($n\in\S$) or else an edge of $\De$, and where the domains and ranges of the sub-terms $X_q = x_{q,1}\op\cdots\op x_{q,k_q}$ are compatible as appropriate.  Indeed, this is clear if~$w$ is itself of the form \ref{T1} or \ref{T2} as above.  Otherwise, $w$ is of the form $s\circ t$ or $s\op t$, as in \ref{T3} and~\ref{T4}, and inductively $s$ and $t$ are both of the form \eqref{e:term}.  It immediately follows that $s\circ t$ is of the desired form.  To see that $s\op t$ also is, write $s = X_1\circ\cdots\circ X_m$ and $t = Y_1\circ\cdots\circ Y_n$, as in \eqref{e:term}.  Replacing $s$ by $s\circ\io_r\circ\cdots\circ\io_r$ ($r=\br(s)$), or  $t$ by $t\circ\io_r\circ\cdots\circ\io_r$ ($r=\br(t)$), if necessary, we may assume that $m=n$.  We then have $s\op t = (X_1\circ\cdots\circ X_n)\op(Y_1\circ\cdots\circ Y_n) = (X_1\op Y_1)\circ\cdots\circ(X_n\op Y_n)$.

In fact, every term from $\De^\oast$ is equal either to an empty path or to a term of the form
\begin{equation}\label{e:term2}
X_1\circ\cdots\circ X_m \qquad\text{where $m\geq1$ and each $X_i = \io_{p_i}\op x_i \op\io_{q_i}$ for some $p_i,q_i\in\N$ and $x_i\in\De$.}
\end{equation}
(In the above, if $\S=\PP$, but $p_i=0$ say, then we interpret $X_i = \io_0\op x_i \op\io_{q_i} = x_i \op\io_{q_i}$, and so on.)  Indeed, by \eqref{e:term}, it suffices to show that this is true of every term of the form ${w=y_1\op\cdots\op y_k}$, where each $y_i$ is either an empty path or else an edge of $\De$.  We now proceed by induction on~$\mu=\mu(w)$, defined to be the number of $i\in\{1,\ldots,k\}$ such that $y_i\in\De$.  For each $i$, write $r_i=\br(y_i)$.  If $\mu=0$, then $w=\io_{r_1}\op\cdots\op\io_{r_k}=\io_{r_1+\cdots+r_k}$.  If $\mu\geq1$, say with $y_i\in\De$, then with $d=\bd(y_i)$, we have
\begin{align*}
w &= (y_1\circ\io_{r_1})\op\cdots\op (y_{i-1}\circ\io_{r_{i-1}}) \op (\io_d\circ y_i) \op (y_{i+1}\circ\io_{r_{i+1}})\op\cdots\op(y_k\circ\io_{r_k}) \\
&= (y_1\op\cdots\op y_{i-1}\op\io_d\op y_{i+1}\op\cdots\op y_k) \circ (\io_{r_1}\op\cdots\op\io_{r_{i-1}} \op y_i \op \io_{r_{i+1}}\op\cdots\op\io_{r_k}) \\
&= (y_1\op\cdots\op y_{i-1}\op\io_d\op y_{i+1}\op\cdots\op y_k) \circ (\io_{r_1+\cdots+r_{i-1}} \op y_i \op \io_{r_{i+1}+\cdots+r_k}).
\end{align*}
Since $\mu(y_1\op\cdots\op y_{i-1}\op\io_d\op y_{i+1}\op\cdots\op y_k)=\mu(w)-1$, the inductive assumption completes the proof.

\subsection{First main result: category presentations}\label{ss:pres}

Our aim in this section is to prove Theorem \ref{t:pres}, which shows how to build a presentation for a category~$\C$ (with object set $\S$, being $\N$ or $\PP$) out of presentations for the endomorphism monoids $\C_n$ under certain natural assumptions stated below.  In what follows, we write $\iob_n\in\C$ for the identity at $n\in\S$.  The reason for this (and other) over-line notation will become clear shortly.

The first assumption concerns the concerns the connectivity of $\C$:

\begin{ass}\label{a:1}
We assume that there is an integer $d\geq1$ such that for all $m,n\in\S$,
\[
\C_{m,n} \not=\emptyset \iff m\equiv n\Mod d.
\]
\end{ass}

In most of our applications, we will have $d=1$, so that all hom-sets are non-empty; in some cases we have $d=2$; in a small number of exceptional cases, we have a more complicated arrangement (see Section \ref{ss:pres3}).  The next assumption asserts the existence of certain one-sided units, and is stated in terms of the integer $d$ from Assumption \ref{a:1}:

\begin{ass}\label{a:2}
We assume that for each $n\in\S$ there exist $\lamb_n\in\C_{n,n+d}$ and $\rhob_n\in\C_{n+d,n}$ such that
\[
\lamb_n\rhob_n = \iob_n.
\]
\end{ass}

For each $m,n\in\S$ with $m\leq n$ and $m\equiv n\Mod d$, we define
\begin{equation}\label{e:lambrhob}
\lamb_{m,n} = \lamb_m \lamb_{m+d}\cdots\lamb_{n-d} \in\C_{m,n} \AND \rhob_{n,m} = \rhob_{n-d}\cdots\rhob_{m+d}\rhob_m\in\C_{n,m},
\end{equation}
noting that $\lamb_{m,n}\rhob_{n,m}=\iob_m$.  (When $m=n$, we interpret these empty products as $\lamb_{m,m}=\rhob_{m,m}=\iob_m$.)  For $m,n$ as above, we define maps
\begin{equation}\label{e:RL}
\R_{m,n}:\C_{m,n}\to\C_n:a\mt\rhob_{n,m}a \AND \L_{n,m}:\C_{n,m}\to\C_n:a\mt a\lamb_{m,n}.
\end{equation}
For any $a\in\C_{m,n}$ and $b\in\C_{n,m}$, we have
\begin{equation}\label{e:lR}
a=\iob_ma=\lamb_{m,n}\rhob_{n,m}a=\lamb_{m,n}\cdot\R_{m,n}(a) \ANDSIM b=\L_{n,m}(b)\cdot\rhob_{n,m}.
\end{equation}

Our next assumption essentially just fixes notation for presentations for the endomorphism monoids in~$\C$.

\begin{ass}\label{a:3}
For each $n\in\S$, we assume that the endomorphism monoid $\C_n$ has presentation $\pres{X_n}{R_n}$ via $\phi_n:X_n^*\to\C_n$.  Further, we assume the alphabets $X_n$ ($n\in\S$) are pairwise disjoint, and denote by $\io_n$ the empty word over $X_n$.
\end{ass}

Now let $\Ga$ be the digraph with vertex set $\S$, and edge set 
\begin{equation}\label{e:Ga}
L\cup R\cup X \qquad\text{where}\qquad L=\set{\lam_n}{n\in\S} \COMMA R=\set{\rho_n}{n\in\S} \COMMA X = \bigcup_{n\in\S}X_n,
\end{equation}
with sources and targets given by
\begin{equation}\label{e:st}
\bd(x)=\br(x)=\bd(\lam_n)=\br(\rho_n)=n \AND \br(\lam_n)=\bd(\rho_n)=n+d \qquad\text{for all $n\in\S$ and $x\in X_n$.}
\end{equation}
As usual, we identify $\Ga$ with its edge set: $\Ga\equiv L\cup R\cup X$.  We have a morphism 
\begin{equation}\label{e:phi}
\phi:\Ga^*\to\C \qquad\text{given by}\qquad \lam_n\phi = \lamb_n \COMMA \rho_n\phi = \rhob_n \COMMA x\phi=x\phi_n \qquad\text{for $n\in\S$ and  $x\in X_n$.}
\end{equation}
We extend the over-line notation to words/paths over $\Ga$, writing $\ol w=w\phi$ for all $w\in\Ga^*$.

For $m,n\in\S$ with $m\leq n$ and $m\equiv n\Mod d$, we define the words/paths
\[
\lam_{m,n} = \lam_m \lam_{m+d}\cdots\lam_{n-d} \in\Ga^*_{m,n} \AND \rho_{n,m} = \rho_{n-d}\cdots\rho_{m+d}\rho_m\in\Ga^*_{n,m},
\]
where again we interpret $\lam_{m,m}=\rho_{m,m}=\io_m$.  Note that $\lam_n=\lam_{n,n+d}$ and $\rho_n=\rho_{n+d,n}$ for all $n$, and that $\lam_{l,m}\lam_{m,n}=\lam_{l,n}$ and $\rho_{n,m}\rho_{m,l}=\rho_{n,l}$ for appropriate $l,m,n$.

\begin{lemma}\label{l:surj}
The morphism $\phi:\Ga^*\to\C$ is surjective.
\end{lemma}

\pf
Let $a\in\C$, and write $m=\bd(a)$ and $n=\br(a)$; we assume that $m\leq n$, the other case being symmetrical.  By~\eqref{e:lR}, we have $a=\lamb_{m,n}\cdot\R_{m,n}(a)$.  Since $\R_{m,n}(a)\in\C_n$, and since $\phi_n$ is surjective (cf.~Assumption~\ref{a:3}), we have $\R_{m,n}(a)=\ol w$ for some $w\in X_n^*$.  Thus, $a = \lamb_{m,n}\ol w = (\lam_{m,n}w)\phi$.
\epf

The final assumption constructs an appropriate set of relations, and is again stated in terms of the integer~$d$ from Assumption \ref{a:1}:

%\newpage

\begin{ass}\label{a:4}
We assume that $\Om\sub\Ga^*\times\Ga^*$ is a set of relations such that, writing ${\sim}$ for the congruence $\Om^\sharp$ on $\Ga^*$, the following all hold:
\ben
\item \label{a41} For every relation $(u,v)\in\Om$, we have $\ol u=\ol v$.
\item \label{a42} For all $n\in\S$, $\Om$ contains the relation $\lam_n\rho_n=\io_n$, and a relation of the form $\rho_n\lam_n=w_n$ for some word $w_n\in X_{n+d}^*$.
\item \label{a43} For all $n\in\S$, $R_n\sub\Om$.
\item \label{a44} For all $n\in\S$, there are mappings
\[
X_n\to X_{n+d}^*:x\mt x_+ \AND X_n\to X_{n+d}^*:x\mt x^+,
\]
and $\Om$ contains the relations
\[
x\lam_n = \lam_nx_+ \AND \rho_nx = x^+\rho_n \qquad\text{for all $x\in X_n$.}
\]
\item \label{a45} For all $n\in\S$, and for all $w\in X_{n+d}^*$, we have $\lam_nw\rho_n\sim u$ for some $u\in X_n^*$.
\een
\end{ass}

\begin{rem}\label{r:a4}
In practice, there could be several choices of the one-sided units 
$\lamb_n,\rhob_n$ 
%$\lamb_n\in\C_{n,n+d}$ and $\rhob_n\in\C_{n+d,n}$ 
from Assumption~\ref{a:2}, leading to different relations of the form \ref{a42} and \ref{a44} above.  In the applications in Sections~\ref{s:DC} and \ref{s:BC}, we have made choices that we believe lead to the simplest and most convenient versions of these relations.  See also Remarks \ref{r:a4'} and \ref{r:a4''}.
\end{rem}

For the remainder of Section \ref{ss:pres}, we assume that Assumptions \ref{a:1}--\ref{a:4} all hold, and we continue to write ${\sim}=\Om^\sharp$.  Our ultimate goal is to show that $\C$ has presentation $\pres\Ga\Om$ via $\phi$; see Theorem \ref{t:pres}.  

In what follows, we extend the maps from Assumption \ref{a:4}\ref{a44} to morphisms
\[
X_n^*\to X_{n+d}^*:w\mt w_+ \AND X_n^*\to X_{n+d}^*:w\mt w^+.
\]
It quickly follows that $w\lam_n \sim \lam_nw_+$ and $\rho_nw \sim w^+\rho_n$ for all $w\in X_n^+$.

\begin{rem}\label{r:a45}
When working with a specific category $\C$, item \ref{a45} in Assumption \ref{a:4} generally involves the most work to verify.  One way to do so is to show that for any $n\in\S$ and $w\in X_{n+d}^*$, we have $w_nww_n\sim w_nu^+w_n$ or $w_nww_n\sim w_nu_+w_n$ for some $u\in X_n^*$; here $w_n$ is the word from Assumption~\ref{a:4}\ref{a42}.  Indeed, suppose $w_nww_n\sim w_nu^+w_n$ holds (with the other case being similar).  First note that item \ref{a42} gives
\[
\lam_n = \io_n\lam_n \sim \lam_n\rho_n\lam_n \sim \lam_nw_n \ANDSIM \rho_n \sim w_n\rho_n.
\]
Combining these with \ref{a42}, \ref{a44} and $w_nww_n\sim w_nu^+w_n$, we obtain
\[
\lam_nw\rho_n \sim \lam_n w_nww_n\rho_n \sim \lam_nw_nu^+w_n\rho_n \sim \lam_nu^+\rho_n \sim \lam_n\rho_nu\sim\io_nu=u.
\]
\end{rem}

We now prove a sequence of lemmas, building up to Lemma \ref{l:w'2}, which gives a set of normal forms that are crucial in the proof of the main result.

\begin{lemma}\label{l:lr2}
If $m,n\in\S$ are such that $m\leq n$ and $m\equiv n\Mod d$, then 
\ben
\item \label{lr21} $\lam_{m,n}\rho_{n,m}\sim\io_m$,
\item \label{lr22} $\rho_{n,m}\lam_{m,n}\sim w$ for some $w\in X_n^*$.
\een
\end{lemma}

\pf
We just prove \ref{lr22}, by induction on $n$; \ref{lr21} is similar but easier.  The claim being clear if $n=m$ (take $w=\io_n$), we assume that $n>m$.  Then
\begin{align*}
\rho_{n,m}\lam_{m,n} = \rho_{n-d}\rho_{n-d,m}\lam_{m,n-d}\lam_{n-d} &\sim \rho_{n-d}v\lam_{n-d} &&\text{by induction, for some $v\in X_{n-d}^*$}\\
&\sim \rho_{n-d}\lam_{n-d}v_+ &&\text{by Assumption \ref{a:4}\ref{a44}}\\
&\sim w_{n-d}v_+ &&\text{by Assumption \ref{a:4}\ref{a42},}
\end{align*}
so we take $w=w_{n-d}v_+$.
\epf

For the next proof, given $n\in\S$ and $k\in\PP$, we define a mapping $X_n^*\to X_{n+kd}^*:w\mt w^{+^k}$ in the obvious way, by iteratively composing all the individual $w\mt w^+$ mappings ($X_n^*\to X_{n+d}^*$, $X_{n+d}^*\to X_{n+2d}^*$, and so on).

\begin{lemma}\label{l:w'}
Let $w\in\Ga^*$, with $m=\bd(\ol w)$ and $n=\br(\ol w)$.  
\ben
\item \label{w'1} If $m\leq n$, then $w\sim \lam_{m,n}w'$ for some $w'\in X_n^*$.
\item \label{w'2} If $m\geq n$, then $w\sim w'\rho_{m,n}$ for some $w'\in X_m^*$.
\een
\end{lemma}

\pf
We prove the lemma by induction on $k$, the length of $w$.  If $k=0$, then $m=n$ and we have $w=\lam_{m,n}=\rho_{m,n}=\io_m$, and \ref{w'1} and \ref{w'2} both hold with $w'=\io_m$.  Now suppose $k\geq1$, and inductively assume that the lemma holds for words of length less than $k$.  Write $w=x_1\cdots x_k$, where each $x_i\in\Ga$.  For simplicity, we write $u=x_1\cdots x_{k-1}$ and $x=x_k$.  Note that $\bd(u)=m$ and $\br(x)=n$.  We also write $q=\br(u)=\bd(x)$.  By induction, since $u$ has length $k-1$, one of the following holds:
\bena
\item \label{a} $m\leq q$, and $u\sim\lam_{m,q}s$ for some $s\in X_q^*$, or
\item \label{b} $m>q$, and $u\sim t\rho_{m,q}$ for some $t\in X_m^*$.  
\een
\pfcase1
First suppose $x\in X$, so that $q=n$.  
If \ref{a} holds then $w=ux\sim\lam_{m,n}sx$, with $sx\in X_n^*$.  
If~\ref{b} holds then, writing $m=n+kd$, we have $w=ux\sim t\rho_{m,n}x \sim tx^{+^k}\rho_{m,n}$, by Assumption \ref{a:4}\ref{a44}, with $tx^{+^k}\in X_m^*$.

\pfcase2
Next suppose $x\in L$, so that $q=n-d$ and $x=\lam_{n-d}$.  
If \ref{a} holds then from $s\in X_{n-d}^*$, it follows from Assumption \ref{a:4}\ref{a44} that $w = u\lam_{n-d} \sim \lam_{m,n-d}s\lam_{n-d} \sim \lam_{m,n-d}\lam_{n-d}s_+ = \lam_{m,n}s_+$, with $s_+\in X_n^*$.
If \ref{b} holds then, writing $m=n+kd$, and using Assumption \ref{a:4}\ref{a42} and \ref{a44}, we have $w = u\lam_{n-d} \sim t\rho_{m,n-d}\lam_{n-d} = t\rho_{m,n}\rho_{n-d}\lam_{n-d} \sim t\rho_{m,n}w_{n-d} \sim tw_{n-d}^{+^k}\rho_{m,n}$, with $tw_{n-d}^{+^k}\in X_m^*$.

\pfcase3
Finally, suppose $x\in R$, so that $q=n+d$ and $x=\rho_n$.
If \ref{a} holds then either
\bit
\item $m=q$, and $u\sim\lam_{m,m}s=\io_ms=s$, and so $w=u\rho_n\sim s\rho_n = s\rho_{n+d,n} = s\rho_{m,n}$, with $s\in X_m^*$, or
\item $m<q$, and by Assumption \ref{a:4}\ref{a45}, $w = u\rho_n \sim \lam_{m,n+d}s\rho_n = \lam_{m,n}\lam_ns\rho_n \sim \lam_{m,n}u$ for some $u\in X_n^*$.
\eit
If \ref{b} holds then $w=u\rho_n\sim t\rho_{m,n+d}\rho_n = t\rho_{m,n}$.
\epf

The next lemma strengthens the previous one, and is the main technical result we need.  The statement refers to the maps $\R_{m,n}$ and $\L_{m,n}$ defined in \eqref{e:RL}.

\begin{lemma}\label{l:w'2}
Let $w\in\Ga^*$, with $m=\bd(\ol w)$ and $n=\br(\ol w)$.  
\ben
\item \label{w'21} If $m\leq n$, then $w\sim \lam_{m,n}w'$ for some $w'\in X_n^*$ with $\ol w'\in\im(\R_{m,n})$.
\item \label{w'22} If $m\geq n$, then $w\sim w'\rho_{m,n}$ for some $w'\in X_m^*$ with $\ol w'\in\im(\L_{m,n})$.
\een
\end{lemma}

\pf
We assume that $m\leq n$, the other case being symmetrical.  First note that
\begin{align*}
w &\sim \lam_{m,n}v &&\text{by Lemma \ref{l:w'}, for some $v\in X_n^*$}\\
&= \io_m\lam_{m,n}v \sim \lam_{m,n}\rho_{n,m}\lam_{m,n}v &&\text{by Lemma \ref{l:lr2}\ref{lr21}.}
\end{align*}
By Lemma \ref{l:lr2}\ref{lr22}, we have $\rho_{n,m}\lam_{m,n}\sim u$ for some $u\in X_n^*$.  Then with $w'=uv$, it follows from the above calculations that $w\sim\lam_{m,n}w'$.  We also have $\ol w' = \ol{uv} = \rhob_{n,m}\lamb_{m,n}\ol v = \R_{m,n}(\lamb_{m,n}\ol v)$.
\epf

We now have all we need to prove our first main result:

\begin{thm}\label{t:pres}
With notation as above, and subject to Assumptions \ref{a:1}--\ref{a:4}, the category $\C$ has presentation $\pres\Ga\Om$ via $\phi$.
\end{thm}

\pf
We showed in Lemma \ref{l:surj} that $\phi$ is surjective.  It remains to show that $\ker(\phi)=\Om^\sharp$.  First note that Assumption \ref{a:4}\ref{a41} says $\Om\sub\ker(\phi)$; since $\ker(\phi)$ is a congruence it follows that $\Om^\sharp\sub\ker(\phi)$.  

For the reverse containment, suppose $(u,v)\in\ker(\phi)$, meaning that $u,v\in\Ga^*$ and $\ol u=\ol v$.  Write ${m=\bd(\ol u)=\bd(\ol v)}$ and $n=\br(\ol u)=\br(\ol v)$, and assume that $m\leq n$, the other case being symmetrical.  By Lemma~\ref{l:w'2}\ref{w'21}, we have $u\sim\lam_{m,n}u'$ and $v\sim\lam_{m,n}v'$ for some $u',v'\in X_n^*$ with $\ol u',\ol v'\in\im(\R_{m,n})$, say $\ol u'=\R_{m,n}(a)$ and $\ol v'=\R_{m,n}(b)$, where $a,b\in\C_{m,n}$.  Then by \eqref{e:lR} we have
\[
a = \lamb_{m,n}\cdot\R_{m,n}(a) = \lamb_{m,n}\ol u' = \ol u = \ol v = \lamb_{m,n}\ol v' = \lamb_{m,n}\cdot\R_{m,n}(b) = b.
\]
But then $\ol u' = \R_{m,n}(a) = \R_{m,n}(b) = \ol v'$.  Since $u',v'\in X_n^*$, it follows that $(u',v')\in\ker(\phi_n)=R_n^\sharp$; cf.~Assumption \ref{a:3}.  By Assumption \ref{a:4}\ref{a43}, it follows that $u'\sim v'$.  Putting all of this together, we deduce that $u\sim\lam_{m,n}u'\sim\lam_{m,n}v'\sim v$.
\epf

\begin{rem}
Theorem \ref{t:pres} shows how to build a presentation for the category $\C$ out of presentations for its endomorphism monoids $\C_n=\C_{n,n}$ (under appropriate assumptions).  
It would be interesting to study the extent to which one could work in the opposite direction: i.e., begin with a presentation for~$\C$, and deduce presentations for the monoids $\C_n$.  It does not seem this is likely to work in general, however.  For example, consider the category $\B$ of all binary relations $\{1,\ldots,m\}\to\{1,\ldots,n\}$, $m,n\in\N$.  Then~$\B$ can be generated by a relatively simple set of relations akin to the partitions we use in Section \ref{s:DC}; see \cite{DP2013}.  On the other hand, it is known that the minimal-size generating sets of the endomorphism monoids~$\B_n$ grow super-exponentially with $n$; see \cite[Corollary 3.1.8]{HMSW2020}.
\end{rem}

\subsection{Second main result: tensor category presentations}\label{ss:pres2}

We now show how to rewrite the presentation from Theorem \ref{t:pres} into a tensor presentation, again under certain natural assumptions, stated below.  For the duration of Section \ref{ss:pres2}, we fix a category~$\C$ (with object set $\S$, being one of $\N$ or $\PP$) satisfying Assumptions \ref{a:1}--\ref{a:4}.  We also keep all the notation of Section~\ref{ss:pres}, in particular the presentation $\pres\Ga\Om$ from Theorem~\ref{t:pres}, including the surmorphism ${\phi:\Ga^*\to\C:w\mt\ol w}$.

\begin{ass}\label{a:5}
We assume that $\C$ is a (strict) tensor category over $\S$.
\end{ass}

As in Section \ref{ss:cat}, we denote the tensor operation on $\C$ by $\op$.  

\begin{ass}\label{a:6}
We assume that
\bit
\item $\De$ is a digraph on vertex set $\S$, 
\item $\Xi\sub\De^\oast\times\De^\oast$ is a set of relations over $\De$, and 
\item $\Phi:\De^\oast\to\C$ is a morphism.  For $w\in\De^\oast$ we write $\ul w=w\Phi$.
\eit
Writing $\approx$ for the congruence $\Xi\tsharp$ on $\De^\oast$, we also assume that the following all hold:
\ben
\item \label{a61}  For every relation $(u,v)\in\Xi$, we have $\ul u=\ul v$.
\item \label{a62}  There is a mapping $\Ga\to\De^\oast:x\mt\wh x$ such that $\ulwh{x}=\ol x$ (i.e., $\wh x\Phi=x\phi$) for all $x\in\Ga$.  As $\Ga^*$ is freely generated by $\Ga$, we can extend this mapping to a morphism  $\Ga^*\to\De^\oast:w\mt\wh w$.  (It quickly follows that $\ulwh{w}=\ol w$ for all $w\in\Ga^*$.)
\item \label{a63}  For a generator $x\in\De$, and for natural numbers $m,n\in\N$, we define the term
\[
x_{m,n} = \io_m\op x\op\io_n\in\De^\oast.
\]
(Again, when $\S=\PP$, we interpret $x_{0,n}=x\op\io_n$, and so on.)  We assume that for every such $x,m,n$, we have $x_{m,n}\approx\wh w$ for some $w\in\Ga^*$.
\item \label{a64}  For every relation $(u,v)\in\Om$, we have $\wh u\approx\wh v$.
\een
\end{ass}

Here is our second main result:

\begin{thm}\label{t:pres2}
With notation as above, and subject to Assumptions \ref{a:1}--\ref{a:6}, the category $\C$ has tensor presentation $\pres\De\Xi$ via $\Phi$.
\end{thm}

\pf
To show that $\Phi$ is surjective, suppose $a\in\C$.  Since $\phi$ is surjective, we have $a=w\phi$ for some $w\in\Ga^*$.  But then $a=\wh w\Phi$ by Assumption \ref{a:6}\ref{a62}.

As in the proof of Theorem \ref{t:pres}, Assumption \ref{a:6}\ref{a61} gives $\Xi\tsharp\sub\ker(\Phi)$.  For the reverse inclusion, suppose $(u,v)\in\ker(\Phi)$, meaning that $u,v\in\De^\oast$ and $\ul u=\ul v$.  First writing $u$ and $v$ in the form \eqref{e:term2}, and then applying Assumption \ref{a:6}\ref{a63}, we have $u\approx\wh s$ and $v\approx\wh t$ for some $s,t\in\Ga^*$.  Assumption \ref{a:6}\ref{a61} and~\ref{a62} then give $\ol s=\ulwh{s}=\ul u=\ul v=\ulwh{t}=\ol t$.  It follows from Theorem \ref{t:pres} that there is a sequence
\[
s = w_0 \to w_1 \to \cdots \to w_k = t,
\]
where for each $1\leq i\leq k$, $w_i\in\Ga^*$ is obtained from $w_{i-1}$ by a single application of a relation from $\Om$.  It then follows from Assumption~\ref{a:6}\ref{a64} that $\wh s = \wh w_0 \approx \wh w_1 \approx\cdots\approx\wh w_k = \wh t$.  Thus, $u\approx \wh s\approx\wh t\approx v$.
\epf

In Sections \ref{s:DC} and \ref{s:BC}, we apply Theorems \ref{t:pres} and \ref{t:pres2} to a number of concrete categories $\C$ in order to rapidly obtain useful presentations.  Typically, the situation is as follows.
\bit
\item
In each case, Assumptions~\ref{a:1} and~\ref{a:2} are easily checked, and the presentations for the endomorphism monoids $\C_n$ required in Assumption~\ref{a:3} are imported from various sources in the literature \cite{HR2005,JEgrpm,JEgrpm2,KM2006,JE2021,BDP2002,JE2007a,JE2007b,Lavers1997,EL2004,Gilbert2006}.  Conditions \ref{a41}--\ref{a44} of Assumption \ref{a:4} will likewise be easy to check, with only condition~\ref{a45} a little less straightforward.  With these checks done, Theorem \ref{t:pres} then gives us a category presentation $\pres\Ga\Om$ for $\C$.
\item
We then use Theorem \ref{t:pres2} to rewrite $\pres\Ga\Om$ into a tensor presentation $\pres\De\Xi$.  The main work here involves checking items \ref{a62}--\ref{a64} of Assumption \ref{a:6}, and this largely consists of technical/elementary calculations.
\eit
The only exception to the above paradigm occurs in Section \ref{s:BC} when we consider the vine category $\V$, and the transformation categories $\T$ and $\O$.  If $\C$ is any of these, then we have $\C_{m,n}\not=\emptyset$ precisely when $m=0$ or $n\geq1$.  For such categories, we need a slight variation of Theorem~\ref{t:pres2}, and we address this in the next section.

\subsection{A variation on the second main result}\label{ss:pres3}

In this section we prove a result that will allow us to deal with the kinds of categories discussed at the end of Section \ref{ss:pres2}.  

\begin{ass}\label{a:7}
We assume $\C$ is a strict tensor category over $\N$ for which
\[
\C_{m,n} \not=\emptyset \iff m=0 \text{ or } n\geq1.
\]
We assume additionally that $|\C_{0,n}|=1$ for all $n\in\N$, and we write $\ul\U$ for the unique element of $\C_{0,1}$.  It follows that the unique element of $\C_{0,n}$ is $\ul\U^{\op n} = \ul\U\op\cdots\op\ul\U$, for each $n\in\N$.  (We interpret $\ul\U^{\op0}=\io_0$.)
\end{ass}

We define the subcategory $\C^+=\bigcup_{m,n\in\PP}\C_{m,n}$, which has object set $\PP$, and we assume that $\C^+$ satisfies Assumptions \ref{a:1}--\ref{a:4} (with $d=1$, as $\C_{0,1}\not=\emptyset$).  We fix the presentation $\pres\Ga\Om$ for $\C^+$ given in Theorem \ref{t:pres}.  

\begin{ass}\label{a:8}
We assume that
\bit
\item $\De$ is a digraph on vertex set $\N$,
\item $\Xi\sub\De^\oast\times\De^\oast$ is a set of relations over $\De$, and 
\item $\Phi:\De^\oast\to\C$ is a morphism.  For $w\in\De^\oast$ we write $\ul w=w\Phi$.
\eit
Writing $\approx$ for the congruence $\Xi\tsharp$ on $\De^\oast$, we also assume that the following all hold:
\ben
\item \label{a80}  The digraph $\De$ has a unique edge $\U$ with source $0$, and we also have $\br(\U)=1$.
\item \label{a81}  For every relation $(u,v)\in\Xi$, we have $\ul u=\ul v$.
\item \label{a82}  There is a mapping $\Ga\to\De^\oast:x\mt\wh x$ such that $\ulwh{x}=\ol x$ (i.e., $\wh x\Phi=x\phi$) for all $x\in\Ga$.  We extend this to a morphism $\Ga^*\to\De^\oast:w\mt\wh w$.  
\item \label{a83}  For a generator $x\in\De$, and for natural numbers $m,n\in\N$, define $x_{m,n} = \io_m\op x\op\io_n\in\De^\oast$.
We assume that for $(m,x,n)\not=(0,\U,0)$, we have $x_{m,n}\approx\wh w$ for some $w\in\Ga^*$.
\item \label{a84}  For every relation $(u,v)\in\Om$, we have $\wh u\approx\wh v$.
\item \label{a85}  For every $x\in\De\sm\{\U\}$, we have $\U^{\op m}\circ x \approx \U^{\op n}$, where $m=\bd(x)$ and $n=\br(x)$.
\een
\end{ass}

\begin{thm}\label{t:pres3}
With notation as above, and subject to $\C^+$ satisfying Assumptions \ref{a:1}--\ref{a:4}, and $\C$ satisfying Assumptions \ref{a:7} and \ref{a:8}, the category $\C$ has tensor presentation $\pres\De\Xi$ via $\Phi$.
\end{thm}

\pf
To show that $\Phi$ is surjective, suppose $a\in\C$.  If $\bd(a)\not=0$, then we follow the proof of Theorem~\ref{t:pres2} to show that $a=\wh w\Phi$ for some $w\in\Ga^*$.  Otherwise, $a=\ul\U^{\op n}=\U^{\op n}\Phi$ for some $n\in\N$.

Again it remains to show that $\ker(\Phi)\sub{\approx}$, where ${\approx}=\Xi\tsharp$, so suppose $(u,v)\in\ker(\Phi)$, and write $m=\bd(\ul u)=\bd(\ul v)$.  If $m\not=0$, then we follow the proof of Theorem \ref{t:pres2} to show that $u\approx v$.  For the rest of the proof we assume that $m=0$.  The proof will be complete if we can show that $u\approx\U^{\op n}$ for some $n\in\N$, as then $n=\br(\ul\U^{\op n})=\br(\ul u)=\br(\ul v)$, and the same argument will also give $v\approx\U^{\op n}$.  We may assume that $u=X_1\circ\cdots\circ X_k$ has the form \eqref{e:term2}, and we use induction on $k$.  If $k=0$, then $u=\io_0=\U^{\op 0}$, so we now assume that $k\geq1$.  By induction, $X_1\circ\cdots\circ X_{k-1}\approx \U^{\op q}$, where $q=\br(X_1\circ\cdots\circ X_{k-1})=\bd(X_k)$.  Further, we have $X_k = \io_a\op x\op\io_b$ for some $a,b\in\N$ and $x\in\De$.  If $x=\U$ then $q=\bd(X_k)=a+b$ and
\[
u \approx \U^{\op q} \circ(\io_a\op\U\op\io_b) = (\U^{\op a}\op\io_0\op\U^{\op b})\circ(\io_a\op\U\op\io_b) = (\U^{\op a}\circ\io_a)\op(\io_0\circ\U)\op(\U^{\op b}\circ\io_b) = \U^{\op(a+b+1)}.
\]
Otherwise, writing $c=\bd(x)$ and $d=\br(x)$, we have $q=a+c+b$, and Assumption~\ref{a:8}\ref{a85} gives
\begin{align*}
u &\approx (\U^{\op a}\op\U^{\op c}\op\U^{\op b}) \circ(\io_a\op x\op\io_b) \\
&= (\U^{\op a}\circ\io_a)\op(\U^{\op c}\circ x)\op(\U^{\op b}\circ\io_b) \approx \U^{\op a}\op\U^{\op d}\op\U^{\op b} = \U^{\op(a+d+b)}.  \qedhere
\end{align*}
\epf

\section{Diagram categories}\label{s:DC}

In this section we apply the theory developed in Section \ref{s:C} to a number of diagram categories, namely the partition category $\P$ (Section \ref{ss:P}), the Brauer category $\B$ (Section \ref{ss:B}) and the Temperley-Lieb category $\TL$ (Section \ref{ss:TL}).  In each case, we give a category presentation and a tensor presentation, based respectively on Theorems \ref{t:pres} and \ref{t:pres2}.  In Section \ref{ss:L} we show how these lead to presentations for linear versions of $\P$, $\B$ and $\TL$.  We begin by reviewing the relevant definitions and fixing notation.

\subsection{Preliminaries}\label{ss:preDC}

For $n\in\N$ we define the set $[n]=\{1,\ldots,n\}$, interpreting $[0]=\emptyset$.  For any subset $A$ of $\N$, we fix two disjoint copies of~$A$, namely $A'=\set{a'}{a\in A}$ and $A''=\set{a''}{a\in A}$.  The partition category $\P$ consists of all set partitions $\al$ of $[m]\cup[n]'$, for all $m,n\in\N$, under a composition defined shortly.  For such a partition~$\al$, we write $\bd(\al)=m$ and $\br(\al)=n$, and we write ${\P_{m,n}=\set{\al\in\P}{\bd(\al)=m,\ \br(\al)=n}}$.  A partition $\al\in\P_{m,n}$ will be identified with any graph on vertex set $[m]\cup[n]'$ whose connected components are the blocks of $\al$.  When drawing such a graph in the plane, a vertex $i\in[m]$ is always drawn at $(i,1)$, a vertex $i'\in[n]'$ at $(i,0)$, and all edges are contained in the rectangle ${\bigset{(x,y)\in\RRR^2}{1\leq x\leq \max(m,n),\ 0\leq y\leq1}}$.  See Figure \ref{f:P_product} for some examples.

To describe the composition operation on $\P$, let $m,n,q\in\N$ and fix some $\al\in\P_{m,n}$ and $\be\in\P_{n,q}$.  Let $\al_\downarrow$ be the graph on vertex set $[m]\cup[n]''$ obtained by renaming each lower vertex~$i'$ of $\al$ to $i''$, let~$\be^\uparrow$ be the graph on vertex set $[n]''\cup[q]'$ obtained by renaming each upper vertex $i$ of $\be$ to $i''$, and let $\Pi(\al,\be)$ be the graph on vertex set $[m]\cup[n]''\cup[q]'$ whose edge set is the union of the edge sets of~$\al_\downarrow$ and $\be^\uparrow$.  We call $\Pi(\al,\be)$ the product graph.  The product/composition $\al\be=\al\circ\be\in\P_{m,q}$ is then defined to be the partition of $[m]\cup[q]'$ for which $x,y\in[m]\cup[q]'$ belong to the same block of~$\al\be$ if and only if $x,y$ are connected by a path in $\Pi(\al,\be)$.  When drawing a product graph $\Pi(\al,\be)$, we draw the vertices from $[m]$, $[n]''$ and $[q]'$ at heights $y=2$, $y=1$ and $y=0$, respectively.  An example calculation is given in Figure \ref{f:P_product}.

Note that the product graph $\Pi(\al,\be)$ may contain ``floating components'', which are contained entirely in the middle row of the graph.  These do not figure in the composition on $\P$, but we will return our attention to them in Section \ref{ss:L} when we treat the linear diagram categories.

\begin{figure}[ht]
\begin{center}
\begin{tikzpicture}[scale=.45]
\begin{scope}[shift={(0,0)}]	
\uvs{1,...,6}
\lvs{1,...,8}
\uarcx14{.6}
\uarcx23{.3}
\uarcx56{.3}
\darc12
\darcx26{.6}
\darcx45{.3}
\darc78
\stline34
\draw(0.6,1)node[left]{$\al=$};
\draw[->](9.5,-1)--(11.5,-1);
\end{scope}
\begin{scope}[shift={(0,-4)}]	
\uvs{1,...,8}
\lvs{1,...,7}
\uarc12
\uarc34
\darc45
\darc67
\stline31
\stline55
\stline87
\draw(0.6,1)node[left]{$\be=$};
\end{scope}
\begin{scope}[shift={(12,-1)}]	
\uvs{1,...,6}
\lvs{1,...,8}
\uarcx14{.6}
\uarcx23{.3}
\uarcx56{.3}
\darc12
\darcx26{.6}
\darcx45{.3}
\darc78
\stline34
\draw[->](9.5,0)--(11.5,0);
\end{scope}
\begin{scope}[shift={(12,-3)}]	
\uvs{1,...,8}
\lvs{1,...,7}
\uarc12
\uarc34
\darc45
\darc67
\stline31
\stline55
\stline87
\end{scope}
\begin{scope}[shift={(24,-2)}]	
\uvs{1,...,6}
\lvs{1,...,7}
\uarcx14{.6}
\uarcx23{.3}
\uarcx56{.3}
\darc14
\darc45
\darc67
\stline21
\draw(7.4,1)node[right]{$=\al\be$};
\end{scope}
\end{tikzpicture}
\caption{Calculating $\al\be=\al\circ\be$, where $\al\in\P_{6,8}$ and $\be\in\P_{8,7}$.  The product graph $\Pi(\al,\be)$ is shown in the middle.}
\label{f:P_product}
\end{center}
\end{figure}

The identity of $\P$ at object $n$ is the partition $\iob_n=\bigset{\{i,i'\}}{i\in[n]}$; see Figure \ref{f:P_gens}.  The endomorphism monoid $\P_n=\P_{n,n}$ is the partition monoid of degree~$n$~\cite{HR2005}.  The units of $\P_n$ are the partitions of the form $\bigset{\{i,(i\pi)'\}}{i\in[n]}$ for some permutation~$\pi\in\SS_n$.  Such a unit will be identified with the permutation~$\pi$ itself, and the group of all such units with the symmetric group $\SS_n$.  

The category $\P$ has the structure of a (strict) tensor category, with the $\op$ operation defined as follows.  Let $m,n,k,l\in\N$, and let $\al\in\P_{m,n}$ and $\be\in\P_{k,l}$.  First, the partition $\al\op\be\in\P_{m+k,n+l}$ contains all the blocks of $\al$.  Additionally, for each block $A\cup B'$ of $\be$, $\al\op\be$ also contains the block $(A+m)\cup(B+n)'$.  Geometrically, $\al\op\be$ is obtained by placing a copy of $\be$ to the right of $\al$, as in Figure \ref{f:P_sum}.  The $\op$ operation, and the subgroups $\SS_n\sub\P_n$ ($n\in\N$), give $\P$ the structure of a PROP, as in \cite[Section 24]{MacLane1965}.

\begin{figure}[ht]
\begin{center}
\begin{tikzpicture}[scale=.45]
\begin{scope}[shift={(0,0)}]	
\uvcs{1,...,6}{red}
\lvcs{1,...,8}{red}
\uarcxc14{.6}{red}
\uarcxc23{.3}{red}
\darcxc23{.3}{red}
\darcxc16{.6}{red}
\darcxc45{.3}{red}
\stlinec57{red}
\stlinec68{red}
\draw(0.6,1)node[left]{${\red\al}=$};
\draw[->](9.5,-1)--(11.5,-1);
\end{scope}
\begin{scope}[shift={(0,-4)}]	
\uvcs{1,...,8}{blue}
\lvcs{1,...,7}{blue}
\uarcc12{blue}
\uarcc34{blue}
\darcc45{blue}
\darcc67{blue}
\stlinec31{blue}
\stlinec55{blue}
\stlinec87{blue}
\draw(0.6,1)node[left]{${\blue\be}=$};
\end{scope}
\begin{scope}[shift={(12,-2)}]	
\uvcs{1,...,6}{red}
\lvcs{1,...,8}{red}
\uarcxc14{.6}{red}
\uarcxc23{.3}{red}
\darcxc23{.3}{red}
\darcxc16{.6}{red}
\darcxc45{.3}{red}
\stlinec57{red}
\stlinec68{red}
\uvcs{7,...,14}{blue}
\lvcs{9,...,15}{blue}
\uarcc78{blue}
\uarcc9{10}{blue}
\darcxc{12}{13}{.2}{blue}
\darcxc{14}{15}{.2}{blue}
\stlinec99{blue}
\stlinec{11}{13}{blue}
\stlinec{14}{15}{blue}
\draw(15.4,1)node[right]{$={\red\al}\op{\blue\be}$};
\end{scope}
\end{tikzpicture}
\caption{Calculating $\al\op\be$, where $\al\in\P_{6,8}$ and $\be\in\P_{8,7}$.}
\label{f:P_sum}
\end{center}
\end{figure}

The category $\P$ has a natural involution $\al\mt\al^*$ obtained by interchanging dashed and un-dashed elements.  Geometrically, $\al^*$ is obtained by reflecting (a graph representing) $\al$ in a horizontal axis.  This gives $\P$ the structure of a so-called regular $*$-category, as defined in \cite[Section 2]{DDE2021}, meaning that the following hold for all $\al,\be\in\P$ with $\br(\al)=\bd(\be)$:
\[
\bd(\al^*)=\br(\al) \COMMA \br(\al^*)=\bd(\al) \COMMA (\al^*)^*=\al \COMMA \al=\al\al^*\al \COMMA (\al\be)^*=\be^*\al^*.
\]
We also have $(\al\op\be)^*=\al^*\op\be^*$ for all~$\al,\be$.

A partition $\al\in\P$ is planar if some graph representing it may be drawn (in the plane, as described above) with no edge crossings.  In Figure \ref{f:P_product}, for example, $\be$ is planar but~$\al$ is not.  The set $\PlP$ of all planar partitions is a subcategory of $\P$, the planar partition category.

A partition $\al\in\P$ is called a Brauer partition if each block has size $2$.  The set $\B$ of all such Brauer partitions is the Brauer category.  Clearly $\B_{m,n}$ is non-empty if and only if $m\equiv n\Mod2$.  The set $\TL=\PlP\cap\B$ of all planar Brauer partitions is the Temperley-Lieb category.  The partition $\al$ pictured in Figure \ref{f:P_sum} belongs to $\TL$.

The categories $\PlP$, $\B$ and $\TL$ are also closed under $\op$ and ${}^*$, so these are all (strict) tensor regular $*$-categories.  The symmetric groups $\SS_n$ ($n\in\N$) are contained in $\B$ but not in $\PlP$ or~$\TL$.  It follows that $\B$ is a PROP, though $\PlP$ and $\TL$ are not.

\subsection{The partition category}\label{ss:P}

We now come to the first of our applications of the general machinery developed in Section \ref{s:C}.  Our goal in this section is to apply Theorems \ref{t:pres} and \ref{t:pres2} to obtain presentations for the partition category~$\P$; see Theorems \ref{t:P1} and \ref{t:P2}, below.  This section can be thought of as a blueprint for those that follow, so our treatment will be fairly detailed.

First note that Assumption \ref{a:1} holds in $\P$, with $d=1$, as all hom-sets $\P_{m,n}$ are non-empty.  For Assumption \ref{a:2}, we take the partitions $\lamb_n\in\P_{n,n+1}$ and $\rhob_n\in\P_{n+1,n}$ shown in Figure \ref{f:P_gens}; clearly $\lamb_n\rhob_n=\iob_n$ for all $n\in\N$.  For Assumption~\ref{a:3}, we require presentations for the partition monoids $\P_n$ ($n\in\N$).  Such presentations are stated in \cite{HR2005}; proofs may be found in \cite{JEgrpm,JEgrpm2}.  First, for $n\in\N$ define an alphabet
\begin{align*}
X_n = S_n\cup E_n\cup T_n \qquad\text{where}\qquad S_n &= \set{\si_{i;n}}{1\leq i<n} ,\\
E_n &= \set{\ve_{i;n}}{1\leq i\leq n} , \\
T_n &= \set{\tau_{i;n}}{1\leq i<n}.
\end{align*}
Define the morphism $\phi_n:X_n^*\to\P_n:w\mt\ol w$, where the partitions $\ol x$ ($x\in X_n$) are shown in Figure~\ref{f:P_gens}.

\begin{figure}[ht]
\begin{center}
\begin{tikzpicture}[scale=.45]
\begin{scope}[shift={(16,-5)}]	
\uvs{1,6}
\lvs{1,6,7}
\udotted16
\ddotted16
\stline11
\stline66
\draw(0.5,1)node[left]{$\lamb_n=$};
\node()at(1,2.5){\tiny$1$};
\node()at(6,2.5){\tiny$n$};
\end{scope}
\begin{scope}[shift={(16,-10)}]	
\uvs{1,6,7}
\lvs{1,6}
\udotted16
\ddotted16
\stline11
\stline66
\draw(0.5,1)node[left]{$\rhob_n=$};
\node()at(1,2.5){\tiny$1$};
\node()at(6,2.5){\tiny$n$};
\end{scope}
\begin{scope}[shift={(16,0)}]	
\uvs{1,6}
\lvs{1,6}
\udotted16
\ddotted16
\stline11
\stline66
\draw(0.5,1)node[left]{$\iob_n=$};
\node()at(1,2.5){\tiny$1$};
\node()at(6,2.5){\tiny$n$};
\end{scope}
\begin{scope}[shift={(0,0)}]	
\uvs{1,4,5,6,7,10}
\lvs{1,4,5,6,7,10}
\udotted14
\ddotted14
\udotted7{10}
\ddotted7{10}
\stline11
\stline44
\stline56
\stline65
\stline77
\stline{10}{10}
\draw(0.5,1)node[left]{$\sib_{i;n}=$};
\node()at(1,2.5){\tiny$1$};
\node()at(5,2.5){\tiny$i$};
\node()at(10,2.5){\tiny$n$};
\end{scope}
\begin{scope}[shift={(0,-10)}]	
\uvs{1,4,5,6,7,10}
\lvs{1,4,5,6,7,10}
\udotted14
\ddotted14
\udotted7{10}
\ddotted7{10}
\stline11
\stline44
\stline55
\stline66
\stline77
\stline{10}{10}
\uarcx56{.3}
\darcx56{.3}
\draw(0.5,1)node[left]{$\taub_{i;n}=$};
\node()at(1,2.5){\tiny$1$};
\node()at(5,2.5){\tiny$i$};
\node()at(10,2.5){\tiny$n$};
\end{scope}
\begin{scope}[shift={(0,-5)}]	
\uvs{1,4,5,6,10}
\lvs{1,4,5,6,10}
\udotted14
\ddotted14
\udotted6{10}
\ddotted6{10}
\stline11
\stline44
\stline66
\stline{10}{10}
\draw(0.5,1)node[left]{$\veb_{i;n}=$};
\node()at(1,2.5){\tiny$1$};
\node()at(5,2.5){\tiny$i$};
\node()at(10,2.5){\tiny$n$};
\end{scope}
\end{tikzpicture}
\caption{Partition generators $\ol x\in\P$ ($x\in\Ga$), as well as $\iob_n$.}
\label{f:P_gens}
\end{center}
\end{figure}

Now let $R_n$ be the following set of relations over $X_n$ (here, and in other such lists of relations, the subscripts range over all meaningful values, subject to any stated constraints):
\begin{align}
\label{P1}\tag*{(P1)}  \si_{i;n}^2 = \iota_n , \qquad\ \ve_{i;n}^2 = \ve_{i;n} , &&&  \tau_{i;n}^2 = \tau_{i;n} = \tau_{i;n}\si_{i;n} = \si_{i;n}\tau_{i;n},  \\
\label{P2}\tag*{(P2)}  \si_{i;n}\ve_{i;n} = \ve_{i+1;n}\si_{i;n} , &&& \ve_{i;n}\ve_{i+1;n}\si_{i;n} = \ve_{i;n}\ve_{i+1;n}, \\
\label{P3}\tag*{(P3)}  \ve_{i;n}\ve_{j;n} = \ve_{j;n}\ve_{i;n}  , &&& \tau_{i;n}\tau_{j;n} = \tau_{j;n}\tau_{i;n},  \\
\label{P4}\tag*{(P4)}  \si_{i;n}\si_{j;n} = \si_{j;n}\si_{i;n} , &&& \si_{i;n}\tau_{j;n} = \tau_{j;n}\si_{i;n},  &&\text{if $|i-j|>1$,}   \\
\label{P5}\tag*{(P5)}  \si_{i;n}\si_{j;n}\si_{i;n} = \si_{j;n}\si_{i;n}\si_{j;n} , &&& \si_{i;n}\tau_{j;n}\si_{i;n} = \si_{j;n}\tau_{i;n}\si_{j;n},  &&\text{if $|i-j|=1$},   \\
\label{P6}\tag*{(P6)}  \si_{i;n}\ve_{j;n} = \ve_{j;n}\si_{i;n} , &&& \tau_{i;n}\ve_{j;n} = \ve_{j;n}\tau_{i;n}, &&\text{if $j\not=i,i+1$},   \\
\label{P7}\tag*{(P7)}  \tau_{i;n}\ve_{j;n}\tau_{i;n} = \tau_{i;n} , &&& \ve_{j;n}\tau_{i;n}\ve_{j;n} = \ve_{j;n}, &&\text{if $j=i,i+1$.}  
\end{align}

\begin{thm}[cf.~\cite{HR2005,JEgrpm,JEgrpm2}]\label{t:Pn}
For any $n\in\N$, the partition monoid $\P_n$ has presentation $\pres{X_n}{R_n}$ via $\phi_n$.  \epfres
\end{thm}

Now let $\Ga\equiv L\cup R\cup X$ be the digraph over vertex set $\N$ as in \eqref{e:Ga} and \eqref{e:st}, and let $\phi:\Ga^*\to\P$ be the morphism given in \eqref{e:phi}.  Let $\Om$ be the set of relations over $\Ga$ consisting of $\bigcup_{n\in\N}R_n$ and additionally:
\begin{align}
\label{P8}\tag*{(P8)}  \lam_n\rho_n = \iota_n , &&& \rho_n\lam_n = \ve_{n+1;n+1}, \\
\label{P9}\tag*{(P9)}  \th_{i;n}\lam_n = \lam_n\th_{i;n+1} , &&& \rho_n\th_{i;n} = \th_{i;n+1}\rho_n, &&\text{for $\th\in\{\si,\ve,\tau\}$.}  
\end{align}
In the language of Assumption \ref{a:4}\ref{a42}, we have $w_n=\ve_{n+1;n+1}$ for all $n\in\N$.  The maps in Assumption~\ref{a:4}\ref{a44} are given by $\th_{i;n}^+ = (\th_{i;n})_+ = \th_{i;n+1}$, for $\th\in\{\si,\ve,\tau\}$.

Here is the first main result of this section, expressed in terms of the above notation.  For the proof, it will be convenient to define an embedding
\[
\P_n\to\P_{n+1} : \al\mt\al^+ = \al\op\iob_1 \qquad\text{for each $n\in\N$.}
\]
So $\al^+\in\P_{n+1}$ is obtained from $\al\in\P_n$ by adding the transversal $\{n+1,(n+1)'\}$.  Our re-use of the~${}^+$ notation stems from the fact that $\ol x^+ = \ol{x^+}$ for all $x\in X_n$, and hence $\ol w^+=\ol{w^+}$ for all $w\in X_n^*$.

\begin{thm}\label{t:P1}
The partition category $\P$ has presentation~$\pres\Ga\Om$ via~$\phi$.  
\end{thm}

\pf
In order to apply Theorem \ref{t:pres}, all that remains to check is that items \ref{a41} and \ref{a45} of Assumption \ref{a:4} hold.  The first is easily checked diagrammatically.  In fact, given Theorem \ref{t:Pn}, we only need to check \ref{P8} and \ref{P9}; see Figure \ref{f:rel_P} for (part of) the latter.  To verify Assumption \ref{a:4}\ref{a45}, we use Remark~\ref{r:a45}, which says it suffices to show the following, where for simplicity we write $\ve=\ve_{n+1;n+1}$:
\bit
\item For all $n\in\N$ and  $w\in X_{n+1}^*$,  we have  $\ve w\ve  \sim \ve u^+\ve $  for some $u\in X_n^*$.
\eit
To show this, let $w\in X_{n+1}^*$, and put $\al = \ol{\ve w\ve }\in\P_{n+1}$.  Note that $\al$ contains the blocks $\{n+1\}$ and $\{n+1\}'$.  Let $\be\in\P_n$ be the partition obtained from $\al$ by deleting these two blocks; so~$\al=\veb\be^+\veb$.  By Theorem \ref{t:Pn}, we have $\be=\ol u$ for some $u\in X_n^*$, and we note that $\be^+=\ol u^+=\ol{u^+}$.  But then $\ol{\ve w\ve }=\al=\veb\be^+\veb = \ol{\ve u^+\ve }$.  It follows from Theorem~\ref{t:Pn} (and $R_{n+1}\sub\Om$) that $\ve w\ve \sim\ve u^+\ve $.
\epf

\begin{figure}[ht]
\begin{center}
\begin{tikzpicture}[scale=.43]
\begin{scope}[shift={(0,0)}]	
\uvs{1,4,5,6,7,10}
\lvs{1,4,5,6,7,10}
\udotted14
\ddotted14
\udotted7{10}
\ddotted7{10}
\stline11
\stline44
\stline56
\stline65
\stline77
\stline{10}{10}
\node()at(1,2.5){\tiny$1$};
\node()at(5,2.5){\tiny$i$};
\node()at(10,2.5){\tiny$n$};
\end{scope}
\begin{scope}[shift={(0,-2)}]	
\uvs{1,4,5,6,7,10}
\lvs{1,4,5,6,7,10,11}
\udotted14
\ddotted14
\udotted7{10}
\ddotted7{10}
\stline11
\stline44
\stline55
\stline66
\stline77
\stline{10}{10}
\end{scope}
\begin{scope}[shift={(14,-1)}]	
\uvs{1,4,5,6,7,10}
\lvs{1,4,5,6,7,10,11}
\udotted14
\ddotted14
\udotted7{10}
\ddotted7{10}
\stline11
\stline44
\stline56
\stline65
\stline77
\stline{10}{10}
\node()at(1,2.5){\tiny$1$};
\node()at(5,2.5){\tiny$i$};
\node()at(10,2.5){\tiny$n$};
\node()at(-1,1){$=$};
\node()at(13,1){$=$};
\end{scope}
\begin{scope}[shift={(28,0)}]	
\uvs{1,4,5,6,7,10}
\lvs{1,4,5,6,7,10}
\udotted14
\ddotted14
\udotted7{10}
\ddotted7{10}
\stline11
\stline44
\stline55
\stline66
\stline77
\stline{10}{10}
\node()at(1,2.5){\tiny$1$};
\node()at(5,2.5){\tiny$i$};
\node()at(10,2.5){\tiny$n$};
\end{scope}
\begin{scope}[shift={(28,-2)}]	
\uvs{1,4,5,6,7,10,11}
\lvs{1,4,5,6,7,10,11}
\udotted14
\ddotted14
\udotted7{10}
\ddotted7{10}
\stline11
\stline44
\stline56
\stline65
\stline77
\stline{10}{10}
\stline{11}{11}
\end{scope}
\end{tikzpicture}
\caption{Part of relation \ref{P9}: $\sib_{i;n}\lamb_n=\lamb_n\sib_{i;n+1}$.}
\label{f:rel_P}
\end{center}
\end{figure}

\begin{rem}\label{r:P1}
The proof of Theorem \ref{t:P1} involved an application of Theorem \ref{t:pres}.  Various ingredients in the proof of Theorem \ref{t:pres} have natural diagrammatic meanings when $\C=\P$.  For example, Figure~\ref{f:lrP} shows the elements $\lamb_{m,n}$ and $\rhob_{n,m}$ ($0\leq m\leq n$), and the mappings $\R_{m,n}:\P_{m,n}\to\P_n$ and ${\L_{n,m}:\P_{n,m}\to\P_n}$ defined in \eqref{e:RL}.
\end{rem}

\begin{figure}[ht]
\begin{center}
\begin{tikzpicture}[scale=.43]
\begin{scope}[shift={(-14,6)}]	
\uvs{1,4}
\lvs{1,4,5,8}
\stline11
\stline44
\udotted14
\ddotted14
\ddotted58
\draw(0.5,1)node[left]{$\lamb_{m,n}=$};
\node()at(1,2.5){\tiny$1$};
\node()at(4,2.5){\tiny$m$};
\node()at(1,-.5){\tiny$1$};
\node()at(8,-.5){\tiny$n$};
\end{scope}
\begin{scope}[shift={(-14,0)}]	
\uvs{1,4,5,8}
\lvs{1,4}
\stline11
\stline44
\udotted14
\udotted58
\ddotted14
\draw(0.5,1)node[left]{$\rhob_{n,m}=$};
\node()at(1,2.5){\tiny$1$};
\node()at(8,2.5){\tiny$n$};
\node()at(1,-.5){\tiny$1$};
\node()at(4,-.5){\tiny$m$};
\end{scope}
\begin{scope}[shift={(0,6)}]	
\fill[blue!20](1,0)--(8,0)--(4,2)--(1,2)--(1,0);
\uvs{1,4}
\lvs{1,8}
\udotted14
\ddotted18
\draw(0.5,1)node[left]{$\al=$};
\node()at(1,2.5){\tiny$1$};
\node()at(4,2.5){\tiny$m$};
\node()at(1,-.5){\tiny$1$};
\node()at(8,-.5){\tiny$n$};
\draw[|-{latex}](4.5,-1.5)--(4.5,-2.5);
\end{scope}
\begin{scope}[shift={(0,0)}]	
\fill[blue!20](1,0)--(8,0)--(4,2)--(1,2)--(1,0);
\uvs{1,4,5,8}
\lvs{1,8}
\udotted14
\udotted58
\ddotted18
\draw(0.5,1)node[left]{\small $\R_{m,n}(\al)=$};
\node()at(1,2.5){\tiny$1$};
\node()at(4,2.5){\tiny$m$};
\node()at(8,2.5){\tiny$n$};
\node()at(1,-.5){\tiny$1$};
\node()at(8,-.5){\tiny$n$};
\end{scope}
\begin{scope}[shift={(14,6)}]	
\fill[blue!20](1,0)--(4,0)--(8,2)--(1,2)--(1,0);
\uvs{1,8}
\lvs{1,4}
\udotted18
\ddotted14
\draw(0.5,1)node[left]{$\al=$};
\node()at(1,2.5){\tiny$1$};
\node()at(8,2.5){\tiny$n$};
\node()at(1,-.5){\tiny$1$};
\node()at(4,-.5){\tiny$m$};
\draw[|-{latex}](4.5,-1.5)--(4.5,-2.5);
\end{scope}
\begin{scope}[shift={(14,0)}]	
\fill[blue!20](1,0)--(4,0)--(8,2)--(1,2)--(1,0);
\uvs{1,8}
\lvs{1,4,5,8}
\udotted18
\ddotted14
\ddotted58
\draw(0.5,1)node[left]{\small $\L_{n,m}(\al)=$};
\node()at(1,2.5){\tiny$1$};
\node()at(8,2.5){\tiny$n$};
\node()at(4,-.5){\tiny$m$};
\node()at(1,-.5){\tiny$1$};
\node()at(8,-.5){\tiny$n$};
\end{scope}
\end{tikzpicture}
\caption{The partitions $\lamb_{m,n}$ and $\rhob_{n,m}$ ($0\leq m\leq n$), and the mappings $\R_{m,n}:\P_{m,n}\to\P_n$ and $\L_{n,m}:\P_{n,m}\to\P_n$; cf.~Remark \ref{r:P1}.}
\label{f:lrP}
\end{center}
\end{figure}

\begin{rem}\label{r:i+1}
Certain redundancies exist in the presentation $\pres\Ga\Om$ from Theorem \ref{t:P1}.  For example, only the $j=i+1$ case of \ref{P7} is needed.  Indeed, using this and \ref{P1} and \ref{P2}, and writing $\th_i=\th_{i;n}$, we have $\tau_i\ve_i\tau_i \sim \tau_i\si_i\ve_i\tau_i \sim \tau_i\ve_{i+1}\si_i\tau_i \sim \tau_i\ve_{i+1}\tau_i \sim \tau_i$.  The other part of \ref{P7} is similar.
\end{rem}

\begin{rem}\label{r:a4'}
As in Remark \ref{r:a4}, we could have made different choices for the one-sided units ${\lamb_n,\rhob_n\in\P}$ to those shown in Figure \ref{f:P_gens}.  For example, out of many other possibilities, we could instead have taken $\lamb_n$ and $\rhob_n$ to be of the respective forms:
\[
\begin{tikzpicture}[scale=.45]
\uvs{1,6,7}
\lvs{1,6,7,8}
\udotted16
\ddotted16
\stline11
\stline66
\stline77
\stline78
\draw(0.5,1)node[left]{$\lamb_n=$};
\node()at(1,2.5){\tiny$1$};
\node()at(7,2.5){\tiny$n$};
\end{tikzpicture}
\begin{tikzpicture}[scale=.45]
\uvs{1,6,7,8}
\lvs{1,6,7}
\udotted16
\ddotted16
\stline11
\stline66
\stline77
\stline87
\draw(0.5,1)node[left]{\AND$\rhob_n=$};
\node()at(1,2.5){\tiny$1$};
\node()at(7,2.5){\tiny$n$};
%\draw(9,0)node[left]{.};
%
\end{tikzpicture}
\]
This would result in a change to the second part of relation \ref{P8}, which would become $\rho_n\lam_n = \tau_{n;n+1}$.  The first part of \ref{P8} would remain the same as above, as would both parts of \ref{P9}.  The resulting presentation for $\P$ does not seem to be any less ``natural'' than that given in Theorem \ref{t:P1}, though of course the above proof of the theorem would have to be modified accordingly.
\end{rem}

Now that we have the presentation $\pres\Ga\Om$ for $\P$, we wish to use Theorem \ref{t:pres2} to transform it into a tensor presentation.  We begin by defining the digraph $\De$ on vertex set $\N$ with four edges: 
\[
X:2\to2 \COMMA D:2\to2 \COMMA U:1\to0 \COMMA \U:0\to1.
\]
Define the morphism $\Phi:\De^\oast\to\P:w\mt\ul w$, where the partitions $\ul x$ ($x\in\De$) are shown in Figure \ref{f:P_Tgens}.  It will also be convenient to write $I=\io_1$ for the empty path $1\to1$, and $\ul I=\iob_1$ for the identity partition from $\P_1$; the latter is also shown in Figure \ref{f:P_Tgens}.  

\begin{figure}[ht]
\begin{center}
\begin{tikzpicture}[scale=.45]
\begin{scope}[shift={(14,0)}]	
\uv1
\draw(0.5,1)node[left]{$\ul U=$};
\end{scope}
\begin{scope}[shift={(20,0)}]	
\lv1
\draw(0.5,1)node[left]{$\ul \U=$};
\end{scope}
\begin{scope}[shift={(0,0)}]	
\lv1
\lv2
\uv1
\uv2
\stline12
\stline21
\draw(0.5,1)node[left]{$\ul X=$};
\end{scope}
\begin{scope}[shift={(7,0)}]	
\lv1
\lv2
\uv1
\uv2
\stline11
\stline22
\darc12
\uarc12
\draw(0.5,1)node[left]{$\ul D=$};
\end{scope}
\begin{scope}[shift={(26,0)}]	
\uv1
\lv1
\stline11
\draw(0.5,1)node[left]{$\ul I=$};
\end{scope}
\end{tikzpicture}
\caption{Partition generators $\ul x\in\P$ ($x\in\De$), as well as $\ul I$.}
\label{f:P_Tgens}
\end{center}
\end{figure}

Let $\Xi$ be the set of the following relations over $\De$, remembering that~$I=\io_1$:
\begin{gather}
\tag*{(P1)$'$} \label{P1'} X\circ X=\io_2 \COMMA \U\circ U = \io_0  , \\
\tag*{(P2)$'$} \label{P2'} D\circ D=D = D\circ X=X\circ D \COMMA (D\op I)\circ(I\op D)=(I\op D)\circ(D\op I),\\
\tag*{(P3)$'$} \label{P3'} (X\op I)\circ(I\op X)\circ(X\op I)=(I\op X)\circ(X\op I)\circ(I\op X),\\
\tag*{(P4)$'$} \label{P4'} (X\op I)\circ(I\op D)\circ(X\op I)=(I\op X)\circ(D\op I)\circ(I\op X),\\
\tag*{(P5)$'$} \label{P5'} X\circ(I\op U)=U\op I \COMMA  (I\op\U)\circ X=\U\op I,\\
\tag*{(P6)$'$} \label{P6'} (I\op \U)\circ D\circ(I\op U) = I \COMMA  D\circ(I\op U\op\U)\circ D=D.
\end{gather}
Note that by Lemma \ref{l:ab} we automatically have the additional relations
\begin{equation}\label{e:UU}
U\circ\U = U\op\U = \U\op U.
\end{equation}
Several other tensor presentations in the paper involve generators with $\br(U)=\bd(\U)=0$, so \eqref{e:UU} holds in all of those as well.

The second main result of this section is the following, expressed in terms of the above notation:

\begin{thm}\label{t:P2}
The partition category $\P$ has tensor presentation $\pres\De\Xi$ via $\Phi$.
\end{thm}

\pf
To apply Theorem \ref{t:pres2}, it remains to verify the conditions of Assumption \ref{a:6}.  Condition~\ref{a61} is easy to check diagrammatically; see Figure~\ref{f:P5'} for~\ref{P4'}.  

\begin{figure}[ht]
\begin{center}
\begin{tikzpicture}[scale=.45]
\begin{scope}[shift={(0,4)}]	
\uvs{1,2,3}
\lvs{1,2,3}
\stline12
\stline21
\stline33
\end{scope}
\begin{scope}[shift={(0,2)}]	
\uvs{1,2,3}
\lvs{1,2,3}
\stline11
\stline22
\stline33
\uarc23
\darc23
\end{scope}
\begin{scope}[shift={(0,0)}]	
\uvs{1,2,3}
\lvs{1,2,3}
\stline12
\stline21
\stline33
\end{scope}
\begin{scope}[shift={(7,2)}]	
\uvs{1,2,3}
\lvs{1,2,3}
\stline11
\stline22
\stline33
\uarc13
\darc13
\draw(-1.5,1)node{$=$};
\draw(5.5,1)node{$=$};
\end{scope}
\begin{scope}[shift={(14,4)}]	
\uvs{1,2,3}
\lvs{1,2,3}
\stline11
\stline23
\stline32
\end{scope}
\begin{scope}[shift={(14,2)}]	
\uvs{1,2,3}
\lvs{1,2,3}
\stline11
\stline22
\stline33
\uarc12
\darc12
\end{scope}
\begin{scope}[shift={(14,0)}]	
\uvs{1,2,3}
\lvs{1,2,3}
\stline11
\stline23
\stline32
\end{scope}
\end{tikzpicture}
\caption{Relation \ref{P4'}: $(\ul X\op \ul I)\circ(\ul I\op \ul D)\circ(\ul X\op \ul I)=(\ul I\op \ul X)\circ(\ul D\op \ul I)\circ(\ul I\op \ul X)$.}
\label{f:P5'}
\end{center}
\end{figure}

Conditions \ref{a62}--\ref{a64} involve a morphism $\Ga^*\to\De^\oast:w\mt\wh w$, which we define by
\begin{align*}
\wh\si_{i;n} &= \io_{i-1}\op X\op\io_{n-i-1}, & \wh\ve_{i;n} &= \io_{i-1}\op U\op\U\op\io_{n-i}, & \wh\lam_n &= \io_n\op\U,\\
\wh\tau_{i;n} &= \io_{i-1}\op D\op\io_{n-i-1}, &&& \wh\rho_n &= \io_n\op U.
\end{align*}
Condition \ref{a62} says that $\wh x\Phi=x\phi$ for all $x\in\Ga$, and this is easily verified diagrammatically.  Conditions~\ref{a63} and~\ref{a64} are verified in Lemmas \ref{l:a63P} and \ref{l:a64P}, respectively.  The first of these refers to the terms $x_{m,n} = \io_m\op x\op\io_n$ ($x\in\De$, $m,n\in\N$).  For the rest of the proof we write ${\approx}=\Xi\tsharp$.

\begin{lemma}\label{l:a63P}
For any $m,n\in\N$, and for any $x\in\De$, we have $x_{m,n} \approx\wh w$ for some $w\in\Ga^*$.
\end{lemma}

\pf
This is clear for $x=X$ or $D$, since $X_{m,n} = \wh\si_{m+1;m+n+2}$ and $D_{m,n} = \wh\tau_{m+1;m+n+2}$.  By symmetry, it remains to consider the case of $x=U$.  For this we claim that 
\begin{equation}\label{e:Umn}
U_{m,n} \approx \wh\si_{m+1;m+n+1} \circ\cdots\circ \wh\si_{m+n;m+n+1} \circ\wh\rho_{m+n}.
\end{equation}
We prove this (for any $m$) by induction on $n$.  For $n=0$ we have $U_{m,n} = \io_m \op U \op \io_0 = \io_m\op U = \wh\rho_m$, which agrees with \eqref{e:Umn} since the product of $\wh\si$'s is empty when $n=0$.  For $n\geq1$, \ref{P5'} gives
\begin{align*}
\wh\si_{m+1;m+n+1}\circ U_{m+1;n-1} &= (\io_m\op X\op\io_{n-1}) \circ (\io_m\op I\op U\op\io_{n-1}) \\
&= \io_m\op(X\circ(I\op U))\op\io_{n-1} \approx \io_m\op U\op I\op\io_{n-1} = U_{m,n},
\end{align*}
and we then apply the inductive assumption to $U_{m+1;n-1}$.
\epf

\begin{lemma}\label{l:a64P}
For any relation $(u,v)\in\Om$, we have $\wh u\approx\wh v$.
\end{lemma}

\pf
Most of the commuting relations from $\Om$ follow immediately from the tensor category axioms.  For example, consider the second part of \ref{P4}.  If $j\geq i+2$, then
\begin{align*}
\wh\si_{i;n}\circ\wh\tau_{j;n} &= (\io_{i-1}\op X\op\io_{j-i-2}\op\io_2\op\io_{n-j-1} ) \circ ( \io_{i-1}\op \io_2\op\io_{j-i-2}\op D\op\io_{n-j-1}) \\
& = \io_{i-1}\op X\op\io_{j-i-2}\op D\op\io_{n-j-1}  ,
\end{align*}
and similarly $\wh\tau_{j;n}\circ\wh\si_{i;n} = \io_{i-1}\op X\op\io_{j-i-2}\op D\op\io_{n-j-1}$.  The $j\leq i-2$ case is virtually identical, as is the first part of \ref{P4}, both parts of \ref{P6}, the first part of \ref{P3}, and the $|i-j|\geq2$ case of the second.  The $|i-j|=1$ case follows from \ref{P2'}, as
\begin{align}
\nonumber \wh\tau_{i;n} \circ \wh\tau_{i+1;n} &= \io_{i-1}\op((D\op I)\circ(I\op D))\op\io_{n-i-2} \\
\label{e:DIID} \AND \wh\tau_{i+1;n} \circ \wh\tau_{i;n} &= \io_{i-1}\op((I\op D)\circ(D\op I))\op\io_{n-i-2}.
\end{align}
For the first part of \ref{P9}, we use Lemma \ref{l:tc}\ref{tc3}, with $a=\wh\th_{i;n}$, $b=\io_n$ and $c=\U$: 
\[
\wh\th_{i;n}\circ\wh\lam_n = \wh\th_{i;n}\circ(\io_n\op\U) 
= (\wh\th_{i;n}\circ\io_n)\op\U 
= (\io_n\circ\wh\th_{i;n})\op(\U\circ I) 
= (\io_n\op\U)\circ(\wh\th_{i;n}\op I) 
= \wh\lam_n\circ\wh\th_{i;n+1}.
\]
The second part is similar.

For every other relation from $\Om$, we adopt the following basic pattern.  We first use the tensor category axioms to write $\wh u=\io_k\op s\op\io_l$ and $\wh v=\io_k\op t\op\io_l$ for some $k,l\in\N$, and where $s,t\in\De^\oast$ are relatively simple, as in \eqref{e:DIID}, and then use relations \ref{P1'}--\ref{P6'} to show that $s\approx t$.  Table \ref{tab:P} shows the required calculations for the remaining relations from $\Om$.  (For the \ref{P7} entry, recall from Remark~\ref{r:i+1} that it suffices to consider only $j=i+1$.)
\begin{table}[ht]
\begin{center}
\begin{tabular}{|c|c|l|}
\hline
Relation & Label & Reduced relation $s\approx t$ \\
\hline
\ref{P1} & (a) & $X\circ X\approx \io_2$ \\
             & (b) & $(U\op\U)\circ(U\op\U)\approx U\op\U$ \\
             & (c) & $D\circ D\approx D\approx D\circ X\approx X\circ D$ \\
\hline
\ref{P2} & (d) & $X\circ(U\op\U\op I) \approx (I\op U\op\U)\circ X$  \\
             & (e) & $(U\op\U\op I)\circ(I\op U\op\U)\circ X\approx (U\op\U\op I)\circ(I\op U\op\U)$ \\
\hline
\ref{P5} & (f) & $(X\op I)\circ(I\op X)\circ(X\op I)\approx(I\op X)\circ(X\op I)\circ(I\op X)$ \\
             & (g) & $(X\op I)\circ(I\op D)\circ(X\op I)\approx(I\op X)\circ(D\op I)\circ(I\op X)$ \\
\hline
\ref{P7} & (h) & $D\circ(I\op U\op\U)\circ D \approx D$ \\
             & (i) & $(I\op U\op\U)\circ D\circ(I\op U\op\U) \approx I\op U\op\U$ \\
\hline
\ref{P8} & (j) & $\U\circ U \approx \io_0$ \\
             & (k) & $U\circ\U \approx U\op\U$ \\
\hline
\end{tabular}
\end{center}
\caption{Reduced relations $s\approx t$ required in the proof of Lemma \ref{l:a64P}; see the text for more details.}
\label{tab:P}
\end{table}

Many of the reduced relations $s\approx t$ in Table \ref{tab:P} are simply contained in $\Om$ itself; this is the case for (a), (c), (f), (g), (h) and (j), while (k) follows straight from \eqref{e:UU}.  This leaves us with (b), (d), (e) and~(i), and we treat these now.

\pfitem{(b)}  For this, we use \ref{P1'} and \eqref{e:UU}: $(U\op\U)\circ(U\op\U) = U\circ\U\circ U\circ\U \approx U\circ\io_0\circ\U = U\circ\U = U\op\U$.

\pfitem{(d)}  For this, we have
\begin{align*}
X\circ(U\op\U\op I) &= X\circ(\U\op U\op I) &&\text{by \eqref{e:UU}}\\
&\approx X\circ(\U\op (X\circ(I\op U))) &&\text{by \ref{P5'}}\\
&= X\circ X\circ (\U\op I\op U) &&\text{by Lemma \ref{l:tc}\ref{tc2}, with $a=\U$, $b=X$ and $c=I\op U$}\\
&\approx \io_2\circ (((I\op \U)\circ X) \op U) &&\text{by \ref{P1'} and \ref{P5'}}\\
&= (I\op \U\op U)\circ X &&\text{by Lemma \ref{l:tc}\ref{tc4}, with $a=I\op\U$, $b=X$ and $c=U$}\\
&= (I\op U\op \U)\circ X &&\text{by \eqref{e:UU}.}
\end{align*}
\firstpfitem{(e)}  Here we must show that $w\circ X\approx w$, where $w=(U\op\U\op I)\circ(I\op U\op\U)$.  First, by \eqref{e:UU} we have
\begin{equation}\label{e:UUUU}
w = (U\op\U\op I)\circ(I\op U\op\U) = ((U\op\U)\circ I)\op(I\circ(U\op\U)) = (U\op\U)\op(U\op\U) = U\op U\op\U\op\U.
\end{equation}
Also note that
\begin{align}
\nonumber (\U\op\U)\circ X = ((\U\circ I)\op\U)\circ X &= \U\circ (I\op\U)\circ X &&\text{by Lemma \ref{l:tc}\ref{tc3}, with $a=c=\U$ and $b=I$}\\
\nonumber &\approx \U\circ (\U\op I) &&\text{by \ref{P5'}}\\
\nonumber &= \U\op (\U\circ I) &&\text{by Lemma \ref{l:tc}\ref{tc2}, with $a=b=\U$ and $c=I$}\\
\label{e:UUX} &= \U\op\U.  
\end{align}
Using \eqref{e:UUUU} and \eqref{e:UUX}, and also Lemma \ref{l:tc}\ref{tc1}, with $a=U\op U$, $b=\U\op\U$ and $c=X$, we then have 
\[
w\circ X = (U\op U\op\U\op \U)\circ X = (U\op U)\op((\U\op\U)\circ X) \approx U\op U\op\U\op\U = w.
\]
\firstpfitem{(i)}  Here we have
\begin{align*}
 (I\op U&\op\U)\circ D\circ(I\op U\op\U) \\
&\approx ((I\op \U\op U)\circ D\circ(I\op U))\op\U &&\text{by \eqref{e:UU} and Lemma \ref{l:tc}\ref{tc3},}\\ & &&\text{\ \ with $a=(I\op U\op\U)\circ D$, $b=I\op U$ and $c=\U$}\\
&= ((I\op \U)\circ D\circ(I\op U))\op U\op\U &&\text{by Lemma \ref{l:tc}\ref{tc4},}\\ & &&\text{\ \ with $a=I\op\U$, $b=D\circ(I\op U)$ and $c=U$}\\
&\approx I\op U\op\U  &&\text{by \ref{P6'}.}  \qedhere
\end{align*}
\epf

\noindent As noted above, now that we have proved Lemmas \ref{l:a63P} and \ref{l:a64P}, the theorem is proved.
\epf

\begin{rem}
Surjectivity of $\Phi$ was also discussed by Martin in \cite[pp.~127-128]{Martin2008}.  
Comes \cite{Comes2020}, in his study of the so-called jellyfish partition categories, stated an alternative (tensor) presentation for~$\P$ in terms of the five generators:
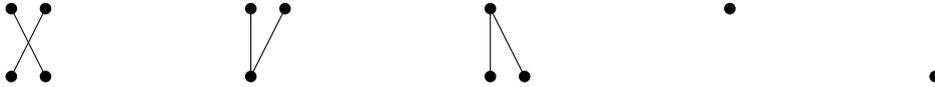
\begin{center}
\begin{tikzpicture}[scale=.45]
\begin{scope}[shift={(33,0)}]	
\uv1
\end{scope}
\begin{scope}[shift={(39,0)}]	
\lv1
\end{scope}
\begin{scope}[shift={(12,0)}]	
\lv1
\lv2
\uv1
\uv2
\stline12
\stline21
\end{scope}
\begin{scope}[shift={(19,0)}]	
\lv1
\uv1
\uv2
\stline11
\stline21
\end{scope}
\begin{scope}[shift={(26,0)}]	
\lv1
\lv2
\uv1
\stline11
\stline12
\end{scope}
\end{tikzpicture}
\end{center}
The proof given by Comes relied on some highly non-trivial results concerning cobordism categories and Frobenius algebras \cite{Kock2004,Abrams1996}.  Theorem \ref{t:P2} leads to an alternative proof of Comes' result, relying on no such heavy machinery; we simply re-write our presentation $\pres\De\Xi$ into his, using Tietze transformations.  \end{rem}

\begin{rem}\label{r:*PROP}
We mentioned in Section \ref{ss:P} that the category $\P$ has a natural involution $\al\mt\al^*$.  The presentation $\pres\De\Xi$ from Theorem \ref{t:P2} could be modified to give an involutory tensor category presentation, in which the involution is part of the signature (as well as $\circ$ and $\op$).  Here we would add the relations $X^*=X$, $D^*=D$, $U^*=\U$ and $\U^*= U$ (or eliminate one of $U$ or $\U$ altogether), and some of relations \ref{P1'}--\ref{P6'} could then be removed.  For example, we can remove ``$=X\circ D$'' from~\ref{P2'}, since using the other parts we have $X\circ D \approx X^*\circ D^* = (D\circ X)^* \approx D^* \approx D$.  Similarly, we would only need to keep one part of \ref{P5'}.

Recall from Section \ref{ss:P} that $\P$ is a so-called PROP, in the sense of \cite[Section 24]{MacLane1965}.  Thus, one could also give a PROP presentation for $\P$, in which case the subgroups (isomorphic to) $\SS_n$ are part of the background ``free'' data.  In this way, $X$ is simply the freely-existing non-trivial element of $\SS_2\sub\De^\oast_2$, and acts according to the PROP laws.  Thus, \ref{P3'} and the first part of \ref{P1'} are part of the free data.  So too are relations \ref{P4'} and \ref{P5'}.  These are less obvious, but follow from the PROP law
\[
(a\op b)\circ f_{n,l}=f_{m,k}\circ(b\op a) \qquad\text{for $a\in \C_{m,n}$ and $b\in \C_{k,l}$.}
\]
Indeed, since $f_{1,1}=X$, $f_{0,1}=\io_1$ and $f_{1,2}=(X\op I)\circ(I\op X)$, these follow from
\[
(U\op I)\circ f_{0,1} = f_{1,1}\circ(I\op U) \AND
(I\op D)\circ f_{1,2} = f_{1,2}\circ(D\op I),
\]
and other such identities.
\end{rem}

\begin{rem}
It is worth observing that the partition category $\P$ is finitely presented as a tensor category, as the sets $\De$ and $\Xi$ from the presentation in Theorem \ref{t:P2} are both finite.  In fact, all of the categories considered in this paper have the same finiteness property; see Theorems \ref{t:B2}, \ref{t:TL2}, \ref{t:PV2}, \ref{t:IB2}, \ref{t:V2} and so on.  The author believes it would be interesting to investigate this phenomenon, and in particular to seek necessary and/or sufficient conditions ensuring that a category $\C$ (satisfying Assumptions \ref{a:1}--\ref{a:6}) has a finite tensor presentation.
\end{rem}

\subsection{The Brauer category}\label{ss:B}

We now turn our attention to the Brauer category $\B$.  The argument here follows the same pattern as that of Section \ref{ss:P}, so we will for the most part abbreviate it.  In order to avoid a build-up of notation, we also re-use symbols such as $\Ga$, $\Om$, $\si$, $\tau$, etc., and indeed in later sections as well.

Since blocks of Brauer partitions have size $2$, Assumption \ref{a:1} holds in $\B$ with $d=2$.  As for Assumption \ref{a:2}, this time we take $\lamb_n\in\B_{n,n+2}$ and $\rhob_n\in\B_{n+2,n}$ to be the partitions shown in Figure~\ref{f:B_gens}.  The presentations for the Brauer monoids $\B_n$ required for Assumption~\ref{a:3} are taken from \cite{KM2006}.  For $n\in\N$ define an alphabet
\[
X_n = S_n\cup T_n \qquad\text{where} \qquad S_n = \set{\si_{i;n}}{1\leq i<n} \AND T_n = \set{\tau_{i;n}}{1\leq i<n}.
\]
Define a morphism $\phi_n:X_n^*\to\B_n:w\mt\ol w$, where the partitions $\ol x$ ($x\in X_n$) are shown in Figure~\ref{f:B_gens}.

\begin{figure}[ht]
\begin{center}
\begin{tikzpicture}[scale=.45]
\begin{scope}[shift={(16,0)}]	
\uvs{1,6}
\lvs{1,6,7,8}
\udotted16
\ddotted16
\stline11
\stline66
\darc78
\draw(0.5,1)node[left]{$\lamb_n=$};
\node()at(1,2.5){\tiny$1$};
\node()at(6,2.5){\tiny$n$};
\end{scope}
\begin{scope}[shift={(16,-5)}]	
\uvs{1,6,7,8}
\lvs{1,6}
\udotted16
\ddotted16
\stline11
\stline66
\uarc78
\draw(0.5,1)node[left]{$\rhob_n=$};
\node()at(1,2.5){\tiny$1$};
\node()at(6,2.5){\tiny$n$};
\end{scope}
\begin{scope}[shift={(0,0)}]	
\uvs{1,4,5,6,7,10}
\lvs{1,4,5,6,7,10}
\udotted14
\ddotted14
\udotted7{10}
\ddotted7{10}
\stline11
\stline44
\stline56
\stline65
\stline77
\stline{10}{10}
\draw(0.5,1)node[left]{$\sib_{i;n}=$};
\node()at(1,2.5){\tiny$1$};
\node()at(5,2.5){\tiny$i$};
\node()at(10,2.5){\tiny$n$};
\end{scope}
\begin{scope}[shift={(0,-5)}]	
\uvs{1,4,5,6,7,10}
\lvs{1,4,5,6,7,10}
\udotted14
\ddotted14
\udotted7{10}
\ddotted7{10}
\stline11
\stline44
\stline77
\stline{10}{10}
\uarcx56{.3}
\darcx56{.3}
\draw(0.5,1)node[left]{$\taub_{i;n}=$};
\node()at(1,2.5){\tiny$1$};
\node()at(5,2.5){\tiny$i$};
\node()at(10,2.5){\tiny$n$};
\end{scope}
\end{tikzpicture}
\caption{Brauer generators $\ol x\in\B$ ($x\in\Ga$).}
\label{f:B_gens}
\end{center}
\end{figure}

Let $R_n$ be the following set of relations over $X_n$:
\begin{align}
\label{B1}\tag*{(B1)} \si_{i;n}^2 = \iota_n ,&&&  \tau_{i;n}^2 = \tau_{i;n}  = \tau_{i;n}\si_{i;n} = \si_{i;n}\tau_{i;n},  \\
\label{B2}\tag*{(B2)} \si_{i;n}\si_{j;n} = \si_{j;n}\si_{i;n} ,&&&  \tau_{i;n}\tau_{j;n} = \tau_{j;n}\tau_{i;n} \COMMa  \si_{i;n}\tau_{j;n} = \tau_{j;n}\si_{i;n},  &&\text{if $|i-j|>1$,}   \\
\label{B3}\tag*{(B3)}  \si_{i;n}\si_{j;n}\si_{i;n} = \si_{j;n}\si_{i;n}\si_{j;n} ,&&&     \si_{i;n}\tau_{j;n}\si_{i;n} = \si_{j;n}\tau_{i;n}\si_{j;n} \COMMa \tau_{i;n}\si_{j;n}\tau_{i;n} = \tau_{i;n},  &&\text{if $|i-j|=1$.}   
\end{align}

\begin{thm}[cf.~\cite{KM2006}]\label{t:Bn}
For any $n\in\N$, the Brauer monoid $\B_n$ has presentation $\pres{X_n}{R_n}$ via $\phi_n$.  \epfres
\end{thm}

Now let $\Ga\equiv L\cup R\cup X$ be the digraph over $\N$ as in \eqref{e:Ga} and \eqref{e:st}, and let $\phi:\Ga^*\to\B$ be the morphism given in \eqref{e:phi}.  Let $\Om$ be the set of relations over $\Ga$ consisting of $\bigcup_{n\in\N}R_n$ and additionally:
\begin{align}
\label{B4}\tag*{(B4)} \lam_n\rho_n = \iota_n ,&&& \rho_n\lam_n = \tau_{n+1;n+2}, \\
\label{B5}\tag*{(B5)} \th_{i;n}\lam_n = \lam_n\th_{i;n+2} ,&&& \rho_n\th_{i;n} = \th_{i;n+2}\rho_n, &&\text{for $\th\in\{\si,\tau\}$.}  
\end{align}

Here is the first main result of this section, expressed in terms of the above notation.  The proof is essentially identical to that of Theorem \ref{t:P1}; the role of $\ve_{n+1;n+1}$ in that proof is played instead by $\tau_{n+1;n+2}$, and we use the map $\B_n\to\B_{n+2}:\al\mt\al^+=\al\op\iob_2$.

\begin{thm}\label{t:B1}
The Brauer category $\B$ has presentation~$\pres\Ga\Om$ via $\phi$.  \epfres
\end{thm}

\begin{rem}\label{r:a4''}
As in Remarks \ref{r:a4} and \ref{r:a4'}, we could have made different choices for the one-sided units $\lamb_n,\rhob_n\in\B$ to those in Figure \ref{f:B_gens}.  However, and in contrast to the situation with the partition category (cf.~Remark \ref{r:a4'}), the choices we have made here seem the most ``natural'', in the sense that any other choice leads to arguably less convenient forms of the relations \ref{B4} and \ref{B5}.
\end{rem}

Next, let $\De$ be the digraph over $\N$ with three edges: 
\[
X:2\to2 \COMMA U:2\to0 \COMMA \U:0\to2.
\]
Define the morphism $\Phi:\De^\oast\to\B:w\mt\ul w$, where the partitions $\ul x$ ($x\in\De$) are shown in Figure \ref{f:B_Tgens}.  It will again be convenient to write $I=\io_1$, and $\ul I=\iob_1\in\B_1$.

\begin{figure}[ht]
\begin{center}
\begin{tikzpicture}[scale=.45]
\begin{scope}[shift={(7,0)}]	
\uv1
\uv2
\uarc12
\draw(0.5,1)node[left]{$\ul U=$};
\end{scope}
\begin{scope}[shift={(14,0)}]	
\lv1
\lv2
\darc12
\draw(0.5,1)node[left]{$\ul \U=$};
\end{scope}
\begin{scope}[shift={(0,0)}]	
\lv1
\lv2
\uv1
\uv2
\stline12
\stline21
\draw(0.5,1)node[left]{$\ul X=$};
\end{scope}
\begin{scope}[shift={(21,0)}]	
\uv1
\lv1
\stline11
\draw(0.5,1)node[left]{$\ul I=$};
\end{scope}
\end{tikzpicture}
\caption{Brauer generators $\ul x\in\B$ ($x\in\De$), as well as $\ul I$.}
\label{f:B_Tgens}
\end{center}
\end{figure}

Let $\Xi$ be the set consisting of the following relations over $\De$:
\begin{gather}
\tag*{(B1)$'$} \label{B1'} X\circ X=\io_2 \COMMA \U\circ U = \io_0   \COMMA X\circ U=U \COMMA \U\circ  X=\U, \\
\tag*{(B2)$'$} \label{B2'} (X\op I)\circ(I\op X)\circ(X\op I)=(I\op X)\circ(X\op I)\circ(I\op X),\\
\tag*{(B3)$'$} \label{B3'} (I\op \U)\circ (U\op I) = I = (\U\op I)\circ(I\op U),\\
\tag*{(B4)$'$} \label{B4'} (X\op I)\circ(I\op U) = (I\op X)\circ(U\op I)  \COMMa \phantom{=} \ (\U\op I)\circ(I\op X) = (I\op \U)\circ(X\op I).  
\end{gather}
The second main result of this section is the following, expressed in terms of the above notation:

\begin{thm}\label{t:B2}
The Brauer category $\B$ has tensor presentation $\pres\De\Xi$ via $\Phi$.
\end{thm}

\pf
The proof follows the same outline as that of Theorem \ref{t:P2}.  We need only verify Assumption~\ref{a:6}.  Condition \ref{a61} is checked diagrammatically, and for the remaining assumptions we use the morphism $\Ga^*\to\De^\oast:w\mt\wh w$, given by
\[
\wh\si_{i;n} = \io_{i-1}\op X\op\io_{n-i-1} \COMMA
\wh\tau_{i;n} = \io_{i-1}\op U\op\U\op\io_{n-i-1} \COMMA
\wh\lam_n = \io_n\op\U  \COMMA
\wh\rho_n = \io_n\op U.
\]
Condition \ref{a62} is again verified diagrammatically.  Condition~\ref{a63} is proved as in Lemma~\ref{l:a63P}: e.g., we have $U_{m,n}\approx\wh\tau_{m+1;m+n+2}\circ\cdots\circ\wh\tau_{m+n;m+n+2}\circ\wh\rho_{m+n}$.  Condition~\ref{a64} is proved as in Lemma \ref{l:a64P}; the details are again mostly very similar, with the only substantial exceptions being the second and third parts of \ref{B3}, so we treat just these.  The second part of \ref{B3} follows from
\begin{align*}
(X\op I)\circ(I\op U\op\U)\circ(X\op I) &= (((X\op I)\circ(I\op U))\op\U)\circ(X\op I) &&\text{by Lemma \ref{l:tc}\ref{tc3}}\\
&\approx (((I\op X)\circ(U\op I))\op\U)\circ(X\op I) &&\text{by \ref{B4'}}\\
&= (I\op X)\circ(U\op I\op\U)\circ(X\op I) &&\text{by Lemma \ref{l:tc}\ref{tc3}}\\
&= (I\op X)\circ(U\op ((I\op\U)\circ(X\op I))) &&\text{by Lemma \ref{l:tc}\ref{tc1}}\\
&\approx (I\op X)\circ(U\op ((\U\op I)\circ(I\op X))) &&\text{by \ref{B4'}}\\
&= (I\op X)\circ(U\op \U\op I)\circ(I\op X) &&\text{by Lemma \ref{l:tc}\ref{tc1}.}
\end{align*}
For the third part of \ref{B3}, we consider only the $j=i+1$ case ($j=i-1$ is similar), which follows from
\begin{align*}
(U\op\U\op I)\circ(I\op X)&\circ(U\op\U\op I) = (U\op((\U\op I)\circ(I\op X))) \circ(U\op\U\op I) &&\hspace{-.2cm}\text{by Lemma \ref{l:tc}\ref{tc1}}\\
&\approx (U\op((I\op \U)\circ(X\op I))) \circ(U\op\U\op I)  &&\hspace{-.5cm}\text{by \ref{B4'}}\\
&= U\op((I\op \U)\circ(X\op I) \circ(U\op\U\op I))  &&\hspace{-.5cm}\text{by Lemma \ref{l:tc}\ref{tc1}}\\
&= U\op((I\op \U)\circ((X\circ U) \op(I\circ (\U\op I))))  \\
&\approx U\op((I\op \U)\circ (U \op(\U\op I)))   &&\hspace{-.5cm}\text{by \ref{B1'}}\\
&= U\op((I\op \U)\circ (\U \op U\op I))   &&\hspace{-.5cm}\text{by \eqref{e:UU}} \\
&= U\op(\U\op((I\op \U)\circ (U\op I)))   &&\hspace{-.5cm}\text{by Lemma \ref{l:tc}\ref{tc2}}\\
&\approx U\op\U\op I   &&\hspace{-.5cm}\text{by \ref{B3'}.}  \qedhere
\end{align*}
\epf

\begin{rem}
As in Remark \ref{r:*PROP}, certain relations could be removed from the presentation by considering the involution and/or PROP structure.  The involution was built into the presentation for~$\B$ given by Lehrer and Zhang in \cite[Theorem 2.6]{LZ2015}; note that $I=\io_1$ is explicitly listed as a generator there, as well as relations such as $I\circ I=I$ and $(I\op I)\circ X=X$.
\end{rem}

%\begin{rem}
%The second part of relation \ref{B3'} is redundant.  Indeed, given the first part, and several other relations from $\Xi$, we have
%\begin{align*}
%(\U\op I)\circ(I\op U) 
%&\approx ((\U\circ X)\op(I\circ I))\op(I\op U)
%\approx (\U\op I)\circ (X\op I)\op(I\op U)\\
%&\approx (\U\op I)\circ (I\op X)\op(U\op I) 
%\approx (I\op \U)\circ (X\op I)\op(U\op I) \\
%&\approx (I\op \U)\circ ((X\circ U)\op (I\circ I))
%\approx (I\op \U)\circ (U\op I)
%\approx I.
%\end{align*}
%\end{rem}

\subsection{The Temperley-Lieb category}\label{ss:TL}

The methods of Sections \ref{ss:P} and \ref{ss:B} can also be used to quickly obtain presentations for the Temperley-Lieb category $\TL$.  Here we let $\Ga$ be the digraph over $\N$ with edge set 
\[
L\cup R\cup \bigcup_{n\in\N}X_n \qquad\text{where}\qquad L=\set{\lam_n}{n\in\N} \COMMA R=\set{\rho_n}{n\in\N} \COMMA X_n = \set{\tau_{i;n}}{1\leq i<n}.
\]
The partitions $\ol x\in\TL$ ($x\in\Ga$) are as already shown in Figure \ref{f:B_gens}.  We define in the usual way the morphism $\phi:\Ga^*\to\TL:w\mt\ol w$, and this time $\Om$ is the following set of relations:
\begin{gather}
\label{TL1}\tag*{(TL1)}  \tau_{i;n}^2 = \tau_{i;n}  \COMMA \tau_{i;n}\tau_{j;n} = \tau_{j;n}\tau_{i;n} \text{ if $|i-j|>1$} \COMMA \tau_{i;n}\tau_{j;n}\tau_{i;n} = \tau_{i;n} \text{ if $|i-j|=1$,}\\
\label{TL2}\tag*{(TL2)} \lam_n\rho_n = \iota_n \COMMA \rho_n\lam_n = \tau_{n+1;n+2} \COMMA \tau_{i;n}\lam_n = \lam_n\tau_{i;n+2} \COMMA \rho_n\tau_{i;n} = \tau_{i;n+2}\rho_n .
\end{gather}
For fixed $n$, relations \ref{TL1} constitute defining relations for the Temperley-Lieb monoid~$\TL_n$; cf.~\cite{JE2021,BDP2002,Jones1983_2}.

\begin{thm}\label{t:TL1}
The Temperley-Lieb category $\TL$ has presentation~$\pres\Ga\Om$ via $\phi$.  \epfres
\end{thm}

Next let $\De$ be the digraph over $\N$ with edges $U:2\to0$ and $\U:0\to2$.  We have a morphism $\Phi:\De^\oast\to\TL:w\mt\ul w$, where the partitions $\ul U$ and $\ul\U$ from $\TL$ are as in Figure~\ref{f:B_Tgens}.  Let $\Xi$ be the set consisting of the following relations over $\De$, where~$I=\io_1$:
\[
\U\circ U = \io_0 \AND (I\op \U)\circ (U\op I) = I = (\U\op I)\circ(I\op U).
\]

\begin{thm}\label{t:TL2}
The Temperley-Lieb category $\TL$ has tensor presentation $\pres\De\Xi$ via $\Phi$.  \epfres
\end{thm}

\begin{rem}
It is not clear to whom Theorem \ref{t:TL2} should be attributed, though it appears to be folklore.  See the discussion in \cite[Section~3.1]{Abramsky2008}, which refers to \cite{FY1989,DP2003}; the results from these papers concern so-called free pivotal categories, but are equivalent to Theorem~\ref{t:TL2}.
\end{rem}

\begin{rem}
As explained in \cite[p.~264.]{Jones1994} and \cite[p.~873]{HR2005}, the monoid $\PlP_n$ of planar partitions of degree $n$ is isomorphic to $\TL_{2n}$, the Temperley-Lieb monoid of degree $2n$.  Although the categories~$\PlP$ and $\TL$ are not similarly isomorphic (as $\TL$ also contains partitions of odd degree), one can derive presentations for $\PlP$ from those for $\TL$ above.  These are in terms of generators $\veb_{i;n}$, $\taub_{i;n}$, $\lamb_n$ and $\rhob_n$ (as pictured in Figure \ref{f:P_gens}), or $\ol U$, $\ol\U$ and $\ol D$ (as pictured in Figure \ref{f:P_Tgens}).  Other diagram categories could be treated similarly \cite{DEG2017,Grood2006,BH2014,MM2014}.
\end{rem}

\subsection{Linear diagram categories}\label{ss:L}

Fix a field $\kk$, and an element $\de\in\kk\sm\{0\}$.  For $m,n\in\N$, we denote by $\P_{m,n}^\de$ the $\kk$-vector space with basis~$\P_{m,n}$: i.e., the set of all formal $\kk$-linear combinations of partitions from $\P_{m,n}$.  We also write
\[
\P^\de=\bigcup_{m,n\in\N}\P_{m,n}^\de\]
for the set of all such combinations.  This set is a category (over $\N$) with composition $\star$ defined as follows.  Consider basis elements $\al\in\P_{m,n}$ and $\be\in\P_{n,k}$.  We write $m(\al,\be)$ for the number of ``floating'' components in the product graph $\Pi(\al,\be)$, as defined in Section \ref{ss:preDC}: i.e., the number of components contained entirely in $[n]''$.  (So $m(\al,\be)=1$ for $\al,\be$ in Figure \ref{f:P_product}.)  The product $\al\star\be\in\P_{m,k}^\de$ is defined by
\[
\al\star\be = \de^{m(\al,\be)} \al\be.
\]
So $\al\star\be$ is a scalar multiple of the basis element $\al\be=\al\circ\be\in\P_{m,k}$.  This composition on basis elements is then extended to all of $\P^\de$ by $\kk$-linearity.  The operation $\op$ and the involution~${}^*$ on $\P$ also extend to corresponding operations on $\P^\de$.

Here we call $\P^\de$ the linear partition category (associated to $\kk$ and $\de$), in order to distinguish it from the partition category $\P$.  Most authors simply call~$\P^\de$ the partition category.  We also have linear versions of the Brauer and Temperley-Lieb categories:~$\B^\de$ and~$\TL^\de$.

The presentations for the diagram categories given in Sections \ref{ss:P}--\ref{ss:TL} may be easily modified to yield presentations for their linear counterparts.  To do so, we use the method of \cite[Section 6]{JEgrpm}; this was originally formulated for algebras, but works virtually unchanged for (tensor) categories.

Let $\C$ be any of $\P$, $\B$ or $\TL$, and let $\Ga$ be the digraph defined in Section \ref{ss:P}, \ref{ss:B} or \ref{ss:TL}, as appropriate.  Consider some path $w\in\Ga^*$.  If $w$ is empty or a single edge, we define $m(w)=0$.  Otherwise, write $w=x_1\cdots x_k$, where $k\geq2$ and each $x_i\in\Ga$.  In the category $\C^\de$, we have $\ol x_1\star\cdots\star\ol x_k = \de^{m(w)} \ol w$, where $m(w) = m(\ol x_1,\ol x_2) + m(\ol{x_1x_2},\ol x_3) +\cdots+ m(\ol{x_1\cdots x_{k-1}},\ol x_k)$.  We may then obtain a presentation for $\C^\de$ by replacing each relation $u=v$ in a presentation for $\C$ (cf.~Theorems~\ref{t:P1},~\ref{t:B1} and~\ref{t:TL1}) by $\de^{m(v)}u=\de^{m(u)}v$:
\bit
\item For $\C=\P$, we replace $\ve_{i;n}^2=\ve_{i;n}$ and $\lam_n\rho_n=\io_n$ by $\ve_{i;n}^2=\de\ve_{i;n}$ and $\lam_n\rho_n=\de\io_n$.
\item For $\C=\B$ and $\TL$, we replace $\tau_{i;n}^2=\tau_{i;n}$ and $\lam_n\rho_n=\io_n$ by $\tau_{i;n}^2=\de\tau_{i;n}$ and $\lam_n\rho_n=\de\io_n$.
\eit

We may similarly obtain tensor presentations for $\P^\de$, $\B^\de$ and $\TL^\de$ from the corresponding presentations for $\P$, $\B$ and $\TL$ (cf.~Theorems \ref{t:P2}, \ref{t:B2} and \ref{t:TL2}).  For all three, this simply amounts to replacing $\U\circ U=\io_0$ by $\U\circ U=\de\io_0$.

\section{Categories of (partial) vines, braids and transformations}\label{s:BC}

Next we consider a number of natural categories of braids and vines, namely the partial vine category~$\PV$ (Section \ref{ss:PV}), the partial braid category $\IB$ (Section \ref{ss:IB}), and the (full) vine category $\V$ (Section~\ref{ss:V}).  Again, we apply the general theory developed in Section \ref{s:C} to obtain presentations for each category.  The connectivity of the category $\V$ is slightly different from the other categories considered so far, requiring the use of Theorem \ref{t:pres3} instead of Theorem \ref{t:pres2}.

In Section \ref{ss:T} we apply the results of Sections \ref{ss:PV}--\ref{ss:V} to quickly obtain presentations for several categories of (partial) transformations/mappings, and in Section \ref{ss:O} obtain analogous results for categories of isotone (order-preserving) mappings.

\subsection{Preliminaries}\label{ss:pre_braids}

We begin by defining the partial vine category $\PV$, following \cite{Artin1947,Lavers1997,JE2007b,ER2023}.  By a string we mean a smooth, tame embedding $\s$ of the unit interval $[0,1]$ into $\RRR^3$, such that:
\bit
\item the $z$-coordinate of $\s(t)$ is $1-t$ for all $t$,
\item $\s(0)=(a,0,1)$ and $\s(1)=(b,0,0)$ for some $a,b\in\PP$.
\eit
For such a string $\s$, we write $I(\s)=a$ and $T(\s)=b$, which codify the initial and terminal points of $\s$.

A partial vine is (a homotopy class of) a tuple $\al=(\s_1,\ldots,\s_k)$ of strings, where $k\geq0$, satisfying: 
\bit
\item $I(\s_1)<\cdots<I(\s_k)$,
\item if $\s_i(t)=\s_j(t)$ for some $t\in[0,1]$, then $\s_i(u)=\s_j(u)$ for all $u\in[t,1]$.
\eit
For such a partial vine $\al$, we define $I(\al)=\big\{I(\s_1),\ldots,I(\s_k)\big\}$ and $T(\al)=\big\{T(\s_1),\ldots,T(\s_k)\big\}$.  Examples of partial vines are shown in Figure \ref{f:vines}.  

For $m,n\in\N$, we write $\PV_{m,n}$ for the set of all partial vines $\al$ with $I(\al)\sub[m]$ and $T(\al)\sub[n]$, and we define the vine category:
\[
\PV = \bigset{(m,\al,n)}{m,n\in\N,\ \al\in\PV_{m,n}},
\]
with domain, range and composition operations as follows.  For $m,n,q\in\N$, and for $\al\in\PV_{m,n}$ and $\be\in\PV_{n,q}$, we define
\[
\bd(m,\al,n)=m \COMMA \br(m,\al,n)=n \COMMA (m,\al,n)\circ(n,\be,q) = (m,\al\be,q),
\]
where the product $\al\be=\al\circ\be\in\PV_{m,q}$ is obtained by:
\bit
\item placing a translated copy of $\be$ below $\al$,
\item scaling so that the resulting object lies in the region $0\leq z\leq 1$,
\item removing any string fragments that do not join top to bottom.
\eit
Figure \ref{f:vines} gives an example product.  To avoid clutter in our notation, we will typically identify an element $(m,\al,n)$ of $\PV$ with the partial vine $\al\in\PV_{m,n}$ itself, but regard $m$ and $n$ as ``encoded'' in~$\al$, writing $\bd(\al)=m$ and $\br(\al)=n$.  In this way, the hom-sets of $\PV$ are the $\PV_{m,n}$ ($m,n\in\N$), and the endomorphism monoids are the partial vine monoids $\PV_n=\PV_{n,n}$ of \cite{JE2007b}.  The units in $\PV_n$ form the usual Artin braid groups $\B_n$ \cite{Artin1947}; since we will no longer refer to the Brauer monoids, we will re-use the symbol $\B$, as is standard.

\begin{figure}[ht]
\begin{center}
\scalebox{.9}{
\begin{tikzpicture}[scale=.666]

\begin{scope}[shift={(0,6)}]
\understringx2440
\overstringx3432
\overstringx44{1.7}2
\overstringx1422
\overstringx{2.3}{1.5}1{.5}
\overstringx2221
\overstringx{1.7}2{2.3}{1.5}
\overstringx3220
\understringx1{.5}10
\understringx2120
\foreach \x in {1,2,3,4} {\fill (\x,4)circle(.1);}
\foreach \x in {1,2,3,4,5} {\fill (\x,0)circle(.1);}
\draw(0.6,2)node[left]{$\al=$};
\draw[->](6.5,-1)--(8.5,-1);
\end{scope}

\begin{scope}[shift={(0,0)}]
\understringx3431
\understringx2431
\understringx4431
\understringx3130
\overstringx5452
\overstringx5220
\foreach \x in {1,2,3,4,5} {\fill (\x,4)circle(.1);}
\foreach \x in {1,2,3} {\fill (\x,0)circle(.1);}
\draw(0.6,2)node[left]{$\be=$};
\end{scope}

\begin{scope}[shift={(9,5)}]
\understringx2440
\overstringx3432
\overstringx44{1.7}2
\overstringx1422
\overstringx{2.3}{1.5}1{.5}
\overstringx2221
\overstringx{1.7}2{2.3}{1.5}
\overstringx3220
\understringx1{.5}10
\understringx2120
\foreach \x in {1,2,3,4} {\fill (\x,4)circle(.1);}
\foreach \x in {1,2,3,4,5} {\fill (\x,0)circle(.1);}
\draw[->](6.5,0)--(8.5,0);
\end{scope}

\begin{scope}[shift={(9,1)}]
\understringx3431
\understringx2431
\understringx4431
\understringx3130
\overstringx5452
\overstringx5220
\foreach \x in {1,2,3,4,5} {\fill (\x,4)circle(.1);}
\foreach \x in {1,2,3} {\fill (\x,0)circle(.1);}
\end{scope}

\begin{scope}[shift={(18,5)},yscale=0.5]
\understringx2440
\overstringx3432
\overstringx1422
\overstringx2221
\overstringx3220
\understringx2120
\end{scope}

\begin{scope}[shift={(18,3)},yscale=0.5]
\understringx2431
\understringx4431
\understringx3130
\end{scope}

\begin{scope}[shift={(18,3)}]
\foreach \x in {1,2,3,4} {\fill (\x,4)circle(.1);}
\foreach \x in {1,2,3} {\fill (\x,0)circle(.1);}
\draw(4.4,2)node[right]{$=\al\be$};
\end{scope}

\end{tikzpicture}
}
\caption{Calculating $\al\be$, where $\al\in\PV_{4,5}$ and $\be\in\PV_{5,3}$.}
\label{f:vines}
\end{center}
\end{figure}

The category $\PV$ is a (strict) tensor category, with $\op$ defined as for diagram categories:  $\al\op\be$ is obtained by placing a translated copy of $\be$ to the right of $\al$.  Unsurprisingly, $\PV$ is a (strict) braided tensor category in the sense of Joyal and Street \cite{JS1993}, and hence a PROB (PROducts and Braids).

We say that $\al\in\PV_{m,n}$ is: a (full) vine if $I(\al)=[m]$; or a partial braid if the strings of $\al$ do not intersect.  The sets $\V$ and $\IB$ of all vines and partial braids are subcategories of $\PV$: the vine category, and the partial braid category.  Endomorphisms in~$\V$ and $\IB$ form the vine monoids $\V_n$ \cite{Lavers1997} and inverse braid monoids $\IB_n$ \cite{EL2004}.  The category $\IB$ was studied in \cite[Section 12]{ER2023}, where it was observed to be an inverse category in the sense of \cite{Kastl1979} and \cite[Section 2.3.2]{CL2002}: for any $\al\in\IB$, the partial braid obtained by reflecting $\al$ in the plane $z=\frac12$ is the unique element $\be$ of $\IB$ satisfyng $\al=\al\be\al$ and~$\be=\be\al\be$.  The categories $\V$ and $\IB$ are both closed under $\op$, and both contain the braid groups $\B_n$ ($n\in\N$), so they are both PROBs.  However, while all hom-sets in $\PV$ and $\IB$ are non-empty, we have $\V_{m,0}=\emptyset$ for $m\geq1$.

\subsection{The partial vine category}\label{ss:PV}

In this section, we apply Theorems \ref{t:pres} and \ref{t:pres2} to obtain presentations for the partial vine category~$\PV$.  

First, note that Assumption \ref{a:1} holds in $\PV$ with $d=1$.  For Assumption \ref{a:2}, we take the partial vines $\lamb_n\in\PV_{n,n+1}$ and $\rhob_n\in\PV_{n+1,n}$ pictured in Figure \ref{f:PV_gens}.  For Assumption \ref{a:3}, we require presentations for the partial vine monoids $\PV_n$ ($n\in\N$), and these are given in \cite{JE2007b}.  For $n\in\N$, we define the alphabet
\begin{align}
\nonumber X_n = S_n\cup S_n^{-1}\cup E_n\cup M_n\cup H_n \qquad\text{where}\qquad 
S_n &= \set{\si_{i;n}}{1\leq i<n} , & M_n &= \set{\mu_{i;n}}{1\leq i< n}, \\
\nonumber S_n^{-1} &= \set{\si_{i;n}^{-1}}{1\leq i<n} , & H_n &= \set{\eta_{i;n}}{1\leq i< n}. \\
\label{e:PVn_gens} E_n &= \set{\ve_{i;n}}{1\leq i\leq n},
\end{align}
We define a morphism $\phi_n:X_n^*\to\PV_n:w\mt\ol w$, where the partial vines $\ol x$ ($x\in X_n$) are also shown in Figure~\ref{f:PV_gens}.

\begin{figure}[ht]
\begin{center}
\begin{tikzpicture}[scale=.42]
\begin{scope}[shift={(0,0)}]	
\udotted14
\ddotted14
\udotted7{10}
\ddotted7{10}
\stline11
\stline44
\ststring65
\ststring56
\stline77
\stline{10}{10}
\draw(0.5,1)node[left]{$\sib_{i;n}=$};
\node()at(1,2.5){\tiny$1$};
\node()at(5,2.5){\tiny$i$};
\node()at(10,2.5){\tiny$n$};
\uvs{1,4,5,6,7,10}
\lvs{1,4,5,6,7,10}
\end{scope}
\begin{scope}[shift={(0,-5)}]	
\udotted14
\ddotted14
\udotted7{10}
\ddotted7{10}
\stline11
\stline44
\ststring56
\ststring65
\stline77
\stline{10}{10}
\draw(0.5,1)node[left]{$\sib_{i;n}^{-1}=$};
\node()at(1,2.5){\tiny$1$};
\node()at(5,2.5){\tiny$i$};
\node()at(10,2.5){\tiny$n$};
\uvs{1,4,5,6,7,10}
\lvs{1,4,5,6,7,10}
\end{scope}
\begin{scope}[shift={(0,-10)}]	
\uvs{1,4,5,6,10}
\lvs{1,4,5,6,10}
\udotted14
\ddotted14
\udotted6{10}
\ddotted6{10}
\stline11
\stline44
\stline66
\stline{10}{10}
\draw(0.5,1)node[left]{$\veb_{i;n}=$};
\node()at(1,2.5){\tiny$1$};
\node()at(5,2.5){\tiny$i$};
\node()at(10,2.5){\tiny$n$};
\end{scope}
\begin{scope}[shift={(15,0)}]	
\udotted14
\ddotted14
\udotted7{10}
\ddotted7{10}
\stline11
\stline44
\stline55
\stline65
\stline77
\stline{10}{10}
\draw(0.5,1)node[left]{$\mub_{i;n}=$};
\node()at(1,2.5){\tiny$1$};
\node()at(5,2.5){\tiny$i$};
\node()at(10,2.5){\tiny$n$};
\uvs{1,4,5,6,7,10}
\lvs{1,4,5,6,7,10}
\end{scope}
\begin{scope}[shift={(15,-5)}]	
\udotted14
\ddotted14
\udotted7{10}
\ddotted7{10}
\stline11
\stline44
\stline56
\stline66
\stline77
\stline{10}{10}
\draw(0.5,1)node[left]{$\etab_{i;n}=$};
\node()at(1,2.5){\tiny$1$};
\node()at(5,2.5){\tiny$i$};
\node()at(10,2.5){\tiny$n$};
\uvs{1,4,5,6,7,10}
\lvs{1,4,5,6,7,10}
\end{scope}
\begin{scope}[shift={(30,0)}]	
\uvs{1,6}
\lvs{1,6,7}
\udotted16
\ddotted16
\stline11
\stline66
\draw(0.5,1)node[left]{$\lamb_n=$};
\node()at(1,2.5){\tiny$1$};
\node()at(6,2.5){\tiny$n$};
\end{scope}
\begin{scope}[shift={(30,-5)}]	
\uvs{1,6,7}
\lvs{1,6}
\udotted16
\ddotted16
\stline11
\stline66
\draw(0.5,1)node[left]{$\rhob_n=$};
\node()at(1,2.5){\tiny$1$};
\node()at(6,2.5){\tiny$n$};
\end{scope}
\end{tikzpicture}
\caption{Partial vine generators $\ol x\in\PV$ ($x\in\Ga$).}
\label{f:PV_gens}
\end{center}
\end{figure}

Now let $R_n$ be the following set of relations over $X_n$:
\begin{align}
\tag*{(PV1)} \label{PV1} & \si_{i;n}\si_{i;n}^{-1} = \si_{i;n}^{-1}\si_{i;n} = \io_n \COMMa \ve_{i;n}^2=\ve_{i;n} \COMMa \ve_{i;n}\ve_{j;n}=\ve_{j;n}\ve_{i;n}, \\
\tag*{(PV2)} \label{PV2} & \mu_{i;n} = \mu_{i;n}^2 = \eta_{i;n}\mu_{i;n} = \si_{i;n}\mu_{i;n} = \eta_{i;n}\si_{i;n} \COMMa \eta_{i;n} = \eta_{i;n}^2 = \mu_{i;n}\eta_{i;n} = \si_{i;n}\eta_{i;n}  = \mu_{i;n}\si_{i;n}, \\
\tag*{(PV3)} \label{PV3} & \mu_{i;n}\mu_{i+1;n}=\mu_{i;n}\si_{i+1;n} \COMMa \mu_{i;n}\eta_{i+1;n} = \mu_{i;n} \COMMa \eta_{i+1;n}\eta_{i;n}=\eta_{i+1;n}\si_{i;n} \COMMa \eta_{i+1;n}\mu_{i;n} = \eta_{i+1;n}, \\
\tag*{(PV4)} \label{PV4} & \mu_{i+1;n}\mu_{i;n} = \mu_{i;n}\mu_{i+1;n}\mu_{i;n} = \mu_{i+1;n}\mu_{i;n}\mu_{i+1;n} \COMMa \eta_{i;n}\eta_{i+1;n} = \eta_{i;n}\eta_{i+1;n}\eta_{i;n} = \eta_{i+1;n}\eta_{i;n}\eta_{i+1;n}, \\
\tag*{(PV5)} \label{PV5} & \mu_{i;n}\ve_{i+1;n}=\mu_{i;n} \COMMa \ve_{i+1;n}\mu_{i;n}=\ve_{i+1;n} \COMMa \eta_{i;n}\ve_{i;n}=\eta_{i;n} \COMMa \ve_{i;n}\eta_{i;n}=\ve_{i;n}, \\
\tag*{(PV6)} \label{PV6} & \mu_{i+1;n}\si_{i;n}=\si_{i;n}\si_{i+1;n}\mu_{i;n}\mu_{i+1;n} \COMMa \eta_{i;n}\si_{i+1;n}=\si_{i+1;n}\si_{i;n}\eta_{i+1;n}\eta_{i;n}, \\
\tag*{(PV7)} \label{PV7} & \si_{i;n}\ve_{i;n}=\ve_{i+1;n}\si_{i;n} \COMMa \si_{i;n}\ve_{i+1;n}=\ve_{i;n}\si_{i;n} \COMMa \si_{i;n}^2\ve_{i;n}=\ve_{i;n}, \\
\tag*{(PV8)} \label{PV8} & \si_{i;n}\ve_{i;n}\ve_{i+1;n}=\ve_{i;n}\ve_{i+1;n}=\mu_{i;n}\ve_{i;n}=\eta_{i;n}\ve_{i+1;n},  \\
\tag*{(PV9)} \label{PV9} & \si_{i;n}\si_{j;n} = \si_{j;n}\si_{i;n} \COMMa \mu_{i;n}\mu_{j;n} = \mu_{j;n}\mu_{i;n} \COMMa \eta_{i;n}\eta_{j;n} = \eta_{j;n}\eta_{i;n}, &&\hspace{-3.5cm}\text{if $|i-j|>1$,} \\
\tag*{(PV10)} \label{PV10} & \si_{i;n}\mu_{j;n}=\mu_{j;n}\si_{i;n} \COMMa \si_{i;n}\eta_{j;n}=\eta_{j;n}\si_{i;n}, &&\hspace{-3.5cm}\text{if $|i-j|>1$,} \\
\tag*{(PV11)} \label{PV11} & \si_{i;n}\si_{j;n}\si_{i;n} = \si_{j;n}\si_{i;n}\si_{j;n}, &&\hspace{-3.5cm}\text{if $|i-j|=1$,} \\
\tag*{(PV12)} \label{PV12} & \mu_{i;n}\eta_{j;n} = \eta_{j;n}\mu_{i;n} \COMMa \si_{i;n}\ve_{j;n}=\ve_{j;n}\si_{i;n}, &&\hspace{-3.5cm}\text{if $j\not=i,i+1$,} \\
\tag*{(PV13)} \label{PV13} & \mu_{i;n}\ve_{j;n}=\ve_{j;n}\mu_{i;n} \COMMa \eta_{i;n}\ve_{j;n}=\ve_{j;n}\eta_{i;n}, &&\hspace{-3.5cm}\text{if $j\not=i,i+1$.} 
\end{align}

\begin{thm}[cf.~\cite{JE2007b}]\label{t:PVn}
For any $n\in\N$, the partial vine monoid $\PV_n$ has presentation $\pres{X_n}{R_n}$ via~$\phi_n$.  \epfres
\end{thm}

Now let $\Ga\equiv L\cup R\cup X$ be the digraph as in \eqref{e:Ga} and \eqref{e:st}, and let $\phi:\Ga^*\to\PV$ be the morphism in \eqref{e:phi}.  Let $\Om$ be the set of relations over $\Ga$ consisting of $\bigcup_{n\in\N}R_n$ and additionally:
\begin{align}
\label{PV14}\tag*{(PV14)}  \lam_n\rho_n = \iota_n , &&& \rho_n\lam_n = \ve_{n+1;n+1}, \\
\label{PV15}\tag*{(PV15)}  \th_{i;n}\lam_n = \lam_n\th_{i;n+1} , &&& \rho_n\th_{i;n} = \th_{i;n+1}\rho_n, &&\text{for $\th\in\{\si,\si^{-1},\ve,\mu,\eta\}$.}  
\end{align}
The proof of Theorem \ref{t:P1} is easily adapted to give the following:

\begin{thm}\label{t:PV1}
The partial vine category $\PV$ has presentation~$\pres\Ga\Om$ via~$\phi$.  \epfres
\end{thm}

Next, let~$\De$ be the digraph over $\N$, with edges 
\[
X:2\to2 \COMMA X^{-1}:2\to2 \COMMA V:2\to1 \COMMA U:1\to0 \COMMA \U:0\to1.
\]
Define the morphism $\Phi:\De^\oast\to\PV:w\mt\ul w$, where the partial vines $\ul x$ ($x\in\De$) are shown in Figure~\ref{f:PV_Tgens}.  As usual, we also write $I=\io_1$ for the empty path $1\to1$, and $\ul I=\iob_1$ for the identity braid from $\PV_1$ (also shown in Figure~\ref{f:PV_Tgens}).

\begin{figure}[ht]
\begin{center}
\begin{tikzpicture}[scale=.45]
\begin{scope}[shift={(0,0)}]	
\ststring21
\ststring12
\lv1
\lv2
\uv1
\uv2
\draw(0.5,1)node[left]{$\ul X=$};
\end{scope}
\begin{scope}[shift={(7,0)}]	
\ststring12
\ststring21
\lv1
\lv2
\uv1
\uv2
\draw(0.5,1)node[left]{$\ul X^{-1}=$};
\end{scope}
\begin{scope}[shift={(14,0)}]	
\stline11
\stline21
\lv1
\uv1
\uv2
\draw(0.5,1)node[left]{$\ul V=$};
\end{scope}
\begin{scope}[shift={(21,0)}]	
\uv1
\draw(0.5,1)node[left]{$\ul U=$};
\end{scope}
\begin{scope}[shift={(27,0)}]	
\lv1
\draw(0.5,1)node[left]{$\ul \U=$};
\end{scope}
\begin{scope}[shift={(33,0)}]	
\uv1
\lv1
\stline11
\draw(0.5,1)node[left]{$\ul I=$};
\end{scope}
\end{tikzpicture}
\caption{Partial vine generators $\ul x\in\PV$ ($x\in\De$), as well as $\ul I$.}
\label{f:PV_Tgens}
\end{center}
\end{figure}

Let $\Xi$ be the set of the following relations over $\De$:
\begin{gather}
\tag*{(PV1)$'$} \label{PV1'} X\circ X^{-1} = X^{-1}\circ X=\io_2 \COMMA \U\circ U = \io_0  , \\
\tag*{(PV2)$'$} \label{PV2'} X\circ V = V \COMMA V\circ U = U\op U \COMMA (V\op I)\circ V = (I\op V)\circ V \COMMA (I\op \U)\circ V = I , \\
\tag*{(PV3)$'$} \label{PV3'} (X\op I)\circ(I\op X)\circ(X\op I) = (I\op X)\circ(X\op I)\circ(I\op X),\\
\tag*{(PV4)$'$} \label{PV4'} X\circ(U\op I) = I\op U \COMMA X\circ(I\op U) = U\op I,\\
\tag*{(PV5)$'$} \label{PV5'} (\U\op I)\circ X = I\op \U \COMMA  (I\op\U)\circ X = \U\op I,\\
\tag*{(PV6)$'$} \label{PV6'} (I\op V)\circ X = (X\op I)\circ(I\op X)\circ(V\op I) \COMMA (V\op I)\circ X = (I\op X)\circ(X\op I)\circ(I\op V).
\end{gather}

\begin{thm}\label{t:PV2}
The partial vine category $\PV$ has tensor presentation $\pres\De\Xi$ via $\Phi$.
\end{thm}

\pf
As ever, it remains to check Assumption \ref{a:6}.  Again, condition \ref{a61} is established diagrammatically; see Figure \ref{f:P3'} for the third and fourth parts of \ref{PV2'}.

\begin{figure}[ht]
\begin{center}
\begin{tikzpicture}[scale=.45]
\begin{scope}[shift={(-7,2)}]	
\uvs{1,2,3}
\lvs{1,2}
\stline11
\stline21
\stline32
\end{scope}
\begin{scope}[shift={(-7,0)}]	
\uvs{1,2}
\lvs{1}
\stline11
\stline21
\end{scope}
\begin{scope}[shift={(0,1)}]	
\uvs{1,2,3}
\lvs{1}
\stline11
\stline21
\stline31
\draw(-1.5,1)node{$=$};
\end{scope}
\begin{scope}[shift={(7,2)}]	
\uvs{1,2,3}
\lvs{1,2}
\stline11
\stline22
\stline32
\draw(-1.5,0)node{$=$};
\end{scope}
\begin{scope}[shift={(7,0)}]	
\uvs{1,2}
\lvs{1}
\stline11
\stline21
\end{scope}
\begin{scope}[shift={(18,2)}]	
\uvs{1}
\lvs{1,2}
\stline11
\end{scope}
\begin{scope}[shift={(18,0)}]	
\stline11
\stline21
\uvs{1,2}
\lvs{1}
\end{scope}
\begin{scope}[shift={(24,1)}]	
\uvs{1}
\lvs{1}
\stline11
\draw(-1.5,1)node{$=$};
\end{scope}
\end{tikzpicture}
\caption[blah]{The third and fourth parts of relation \ref{PV2'}: $(\ul V\op \ul I)\circ \ul V=(\ul I\op \ul V)\circ \ul V$ and $(\ul I\op \ul \U)\circ \ul V = \ul I$.}
\label{f:P3'}
\end{center}
\end{figure}

For the remaining conditions from Assumption \ref{a:6}, define a morphism $\Ga^*\to\De^\oast:w\mt\wh w$ by
\begin{align*}
\wh\si_{i;n} &= \io_{i-1}\op X\op\io_{n-i-1}, & \wh\mu_{i;n} &= \io_{i-1}\op V\op\U\op\io_{n-i-1}, & \wh\lam_n &= \io_n\op\U,\\
\wh\si_{i;n}^{-1} &= \io_{i-1}\op X^{-1}\op\io_{n-i-1}, & \wh\eta_{i;n} &= \io_{i-1}\op \U\op V\op\io_{n-i-1}, &  \wh\rho_n &= \io_n\op U. \\
\wh\ve_{i;n} &= \io_{i-1}\op U\op\U\op\io_{n-i},
\end{align*}
Condition \ref{a62} of Assumption \ref{a:6} is verified diagrammatically.  Conditions~\ref{a63} and~\ref{a64} are dealt with in Lemmas~\ref{l:a63PV} and~\ref{l:a64PV}.  As usual, we write ${\approx}=\Xi\tsharp$.

\begin{lemma}\label{l:a63PV}
For any $m,n\in\N$, and for any $x\in\De$, we have $x_{m,n} \approx\wh w$ for some $w\in\Ga^*$.
\end{lemma}

\pf
For $x\in\{X,X^{-1},U,\U\}$ the argument is the same as in Lemma \ref{l:a63P}.  For $x=V$, we first use~\ref{PV1'} to show that $\wh\mu_{m+1;m+n+2}\circ U_{m+1;n} \approx V_{m,n}$.  We then apply the $x=U$ case.
\epf

\begin{lemma}\label{l:a64PV}
For any relation $(u,v)\in\Om$, we have $\wh u\approx\wh v$.
\end{lemma}

\pf
The first part of \ref{PV1} follows from the first part of \ref{PV1'}.  The remaining relations from~$\Om$ involving only the $\si_{i;n}$ and $\ve_{j;n}$ are dealt with in the same way as for Lemma \ref{l:a64P}, apart from the third part of \ref{PV7}, which follows (using \eqref{e:UU}, Lemma \ref{l:tc}\ref{tc2} and \ref{PV4'}) from
\[
X\circ X\circ(U\op\U\op I)
= X\circ X\circ(\U\op U\op I)
= \U\op(X\circ X\circ(U\op I))
\approx \U\op(U\op I)
= U\op\U\op I  .
\]
Relations \ref{PV14} and \ref{PV15} are treated in the same way as \ref{P8} and \ref{P9} in Lemma \ref{l:a64P}.

The commuting relations \ref{PV9}, \ref{PV10} and \ref{PV13} follow immediately from the tensor category axioms.  So too does the $|i-j|>1$ case of the first part of \ref{PV12}; for $j=i-1$, we must show that
\[
(I\op V\op\U)\circ(\U\op V\op I) \approx (\U\op V\op I)\circ(I\op V\op\U).
\]
For this, Lemma \ref{l:tc}\ref{tc2} and the tensor category axioms give
\[
(I\op V\op\U)\circ(\U\op V\op I) = ((I\op V)\circ(\U\op V)) \op (\U\circ I) = \U\op((I\op V)\circ V) \op \U.
\]
Similarly, $(\U\op V\op I)\circ(I\op V\op\U) = \U\op((V\op I)\circ V) \op \U$, at which point we apply \ref{PV2'}.

As in the proof of Lemma \ref{l:a64P} (cf.~Table \ref{tab:P}), the remaining relations from $\Om$ boil down to certain reduced relations $s\approx t$.  Up to symmetry, these are shown in Table \ref{tab:PV}.  

\begin{table}[ht]
\begin{center}
\begin{tabular}{|c|c|l|}
\hline
Relation & Label & Reduced relation $s\approx t$  \\
\hline
\ref{PV2} & (a) & $V\op\U \approx (V\op\U)\circ(V\op\U) \approx (\U\op V)\circ(V\op\U) \approx X\circ(V\op\U) \approx (\U\op V)\circ X$ \\
\hline
\ref{PV3} & (b) & $(V\op\U\op I)\circ(I\op V\op\U) \approx (V\op\U\op I)\circ(I\op X)$   \\
               & (c) & $(V\op\U\op I)\circ(I\op\U\op V) \approx V\op\U\op I$ \\
\hline
\ref{PV4} & (d) & $(I\op V\op\U)\circ(V\op\U\op I) \approx (I\op V\op\U)\circ(V\op\U\op I)\circ(I\op V\op\U) $ \\
               &       & $\phantom{(I\op V\op\U)\circ(V\op\U\op I)} \approx (V\op\U\op I)\circ(I\op V\op\U)\circ(V\op\U\op I)$ \\
\hline
\ref{PV5} & (e) & $(V\op\U)\circ(I\op U\op\U) \approx V\op\U$ \\
               & (f) & $(I\op U\op\U)\circ(V\op\U) \approx I\op U\op\U$ \\
\hline
\ref{PV6} & (g) & $(I\op V\op\U)\circ(X\op I) \approx (X\op I)\circ(I\op X)\circ(V\op\U\op I)\circ(I\op V\op\U)$ \\
\hline
\ref{PV8} & (h) & $(U\op\U\op I)\circ(I\op U\op\U) \approx (V\op\U)\circ(U\op\U\op I) \approx (\U\op V)\circ(I\op U\op\U)$ \\
\hline
\end{tabular}
\end{center}
\caption{Reduced relations $s\approx t$ required in the proof of Lemma \ref{l:a64PV}; see the text for more details.}
\label{tab:PV}
\end{table}

Before we begin with the reduced relations, we make three observations.  For the first, Lemma~\ref{l:tc}\ref{tc3} and \ref{PV2'} give
\begin{equation}\label{e:PV4'}
(V\op\U)\circ V = ((V\circ I)\op\U)\circ V = V\circ(I\op\U)\circ V \approx V\circ I = V
\ANDSIM
(\U\op V)\circ V \approx V.
\end{equation}
For the second, \ref{PV2'} and \ref{PV5'} give
\begin{equation}\label{e:PV3'}
(\U\op I)\circ V \approx (\U\op I)\circ X\circ V \approx (I\op\U)\circ V\approx I.
\end{equation}
For the third, Lemma \ref{l:tc}\ref{tc3} and \ref{PV2'} give
\begin{equation}\label{e:PV2'}
(\U\op\U)\circ V = ((\U\circ I)\op\U)\circ V = \U\circ (I\op\U)\circ V \approx \U\circ I = \U.
\end{equation}
\newpage
\pfitem{(a)}  For this we have
\bit
\item $(V\op\U)\circ(V\op\U) = ((V\op\U)\circ V)\op\U \approx V\op\U$, using Lemma \ref{l:tc}\ref{tc3} and \eqref{e:PV4'};
\item $(\U\op V)\circ(V\op\U) = ((\U\op V)\circ V)\op\U \approx V\op\U$, using Lemma \ref{l:tc}\ref{tc3} and \eqref{e:PV4'};
\item $X\circ(V\op\U) = (X\circ V)\op\U \approx V\op\U$, using Lemma \ref{l:tc}\ref{tc3} and \ref{PV2'};
\item $(\U\op V)\circ X = (\U\op (V\circ I))\circ X = V\circ(\U\op I)\circ X \approx V\circ(I\op \U) = (V\circ I)\op \U = V\op\U$, using~\ref{PV5'} and Lemma \ref{l:tc}\ref{tc2} and \ref{tc3}.
\eit
\firstpfitem{(b)}  This follows from $(V\op\U\op I)\circ(I\op X) = (V\circ I) \op ((\U\op I)\circ X) \approx V\op I\op\U$, where we used~\ref{PV5'} in the last step, and
\begin{align*}
(V\op\U\op I)\circ(I\op V\op\U) &= (V\circ I) \op ((\U\op I)\circ(V\op\U)) \\
&=V \op ((\U\op I)\circ V)\op\U &&\text{by Lemma \ref{l:tc}\ref{tc3}}\\
&\approx V \op I\op\U &&\text{by \eqref{e:PV3'}.}
\end{align*}
\firstpfitem{(c)}  Using Lemma \ref{l:tc}\ref{tc2} and \eqref{e:PV3'}, we have
\begin{align*}
(V\op\U\op I)\circ(I\op\U\op V) = (V\circ I)\op((\U\op I)\circ(\U\op V)) = V\op\U\op((\U\op I)\circ V) \approx V\op\U\op I .
\end{align*}
\firstpfitem{(d)}  First note that
\begin{align*}
(I\op V\op\U)\circ(V\op\U\op I) &= ((I\op V)\circ V)\op(\U\circ(\U\op I)) \\
&= ((I\op V)\circ V)\op(\U\op(\U\circ I)) &&\text{by Lemma \ref{l:tc}\ref{tc2}}\\
&\approx ((V\op I)\circ V)\op\U\op\U &&\text{by \ref{PV2'}.}
\end{align*}
Thus, it suffices to show that $(V\op\U\op I)\circ w\approx w \approx w\circ(I\op V\op\U)$, where $w=((V\op I)\circ V)\op\U\op\U$.  For this we have
\begin{align*}
(V\op\U\op I)\circ w &= (V\op\U\op I)\circ (((V\op I)\circ V)\op\U\op\U) \\
&= ((V\op\U\op I)\circ (V\op I)\circ V)\op\U\op\U &&\text{by Lemma \ref{l:tc}\ref{tc3}}\\
&= (((V\op\U)\circ V)\op (I\circ I))\circ V)\op\U\op\U \\
&\approx ((V\op I)\circ V)\op\U\op\U = w  &&\text{by \eqref{e:PV4'},}
\end{align*}
and
\begin{align*}
w\circ(I\op V\op\U) &= (((V\op I)\circ V)\op\U\op\U)\circ(I\op V\op\U) \\
&= (((V\op I)\circ V)\circ I) \op ((\U\op\U)\circ(V\op\U)) \\
&= ((V\op I)\circ V) \op ((\U\op\U)\circ V)\op\U &&\text{by Lemma \ref{l:tc}\ref{tc3}}\\
&\approx ((V\op I)\circ V) \op \U\op\U = w &&\text{by \eqref{e:PV2'}.}
\end{align*}
\firstpfitem{(e)}  We use \eqref{e:UU} and~\ref{PV1'}:
\[
(V\op\U)\circ(I\op U\op\U) = (V\circ I)\op(\U\circ(U\op\U)) = V\op(\U\circ U\circ\U) \approx V\op(\io_0\circ\U) = V\op\U.
\]
\firstpfitem{(f)}  We use \eqref{e:UU} and \ref{PV2'}:
\[
(I\op U\op\U)\circ(V\op\U) = (I\op \U\op U)\circ(V\op\U) = ((I\op \U)\circ V)\op (U\circ\U) \approx I\op U\op\U.
\]
\firstpfitem{(g)}  First note that
\begin{align*}
(I\op V\op\U)\circ(X\op I) &= ((I\op V)\circ X)\op(\U\circ  I) \\
&\approx ((X\op I)\circ(I\op X)\circ(V\op I)) \op \U &&\text{by \ref{PV6'}}\\
&= (X\op I)\circ(I\op X)\circ(V\op I \op \U) &&\text{by Lemma \ref{l:tc}\ref{tc3},}
\end{align*}
so it remains to observe (using Lemma \ref{l:tc}\ref{tc3} and \eqref{e:PV3'}) that
\begin{align*}
(V\op\U\op I)\circ(I\op V\op\U) &= (V\circ I)\op((\U\op I)\circ(V\op\U)) = V\op((\U\op I)\circ V)\op\U \approx V\op I\op\U.
\end{align*}
\firstpfitem{(h)}  By \eqref{e:UUUU}, the proof of which used only $U\op\U=\U\op U$,  it suffices to show that
\[
(V\op\U)\circ(U\op\U\op I) \approx (\U\op V)\circ(I\op U\op\U) \approx U\op U\op\U\op\U.
\]
For this we have (using \ref{PV2'} and Lemma \ref{l:tc}\ref{tc2})
\begin{align*}
(V\op\U)\circ(U\op\U\op I) = (V\circ U)\op(\U\circ(\U\op I)) \approx (U\op U)\op(\U\op(\U\circ I)) = U\op U\op\U\op\U .
\end{align*}
A similar calculation gives $(\U\op V)\circ(I\op U\op\U) \approx U\op U\op\U\op\U$.
\epf

\noindent This completes the proof of the theorem.
\epf

\begin{rem}\label{r:PV}
As in Remark \ref{r:*PROP}, we could convert the presentation $\pres\De\Xi$ from Theorem \ref{t:PV2} into a PROB presentation.  In doing this, the generators $X$ and $X^{-1}$, and some of the relations from $\Xi$, are part of the ``free'' PROB data, and can hence be removed, specifically the first part of \ref{PV1'} and all of \ref{PV3'}--\ref{PV6'}.  
\end{rem}

\subsection{The partial braid category}\label{ss:IB}

Next we treat the partial braid category $\IB$.  Again Assumption \ref{a:1} holds with $d=1$, and for Assumption~\ref{a:2} we take the partial braids $\lamb_n$ and $\rhob_n$ pictured in Figure \ref{f:PV_gens}.  Several presentations for the inverse braid monoids $\IB_n$ ($n\in\N$) exist \cite{EL2004,JE2007a,Gilbert2006}, but the most convenient one to use for Assumption~\ref{a:3} is due to Gilbert \cite{Gilbert2006}.  For $n\in\N$, let 
\[
X_n = S_n\cup S_n^{-1}\cup E_n \qquad\text{where}\qquad S_n^{\pm1}=\set{\si_{i;n}^{\pm1}}{1\leq i<n} \AND E_n = \set{\ve_{i;n}}{1\leq i\leq n},
\]
and define $\phi_n:X_n^*\to\IB_n:w\mt\ol w$, where the partial braids $\ol\si_{i;n}^{\pm1}$ and $\ol\ve_{i;n}$ are pictured in Figure \ref{f:PV_gens}.  Let $R_n$ be the set of relations
\begin{align}
\tag*{(IB1)} \label{IB1} & \si_{i;n}\si_{i;n}^{-1} = \si_{i;n}^{-1}\si_{i;n} = \io_n \COMMa \ve_{i;n}^2=\ve_{i;n} \COMMa \ve_{i;n}\ve_{j;n}=\ve_{j;n}\ve_{i;n} \COMMa
\si_{i;n}\ve_{j;n}=\ve_{j;n}\si_{i;n} \text{ if $j\not=i,i+1$,} \\
\tag*{(IB2)} \label{IB2} & \si_{i;n}\ve_{i;n}=\ve_{i+1;n}\si_{i;n} \COMMa \si_{i;n}\ve_{i+1;n}=\ve_{i;n}\si_{i;n} \COMMa \si_{i;n}^2\ve_{i;n}=\ve_{i;n} \COMMa \si_{i;n}\ve_{i;n}\ve_{i+1;n}=\ve_{i;n}\ve_{i+1;n},  \\
\tag*{(IB3)} \label{IB3} & \si_{i;n}\si_{j;n} = \si_{j;n}\si_{i;n}  \text{ if $|i-j|>1$}\COMMa
\si_{i;n}\si_{j;n}\si_{i;n} = \si_{j;n}\si_{i;n}\si_{j;n} \text{ if $|i-j|=1$.} 
\end{align}

\begin{thm}[cf.~\cite{Gilbert2006}]\label{t:IBn}
For any $n\in\N$, the inverse braid monoid $\IB_n$ has presentation $\pres{X_n}{R_n}$ via~$\phi_n$.  \epfres
\end{thm}

Now let $\Ga\equiv L\cup R\cup X$ be the digraph as in \eqref{e:Ga} and \eqref{e:st}, and let $\phi:\Ga^*\to\IB$ be the morphism in \eqref{e:phi}.  Let $\Om$ be the set of relations over $\Ga$ consisting of $\bigcup_{n\in\N}R_n$ and additionally:
\begin{align}
\label{IB4}\tag*{(IB4)}  \lam_n\rho_n = \iota_n , &&& \rho_n\lam_n = \ve_{n+1;n+1}, \\
\label{IB5}\tag*{(IB5)}  \th_{i;n}\lam_n = \lam_n\th_{i;n+1} , &&& \rho_n\th_{i;n} = \th_{i;n+1}\rho_n, &&\text{for $\th\in\{\si,\si^{-1},\ve\}$.}  
\end{align}
The proof of the following is again essentially the same as for Theorem \ref{t:P1}:

\begin{thm}\label{t:IB1}
The partial braid category $\IB$ has presentation~$\pres\Ga\Om$ via~$\phi$.  \epfres
\end{thm}

Next, let~$\De$ be the digraph over $\N$, with edges 
\[
X:2\to2 \COMMA X^{-1}:2\to2 \COMMA U:1\to0 \COMMA \U:0\to1,
\]
and define the morphism $\Phi:\De^\oast\to\IB:w\mt\ul w$, where the partial braids $\ul x$ ($x\in\De$) are pictured in Figure \ref{f:PV_Tgens}.  Let $\Xi$ be the set of the following relations over $\De$, where $I=\io_1$:
\begin{gather}
\tag*{(IB1)$'$} \label{IB1'} X\circ X^{-1} = X^{-1}\circ X=\io_2 \COMMA \U\circ U = \io_0 , \\
\tag*{(IB2)$'$} \label{IB2'} (X\op I)\circ(I\op X)\circ(X\op I) = (I\op X)\circ(X\op I)\circ(I\op X),\\
\tag*{(IB3)$'$} \label{IB3'} X\circ(U\op I) = I\op U \COMMa X\circ(I\op U) = U\op I\COMMa (\U\op I)\circ X = I\op \U \COMMa  (I\op\U)\circ X = \U\op I.
\end{gather}
The proof of the following is contained in the proof of Theorem \ref{t:PV2}:

\begin{thm}\label{t:IB2}
The partial braid category $\IB$ has tensor presentation $\pres\De\Xi$ via $\Phi$.  \epfres
\end{thm}

\begin{rem}\label{r:IB2}
For a PROB presentation, only the relation $\U\circ U = \io_0$ is needed (cf.~Remark \ref{r:PV}).
\end{rem}

\subsection{The vine category}\label{ss:V}

We now turn our attention to the (full) vine category $\V$.  Things are a little more complicated here since the connectivity of $\V$ is different to all of the other categories treated so far.  Namely, as we have already observed, $\V_{m,n}\not=\emptyset \iff m=0$ or $n\geq1$.  Thus, we will have to use Theorem \ref{t:pres3} instead of Theorem \ref{t:pres2}.  Accordingly, we define the subcategory
\[
\V^+=\bigcup_{m,n\in\PP}\V_{m,n}.
\]
Assumption \ref{a:1} holds in $\V^+$ with $d=1$.  Next, note that the partial vines $\rhob_n$ used in Sections \ref{ss:PV} and~\ref{ss:IB} (cf.~Figure \ref{f:PV_gens}) do not belong to $\V$.  Thus, for Assumption \ref{a:2}, we take $\lamb_n\in\V_{n,n+1}$ and $\rhob_n\in\V_{n+1,n}$ ($n\in\PP$) to be the (full) vines pictured in Figure \ref{f:V_gens}.  For Assumption \ref{a:3} we take the presentation for the vine monoids $\V_n$ ($n\in\PP$) given by Lavers \cite{Lavers1997}.  For $n\in\PP$, define the alphabet
\begin{align*}
X_n = S_n\cup S_n^{-1}\cup M_n\cup H_n \qquad\text{where}\qquad 
S_n &= \set{\si_{i;n}}{1\leq i<n} , & M_n &= \set{\mu_{i;n}}{1\leq i< n}, \\
S_n^{-1} &= \set{\si_{i;n}^{-1}}{1\leq i<n} , & H_n &= \set{\eta_{i;n}}{1\leq i< n}. 
\end{align*}
Define the morphism $\phi_n:X_n^*\to\V_n:w\mt\ol w$, where the vines $\ol x$ ($x\in X_n$) are also shown in Figure~\ref{f:V_gens}.  Let $R_n$ be the following set of relations over $X_n$:
\begin{align}
\tag*{(V1)} \label{V1} & \si_{i;n}\si_{i;n}^{-1} = \si_{i;n}^{-1}\si_{i;n} = \io_n,  \\
\tag*{(V2)} \label{V2} & \mu_{i;n} = \mu_{i;n}^2 = \eta_{i;n}\mu_{i;n} = \si_{i;n}\mu_{i;n} = \eta_{i;n}\si_{i;n} \COMMa \eta_{i;n} = \eta_{i;n}^2 = \mu_{i;n}\eta_{i;n} = \si_{i;n}\eta_{i;n}  = \mu_{i;n}\si_{i;n}, \\
\tag*{(V3)} \label{V3} & \mu_{i;n}\mu_{i+1;n}=\mu_{i;n}\si_{i+1;n} \COMMa \eta_{i+1;n}\eta_{i;n}=\eta_{i+1;n}\si_{i;n} \COMMa \mu_{i;n}\eta_{i+1;n} = \mu_{i;n} \COMMa \eta_{i+1;n}\mu_{i;n} = \eta_{i+1;n}, \\
\tag*{(V4)} \label{V4} & \mu_{i+1;n}\mu_{i;n} = \mu_{i;n}\mu_{i+1;n}\mu_{i;n} = \mu_{i+1;n}\mu_{i;n}\mu_{i+1;n} \COMMa \eta_{i;n}\eta_{i+1;n} = \eta_{i;n}\eta_{i+1;n}\eta_{i;n} = \eta_{i+1;n}\eta_{i;n}\eta_{i+1;n}, \\
\tag*{(V5)} \label{V5} & \mu_{i+1;n}\si_{i;n}=\si_{i;n}\si_{i+1;n}\mu_{i;n}\mu_{i+1;n} \COMMa \eta_{i;n}\si_{i+1;n}=\si_{i+1;n}\si_{i;n}\eta_{i+1;n}\eta_{i;n}, \\
\tag*{(V6)} \label{V6} & \si_{i;n}\si_{j;n} = \si_{j;n}\si_{i;n} \COMMa \mu_{i;n}\mu_{j;n} = \mu_{j;n}\mu_{i;n} \COMMa \eta_{i;n}\eta_{j;n} = \eta_{j;n}\eta_{i;n}, &&\hspace{-3.5cm}\text{if $|i-j|>1$,} \\
\tag*{(V7)} \label{V7} & \si_{i;n}\mu_{j;n}=\mu_{j;n}\si_{i;n} \COMMa \si_{i;n}\eta_{j;n}=\eta_{j;n}\si_{i;n}, &&\hspace{-3.5cm}\text{if $|i-j|>1$,} \\
\tag*{(V8)} \label{V8} & \si_{i;n}\si_{j;n}\si_{i;n} = \si_{j;n}\si_{i;n}\si_{j;n} ,&&\hspace{-3.5cm}\text{if $|i-j|=1$,} \\
\tag*{(V9)} \label{V9} & \mu_{i;n}\eta_{j;n} = \eta_{j;n}\mu_{i;n},  &&\hspace{-3.5cm}\text{if $j\not=i,i+1$.} 
\end{align}

\begin{thm}[cf.~\cite{Lavers1997}]\label{t:Vn}
For any $n\in\PP$, the vine monoid $\V_n$ has presentation $\pres{X_n}{R_n}$ via~$\phi_n$.  \epfres
\end{thm}

\begin{figure}[ht]
\begin{center}
\begin{tikzpicture}[scale=.42]
\begin{scope}[shift={(0,0)}]	
\udotted14
\ddotted14
\udotted7{10}
\ddotted7{10}
\stline11
\stline44
\ststring65
\ststring56
\stline77
\stline{10}{10}
\draw(0.5,1)node[left]{$\sib_{i;n}=$};
\node()at(1,2.5){\tiny$1$};
\node()at(5,2.5){\tiny$i$};
\node()at(10,2.5){\tiny$n$};
\uvs{1,4,5,6,7,10}
\lvs{1,4,5,6,7,10}
\end{scope}
\begin{scope}[shift={(0,-5)}]	
\udotted14
\ddotted14
\udotted7{10}
\ddotted7{10}
\stline11
\stline44
\ststring56
\ststring65
\stline77
\stline{10}{10}
\draw(0.5,1)node[left]{$\sib_{i;n}^{-1}=$};
\node()at(1,2.5){\tiny$1$};
\node()at(5,2.5){\tiny$i$};
\node()at(10,2.5){\tiny$n$};
\uvs{1,4,5,6,7,10}
\lvs{1,4,5,6,7,10}
\end{scope}
\begin{scope}[shift={(15,0)}]	
\udotted14
\ddotted14
\udotted7{10}
\ddotted7{10}
\stline11
\stline44
\stline55
\stline65
\stline77
\stline{10}{10}
\draw(0.5,1)node[left]{$\mub_{i;n}=$};
\node()at(1,2.5){\tiny$1$};
\node()at(5,2.5){\tiny$i$};
\node()at(10,2.5){\tiny$n$};
\uvs{1,4,5,6,7,10}
\lvs{1,4,5,6,7,10}
\end{scope}
\begin{scope}[shift={(15,-5)}]	
\udotted14
\ddotted14
\udotted7{10}
\ddotted7{10}
\stline11
\stline44
\stline56
\stline66
\stline77
\stline{10}{10}
\draw(0.5,1)node[left]{$\etab_{i;n}=$};
\node()at(1,2.5){\tiny$1$};
\node()at(5,2.5){\tiny$i$};
\node()at(10,2.5){\tiny$n$};
\uvs{1,4,5,6,7,10}
\lvs{1,4,5,6,7,10}
\end{scope}
\begin{scope}[shift={(30,0)}]	
\uvs{1,6}
\lvs{1,6,7}
\udotted16
\ddotted16
\stline11
\stline66
\draw(0.5,1)node[left]{$\lamb_n=$};
\node()at(1,2.5){\tiny$1$};
\node()at(6,2.5){\tiny$n$};
\end{scope}
\begin{scope}[shift={(30,-5)}]	
\uvs{1,6,7}
\lvs{1,6}
\udotted16
\ddotted16
\stline11
\stline66
\stline76
\draw(0.5,1)node[left]{$\rhob_n=$};
\node()at(1,2.5){\tiny$1$};
\node()at(6,2.5){\tiny$n$};
\end{scope}
\end{tikzpicture}
\caption{Vine generators $\ol x\in\V$ ($x\in\Ga$).}
\label{f:V_gens}
\end{center}
\end{figure}

Now let $\Ga\equiv L\cup R\cup X$ be the digraph as in \eqref{e:Ga} and \eqref{e:st}, and let $\phi:\Ga^*\to\V^+$ be the morphism in~\eqref{e:phi}.  Because of the different form of $\rhob_n$ ($n\in\PP$), we have to be a little more careful in constructing the relations from Assumption \ref{a:4}\ref{a44}.  With this in mind, let $\Om$ be the set of relations over~$\Ga$ consisting of $\bigcup_{n\in\PP}R_n$ and additionally:
\begin{align}
\label{V10}\tag*{(V10)}  \lam_n\rho_n = \iota_n , &&& \rho_n\lam_n = \mu_{n;n+1}, \\
\label{V11}\tag*{(V11)}  &&& \th_{i;n}\lam_n = \lam_n\th_{i;n+1}, &&\text{for $\th\in\{\si,\si^{-1},\mu,\eta\}$,}\\
\label{V12}\tag*{(V12)}  &&& \rho_n\th_{i;n} = \th_{i;n+1}\rho_n, &&\text{for $\th\in\{\si,\si^{-1},\mu,\eta\}$ and $i\leq n-2$,}\\
\label{V13}\tag*{(V13)}  &&& \rho_n\th_{n-1;n} = \mu_{n;n+1}\th_{n-1;n+1}\rho_n, &&\text{for $\th\in\{\si,\si^{-1},\mu,\eta\}$.}
\end{align}
In the language of Assumption \ref{a:4}\ref{a44}, the mapping $x\mt x_+$ is still given by $(\th_{i;n})_+=\th_{i;n+1}$, but the $x\mt x^+$ mapping is slightly more complicated.

\begin{thm}\label{t:V1}
The vine category $\V^+$ has presentation~$\pres\Ga\Om$ via~$\phi$.  
\end{thm}

\pf
As usual, the main work is in checking condition \ref{a45} of Assumption \ref{a:4}, and for this we again use Remark \ref{r:a45}.  Writing $\mu=\mu_{n,n+1}$, the proof will be complete if we can show the following:
\bit
\item For all $n\in\PP$ and $w\in X_{n+1}^*$,  we have  $\mu w\mu  \sim \mu u_+\mu $  for some $u\in X_n^*$.
\eit
To do so, let $w\in X_{n+1}^*$, and put $\al=\ol{\mu w\mu}\in\V_{n+1}$.  Also, put $\be=\ol{w\mu}$, and note that $\al=\mub\be\mub$ (since~$\mub^2=\mub$).  Note also that $T(\be)\sub T(\mub)=[n]$, where $T$ is defined in Section \ref{ss:pre_braids}.  Thus, removing string $n+1$ of $\be$ leaves us with a vine from $\V_n$; denote this vine by $\ga$, and write $\ga_+=\ga\op\iob_1\in\V_{n+1}$.  Since $T(\mub)=[n]$, we have $\mub\be=\mub\ga_+$.  By Theorem \ref{t:Vn}, we have $\ga=\ol u$ for some $u\in X_n^*$, and we note that $\ga_+=\ol u_+ = \ol{u_+}$.  But then $\ol{\mu w\mu} = \al = \mub\be\mub = \mub\ga_+\mub = \ol{\mu u_+\mu}$ in $\V_{n+1}$, so it follows from Theorem~\ref{t:Vn} (and $R_{n+1}\sub\Om$)  that $\mu w\mu \sim \mu u_+\mu$.
\epf

Now let~$\De$ be the digraph over $\N$, with edges 
\[
X:2\to2 \COMMA X^{-1}:2\to2 \COMMA V:2\to1 \COMMA \U:0\to1,
\]
and define the morphism $\Phi:\De^\oast\to\V:w\mt\ul w$, where the partial braids $\ul x$ ($x\in\De$) are pictured in Figure \ref{f:PV_Tgens}.  Let $\Xi$ be the set of the following relations over $\De$, where as usual we write~$I=\io_1$:
\begin{gather}
\tag*{(V1)$'$} \label{V1'} X\circ X^{-1} = X^{-1}\circ X=\io_2 \COMMA X\circ V = V , \\
\tag*{(V2)$'$} \label{V2'} (V\op I)\circ V = (I\op V)\circ V \COMMA (I\op \U)\circ V = I , \\
\tag*{(V3)$'$} \label{V3'} (X\op I)\circ(I\op X)\circ(X\op I) = (I\op X)\circ(X\op I)\circ(I\op X),\\
\tag*{(V4)$'$} \label{V4'} (\U\op I)\circ X = I\op \U \COMMA  (I\op\U)\circ X = \U\op I,\\
\tag*{(V5)$'$} \label{V5'} (I\op V)\circ X = (X\op I)\circ(I\op X)\circ(V\op I) \COMMA (V\op I)\circ X = (I\op X)\circ(X\op I)\circ(I\op V).
\end{gather}

\begin{thm}\label{t:V2}
The vine category $\V$ has tensor presentation $\pres\De\Xi$ via $\Phi$.  
\end{thm}

\pf
We prove this by applying Theorem \ref{t:pres3}.  We have already observed that the subcategory~$\V^+$ satisfies Assumptions \ref{a:1}--\ref{a:4}; we fix the presentation $\pres\Ga\Om$ for $\V^+$ from Theorem \ref{t:V1}.  Since $\V$ satisfies Assumption \ref{a:7}, it remains to check Assumption \ref{a:8}.  Conditions \ref{a80} and \ref{a81} are easily verified.  For conditions \ref{a82}--\ref{a84} we use the mapping $\Ga^*\to\De^\oast:w\mt\wh w$ defined by
\begin{align*}
\wh\si_{i;n} &= \io_{i-1}\op X\op\io_{n-i-1}, & \wh\mu_{i;n} &= \io_{i-1}\op V\op\U\op\io_{n-i-1}, & \wh\lam_n &= \io_n\op\U,\\
\wh\si_{i;n}^{-1} &= \io_{i-1}\op X^{-1}\op\io_{n-i-1}, & \wh\eta_{i;n} &= \io_{i-1}\op \U\op V\op\io_{n-i-1}, &  \wh\rho_n &= \io_{n-1}\op V. 
\end{align*}
Condition \ref{a82} is easily checked.  Condition \ref{a83} is given by the next lemma.  From here on we write as usual ${\approx}=\Xi\tsharp$.  

\begin{lemma}\label{l:a73V}
For any $m,n\in\N$, and for any $x\in\De$, if $(m,x,n)\not=(0,\U,0)$ then $x_{m,n} \approx\wh w$ for some $w\in\Ga^*$.
\end{lemma}

\pf
For $x\in\{X,X^{-1},\U\}$ the argument is again essentially the same as in Lemma \ref{l:a63P}.  For $x=V$, we have
\[
V_{m,n} \approx \wh\mu_{m+1;m+n+2}\circ\wh\mu_{m+2;m+n+2}\circ\cdots\circ\wh\mu_{m+n;m+n+2}\circ\wh\rho_{m+n+1}.
\]
To prove this, we first show that $\wh\mu_{m+1;m+n+2}\circ V_{m+1;n-1} \approx V_{m,n}$ and then apply induction.
\epf

Condition \ref{a84} of Assumption \ref{a:8} is given by the following:  

\begin{lemma}\label{l:a74V}
For any relation $(u,v)\in\Om$, we have $\wh u\approx\wh v$.
\end{lemma}

\pf
The proof is essentially contained in the proof of Lemma \ref{l:a64PV}, apart from relations \ref{V10}, \ref{V12} and~\ref{V13}, which involve the different $\wh\rho_n$ terms.  The first two of these are easily verified, and the third boils down to showing that
\[
(I\op V)\circ Z \approx (I\op V\op\U)\circ(Z\op I)\circ(I\op V) \qquad\text{for each $Z=X$, $X^{-1}$, $V\op\U$ and $\U\op V$.}
\]
For this, noting that $\bd(Z)=\br(Z)=2$, we have
\begin{align*}
(I\op V\op\U)\circ(Z\op I)\circ(I\op V) &= (((I\op V)\circ Z)\op(\U\circ I))\circ(I\op V) \\
&= (((I\op V)\circ Z)\op\U)\circ(I\op V) \\
&= (((I\op V)\circ Z\circ(I\op I))\op\U)\circ(I\op V) \\
&= (I\op V)\circ Z\circ(I\op I\op\U)\circ(I\op V) &&\text{by Lemma \ref{l:tc}\ref{tc3}}\\
&= (I\op V)\circ Z\circ((I\circ I)\op ((I\op\U)\circ V)) \\
&\approx (I\op V)\circ Z\circ(I\op I)= (I\op V)\circ Z &&\text{by \ref{V2'}.}  \qedhere
\end{align*}
\epf

All that remains is to check condition \ref{a85} of Assumption \ref{a:8}, and for this we need to show that
\[
(\U\op\U)\circ X \approx \U\op\U \COMMA (\U\op\U)\circ X^{-1} \approx \U\op\U \COMMA (\U\op\U)\circ V\approx \U.
\]
The first and third of these were proved in \eqref{e:UUX} and \eqref{e:PV2'}, using relations contained in \ref{V1'}--\ref{V5'}.  The second follows quickly from the first, together with \ref{V1'}.  
\epf

\subsection{Categories of (partial) transformations}\label{ss:T}

We now show how to use the results of Sections \ref{ss:PV}--\ref{ss:V} to obtain presentations for certain categories of transformations.  Such presentations could also be obtained directly by applying the results of Section~\ref{s:C}, but it is quicker to realise the transformation categories as suitable quotients of $\PV$, $\V$ or $\IB$.

For $m,n\in\N$, write $\PT_{m,n}$ for the set of all partial transformations $[m]\to[n]$: i.e., all functions $A\to [n]$ for $A\sub[m]$.  The partial transformation category is
\[
\PT = \bigset{(m,f,n)}{m,n\in\N,\ f\in\PT_{m,n}}.
\]
For $m,n,q\in\N$, and for $f\in\PT_{m,n}$ and $g\in\PT_{n,q}$, we define
\[
\bd(m,f,n)=m \COMMA \br(m,f,n)=n \COMMA (m,f,n)\circ(n,g,q) = (m,fg,q),
\]
where $fg=f\circ g\in\PT_{m,q}$ is the ordinary relational composition.  As with partial vines, we will avoid clutter by identifying $(m,f,n)\in\PT$ with $f\in\PT_{m,n}$, and writing $\bd(f)=m$ and $\br(f)=n$.

There is a natural surmorphism
\[
\psi:\PV\to\PT:\al\mt\wt\al
\]
determined by the initial and terminal points of the strings of the partial vine $\al$.  For $\al\in\PV_{m,n}$, the partial transformation $\wt\al\in\PT_{m,n}$ has domain $I(\al)$ and image $T(\al)$, and maps $a\mt b$ if and only if $\al$ has a string $\s$ with $I(\s)=a$ and $T(\s)=b$.  For example, with $\al\in\PV_{4,5}$ and $\be\in\PV_{5,3}$ as in Figure~\ref{f:vines}, we have $\wt\al = \left(\begin{smallmatrix}1&2&3&4\\2&4&2&1\end{smallmatrix}\right)\in\PT_{4,5}$ and $\wt\be = \left(\begin{smallmatrix}1&2&3&4&5\\-&3&3&3&2\end{smallmatrix}\right)\in\PT_{5,3}$.

The surmorphism $\psi$ maps the vine category $\V$ and the partial braid category $\IB$ onto the subcategories
\[
\T = \set{f\in\PT}{f\text{ is totally defined}} \AND \I = \set{f\in\PT}{f\text{ is injective}},
\]
which are the full transformation category and the symmetric inverse category, respectively.  Thus, presentations for $\PV$, $\V$ or $\IB$ yield presentations for $\PT$, $\T$ or $\I$ upon adding relations that generate the kernel of $\psi$ (or its restriction to $\V$ or $\IB$).

\begin{lemma}\label{l:ker}
If $\C$ is any of $\PV$, $\V$ or $\IB$, then
\ben
\item \label{ker1} $\ker(\psi\restr_\C)$ is generated as a category congruence by $\bigset{(\ol\si_{i;n}^2,\iob_n)}{n\in\N,\ 1\leq i<n}$,
\item \label{ker2} $\ker(\psi\restr_\C)$ is generated as a tensor category congruence by $\big\{(\ul X\circ\ul X,\iob_2)\big\}$,
\een
\end{lemma}

\pf
We must show that $\ker(\psi\restr_\C)=\xi=\ze$, where
\[
\xi = \bigset{(\ol\si_{i;n}^2,\iob_n)}{n\in\N,\ 1\leq i<n}^\sharp \AND \ze = \big\{(\ul X\circ\ul X,\iob_2)\big\}\tsharp.
\]
Since $\ol\si_{i;n}^2 = \iob_{i-1}\circ(\ul X\circ\ul X)\op\iob_{n-i-2}$, we have $\xi\sub\ze$.  Since $\ul X\psi$ is the transposition $(1,2)\in\SS_2$, we have~${\ze\sub\ker(\psi\restr_\C)}$.  It remains to show that $\ker(\psi\restr_\C)\sub\xi$, so suppose $\al,\be\in\C$ and $\wt\al=\wt\be$.  Write $m=\bd(\al)=\bd(\be)$ and $n=\br(\al)=\br(\be)$, and suppose by symmetry that $m\leq n$.  Let $\lamb_{m,n},\rhob_{n,m}$ be as defined in \eqref{e:lambrhob}, where the latter involves the appropriate $\rhob_i$ (cf.~Figures \ref{f:PV_gens} and \ref{f:V_gens}).  We then have $(\rhob_{n,m}\al)\psi = (\rhob_{n,m}\be)\psi$, with $\rhob_{n,m}\al,\rhob_{n,m}\be\in\C_n$.  Now, the kernel of $\psi\restr_{\C_n}$ is generated (as a semigroup congruence) by $\bigset{(\ol\si_{i;n}^2,\iob_2)}{1\leq i<n}$; see \cite[Lemma~29]{JE2007b}, \cite[Theorem 9]{Lavers1997} and \cite[Proposition 31]{JE2007a} for $\C=\PV$, $\V$ and $\IB$, respectively.  It follows that $\rhob_{n,m}\al \mathrel\xi \rhob_{n,m}\be$, and so $\al = \io_m\al = \lamb_{m,n}\rhob_{n,m}\al \mathrel\xi \lamb_{m,n}\rhob_{n,m}\be = \io_m\be = \be$, as required.
\epf

Thus, we can obtain a category presentation for $\PT$,~$\I$ or~$\T$ by adjoining the relations $\si_{i;n}^2 = \io_n$ to the appropriate presentation $\pres\Ga\Om$ from Theorem \ref{t:PV1}, \ref{t:IB1} or \ref{t:V1}.  Combining this new relation with the first part of \ref{PV1}, this leads to $\si_{i;n}^{-1}=\si_{i;n}$, so that we may remove all generators from $\bigcup_{n\in\N}S_n^{-1}$.  Several relations simplify as a result, but we will not give all the details, as we are more interested in tensor presentations.

We can obtain such a tensor presentation for $\PT$ by adding the relation $X\circ X=\io_2$ to the presentation $\pres\De\Xi$ from Theorem \ref{t:PV2}.  Keeping in mind that this also leads to $X^{-1}=X$, we end up with generators $X$, $V$, $U$ and $\U$, and relations
\begin{gather*}
X\circ X=\io_2 \COMMA (X\op I)\circ(I\op X)\circ(X\op I) = (I\op X)\circ(X\op I)\circ(I\op X), \\
\U\circ U = \io_0  \COMMA X\circ V = V \COMMA V\circ U = U\op U \COMMA (V\op I)\circ V = (I\op V)\circ V \COMMA (I\op \U)\circ V = I , \\
X\circ(U\op I) = I\op U \COMMA (\U\op I)\circ X = I\op \U \COMMA  (I\op V)\circ X = (X\op I)\circ(I\op X)\circ(V\op I) .
\end{gather*}
Note that only the first part of \ref{PV4'} is listed, as the second follows from the first and $X\circ X=\io_2$: \emph{viz.}, ${X\circ(I\op U) \approx X\circ X\circ(U\op I) \approx \io_2\circ(U\op I)=U\op I}$.  Similarly, we only need the first parts of~\ref{PV5'} and \ref{PV6'}.  The above presentation is via the surmorphism determined by
\begin{equation}\label{e:trans}
X\mt\left(\begin{smallmatrix}1&2\\2&1\end{smallmatrix}\right)\in\PT_2 \COMMA
V\mt\left(\begin{smallmatrix}1&2\\1&1\end{smallmatrix}\right)\in\PT_{2,1} \COMMA
U\mt\left(\begin{smallmatrix}1\\-\end{smallmatrix}\right)\in\PT_{1,0} \COMMA
\U\mt(\emptyset)\in\PT_{0,1}.
\end{equation}
We can similarly obtain tensor presentations for $\I$ and $\T$, using Theorems \ref{t:IB2} and \ref{t:V2}:
\bit
\item For $\I$ we have generators $X$, $U$ and $\U$, mapping as in \eqref{e:trans}, and relations
\begin{gather*}
X\circ X=\io_2 \COMMA (X\op I)\circ(I\op X)\circ(X\op I) = (I\op X)\circ(X\op I)\circ(I\op X),\\
\U\circ U = \io_0  \COMMA X\circ(U\op I) = I\op U \COMMA (\U\op I)\circ X = I\op \U  .
\end{gather*}
\item For $\T$ we have generators $X$, $V$ and $\U$, mapping as in \eqref{e:trans}, and relations
\begin{gather*}
X\circ X=\io_2 \COMMA (X\op I)\circ(I\op X)\circ(X\op I) = (I\op X)\circ(X\op I)\circ(I\op X), \\
X\circ V = V \COMMA (V\op I)\circ V = (I\op V)\circ V \COMMA (I\op \U)\circ V = I ,\\
(\U\op I)\circ X = I\op \U \COMMA  (I\op V)\circ X = (X\op I)\circ(I\op X)\circ(V\op I) .
\end{gather*}
\eit

\begin{rem}
All of the above presentations can be simplified by adding in the PROP structure (cf.~Remarks \ref{r:*PROP}, \ref{r:PV} and \ref{r:IB2}).  For example, the category $\I$ then requires only the relation $\U\circ U=\io_0$.
\end{rem}

\subsection{Categories of isotone (partial) transformations}\label{ss:O}

A partial transformation $f\in\PT_{m,n}$ is isotone (order-preserving) if we have $x\leq y\implies xf\leq yf$ for all $x,y\in\dom(f)$.  The set
\[
\PO = \set{f\in\PT}{f \text{ is isotone}}
\]
is a subcategory of $\PT$.  We also have the isotone subcategories of $\T$ and $\I$:
\[
\O = \PO\cap\T \AND \OI = \PO\cap\I.
\]
The results of Section \ref{s:C} can be used to obtain presentations for the categories $\PO$, $\O$ and $\OI$.

For $\PO$ and $\O$, we first note that:
\bit
\item For fixed $n\in\N$, the monoid $\PO_n$ has presentation with generators $E_n\cup M_n\cup H_n$ as in \eqref{e:PVn_gens}, and all relations from \ref{PV1}--\ref{PV13} involving no letters from $S_n\cup S_n^{-1}$; cf.~\cite{Popova1962}.
\item For fixed $n\in\PP$, the monoid $\O_n$ has presentation with generators $M_n\cup H_n$ as in \eqref{e:PVn_gens}, and all relations from \ref{PV1}--\ref{PV13} involving no letters from $E_n\cup S_n\cup S_n^{-1}$; cf.~\cite{Aizenstat1962}.
\eit
Applying the usual techniques, we obtain (via Theorem \ref{t:pres}) a category presentation for $\PO$, which can then be rewritten (using Theorem \ref{t:pres2}) to yield the following:

\begin{thm}\label{t:PO}
The category $\PO$ has tensor presentation with generators
\[
V\mt\left(\begin{smallmatrix}1&2\\1&1\end{smallmatrix}\right)\in\PO_{2,1} \COMMA
U\mt\left(\begin{smallmatrix}1\\-\end{smallmatrix}\right)\in\PO_{1,0} \COMMA
\U\mt(\emptyset)\in\PO_{0,1}
\]
and relations
\[
\epfreseq \U\circ U = \io_0  \COMMa V\circ U = U\op U \COMMa (V\op I)\circ V = (I\op V)\circ V \COMMa (I\op \U)\circ V= I  = (\U\op I)\circ V .  
\]
\end{thm}

For the category $\O$, we first use the above-mentioned presentations for $\O_n$ to obtain (via Theorem~\ref{t:pres}) a presentation for $\O^+=\bigcup_{m,n\in\PP}\O_{m,n}$, and then rewrite it (using Theorem \ref{t:pres3}) to give the following:

\begin{thm}\label{t:O}
The category $\O$ has tensor presentation with generators
\[
V\mt\left(\begin{smallmatrix}1&2\\1&1\end{smallmatrix}\right)\in\O_{2,1} \AND
\U\mt(\emptyset)\in\O_{0,1}
\]
and relations
\[
\epfreseq (V\op I)\circ V = (I\op V)\circ V \AND (I\op \U)\circ V= I  = (\U\op I)\circ V .  
\]
\end{thm}

\begin{rem}
Theorems \ref{t:PO} and \ref{t:O} both include the relation $(I\op \U)\circ V= I  = (\U\op I)\circ V$, whereas \ref{PV2'} has only $(I\op \U)\circ V = I$.  The reason for this is that with the additional generator $X$, the second half of this relation follows from the first (and other parts of \ref{PV1'}--\ref{PV6'}); cf.~\eqref{e:PV3'}.
\end{rem}

For the category $\OI$, we could apply Theorems \ref{t:pres} and \ref{t:pres2} as usual, starting from presentations for the monoids $\OI_n$ as may be found in \cite{Fernandes2001} for example.  However, the relations from \cite{Fernandes2001} are substantially different from \ref{PV1}--\ref{PV13}, so this would require a great deal more work.  In any case, the category $\OI$ is simple enough that we can give a direct proof of the following:

\begin{thm}\label{t:OI}
The category $\OI$ has tensor presentation with generators
\[
U\mt\left(\begin{smallmatrix}1\\-\end{smallmatrix}\right)\in\OI_{1,0} \AND
\U\mt(\emptyset)\in\OI_{0,1}
\]
and the single relation $\U\circ U = \io_0$.
\end{thm}

\pf
Let $\De$ be the digraph over $\N$ with edges $U:1\to0$ and $\U:0\to1$, and let $\Phi:\De^\oast\to\OI$ be the morphism in the statement.  Let $\approx$ be the congruence on $\De^\oast$ generated by the relation $\U\circ U = \io_0$.

To show that $\Phi$ is surjective, let $f\in\OI_{m,n}$, and write $f=\left(\begin{smallmatrix} a_1&\cdots&a_k\\b_1&\cdots&b_k \end{smallmatrix}\right)$, where $a_1<\cdots<a_k$ (and $b_1<\cdots<b_k$).  For convenience, we also define $a_0=b_0=0$, and $a_{k+1}=m+1$ and $b_{k+1}=n+1$.  Then, defining $p_i=a_{i+1}-a_i-1$ and $q_i=b_{i+1}-b_i-1$ for all $0\leq i\leq k$, we have $f=w\Phi$, where
\begin{equation}\label{e:w}
w = (U^{\op p_0}\op\U^{\op q_0}) \op I \op (U^{\op p_1}\op\U^{\op q_1}) \op I \op \cdots \op (U^{\op p_{k-1}}\op\U^{\op q_{k-1}}) \op I \op (U^{\op p_k}\op\U^{\op q_k}).
\end{equation}

Next, ${\approx}\sub\ker(\Phi)$ is easily checked as usual.  For the reverse inclusion, it suffices to show that every term from $\De^\oast$ is $\approx$-equivalent to a term of the form \eqref{e:w}, as this form uniquely determines $w\Phi\in\OI$.  In fact, by \eqref{e:term2} and a simple induction, it is enough to show that for any term $w$ as in~\eqref{e:w} and for $z=\io_c\op x\op \io_d$, with $c,d\in\N$, $x\in\{U,\U\}$ and $\br(w)=\bd(z)$, $w\circ z$ is $\approx$-equivalent to a term of the desired form.  For convenience, we write $w=y_1\op\cdots\op y_l$, where $y_1,\ldots,y_l\in\{I,U,\U\}$ are the symbols in the order they appear in \eqref{e:w}.  For $0\leq i\leq l$, we write $w_i^-=y_1\op\cdots\op y_i$ and $w_i^+=y_{i+1}\op\cdots\op y_l$ (interpreting $w_0^-=w_l^+=\io_0$), so that $w=w_i^-\op w_i^+$ for all $i$.

\pfcase1  Suppose first that $x=U$.  Noting that $\br(w)=\bd(z)=c+d+1>c$, let $0\leq m\leq l$ be maximal so that $\br(w_m^-)=c$.  Since $\br(w_{m+1}^-)>c$, we must have $y_{m+1}=I$ or $\U$, and we have $w=w_m^-\op y_{m+1}\op w_{m+1}^+$ and $\br(w_{m+1}^+)=d$.  If $y_{m+1}=I$, then
\[
w\circ z = (w_m^- \op I \op w_{m+1}^+) \circ (\io_c \op U \op \io_d) = (w_m^-\circ\io_c) \op (I\circ U) \op (w_{m+1}^+\circ \io_d) = w_m^- \op U \op w_{m+1}^+.
\]
If $y_{m+1}=\U$, then we similarly show (using $\U\circ U\approx\io_0$) that $w\circ z \approx w_m^-\op w_{m+1}^+$.  In both cases, the final expressions are either already of the form \eqref{e:w}, or can be transformed into the desired form using the relation $U\op\U = \U\op U$; cf.~\eqref{e:UU}.  

\pfcase2  If $x=\U$, then we again let $0\leq m\leq l$ be maximal so that $\br(w_m^-)=c$, and this time we have $w\circ z = w_m^- \op \U \op w_m^+$, which is again either already of the form \eqref{e:w}, or can be transformed into such form using \eqref{e:UU}.  
\epf

%\begin{rem}
%The above proof can be adapted to treat other categories of mappings determined by mappings having combinations of injectivity, surjectivity, totality, isotonality.
%\end{rem}

\footnotesize
\def\bibspacing{-1.1pt}
\bibliography{biblio}
\bibliographystyle{abbrv}

\end{document}